\RequirePackage{fix-cm}
\documentclass[smallextended]{svjour3}       % onecolumn (second format)
\smartqed  % flush right qed marks, e.g. at end of proof
\usepackage{graphicx}
\usepackage{amssymb,amsmath}
\usepackage{subfig}
\usepackage{adjustbox}
\usepackage{enumitem}
\usepackage{alphalph}

\usepackage{mathptmx}      % use Times fonts if available on your TeX system
\begin{document}

\title{The periodic steady-state solution for queues with Erlang arrivals and service and time-varying periodic transition rates%\thanks{Grants or other notes
%about the article that should go on the front page should be
%placed here. General acknowledgments should be placed at the end of the article.}
}
%\subtitle{Do you have a subtitle?\\ If so, write it here}

\titlerunning{Erlang input-Erlang service}     % if too long for running head

\author{B.H. Margolius         %\and
        %Second Author %etc.
}

%\authorrunning{Short form of author list} % if too long for running head

\institute{B. H. Margolius \at
              Cleveland State University \\
              Tel.: +216-687-4682\\
              %Fax: +123-45-678910\\
              \email{b.margolius@csuohio.edu}     \\
              ORCID: 0000-0002-3069-8804      %  \\
}

\date{Received: date / Accepted: date}
% The correct dates will be entered by the editor

\maketitle

\begin{abstract}
We study a queueing system with Erlang arrivals with $k$ phases and Erlang service with $m$ phases.  Transition rates among phases vary periodically with time.  For these systems, we derive the asymptotic periodic distribution of the level and phase as a function of time within the period.  The asymptotic periodic distribution is analogous to a steady-state distribution for a system with constant rates.  If the time within the period is considered part of the state, then it is a steady-state distribution.  We also obtain waiting time and busy period distributions.  These solutions are expressed as infinite series.  We provide bounds for the error of the estimate obtained by truncating the series.  Examples are provided comparing the solution of the system of ordinary differential equation with a truncated state space to these asymptotic solutions involving remarkably few terms of the infinite series.
% \PACS{PACS code1 \and PACS code2 \and more}
\end{abstract}

\section{Introduction}
\label{intro}

In this paper, we explore several quantities related to the Erlang arrival, Erlang service queue with time-varying periodic transition rates.   The $E_k/E_m/1$ queue is a single server queue.  Arrivals occur in $k$ phases visited sequentially with transitions among phases occurring at rate $\lambda(t)$.  The transition rate is a periodic function of time.  Throughout this paper, we take the length of the period to be one. The service process is also Erlang.  It consists of $m$ phases.  Transitions among phases occur at rate $\mu(t)$ (also periodic), with each phase completed in sequence. 

When service is exponential, the standard deviation of the service time is equal to its expectation.  For Erlang-$m$ service (when rates are constant), with parameter $\mu$, the mean is $\frac{m}{\mu}$ and the variance is $\frac{m}{\mu^2}$.  Of course, similar facts hold for Erlang$-k$ arivals.  This is an advantage when modeling processes for which the variance and standard deviation of the service distribution are not equal.  Erlang-$m$ service or Erlang-$k$ arrivals also lets us track the stage of service or arrival, respectively, of the customer.   These are two advantages cited by Gayon, et al, in choosing Erlang service for modeling a single-item-make-to-stock production system in which items have Erlang production times \cite{Gayon_2009}.   Foh and Zukerman \cite{FohZukerman_2002} used Erlang service to model  random access protocols.    
%\cite{Sandmann_2010} \cite{Ikeda_1991}
Jayasuriya, et al \cite{Jayasuriya_2001} use generalized Erlang service to model channel holding times in a mobile environment.  Kuo and Wang \cite{KuoWang_1997} use an $M/E_m/1$ queue to model a machine repair problem.  Maritas and Xirokostas \cite{Maritas_1977} also study a machine repair problem using Erlang service.  Their model allows for more than one server.  Grassmann \cite{Grassmann_2011} provides additional examples of applications of $E_k/E_m/1$ queues.

Many researchers have studied the  $E_k/E_m/1$ queue, or the simpler $E_k/M/1$ or $M/E_m/1$ queues with constant transtion rates. A traditional approach using generating functions can be found, for example, in Saaty \cite{Saaty}, Kleinrock \cite{Kleinrock} and Medhi \cite{Medhi}.  This is the approach that we use in this paper, extending it to queues with time-varying periodic transition rates.  This paper extends related work applying this approach to other queues with time-varying periodic transition rates. See \cite{margolius_2019} when the generating functions for the queue-length process are scalar, and \cite{margolius_2021} for quasi-birth-death processes (QBD) when the generating functions for the queue-length process are vectors with a component for each phase.  In the 2021 paper, we use a two priority queue with finite waiting room for priority 2 customers as an extended example.

Queues with Erlang arrivals, Erlang service or both, have been analyzed by Smith \cite{smith_1953}, Syski \cite{Syski} and Tak\'{a}cs \cite{TAKACS} using Laplace transform techniques. Tak\'{a}cs studies the waiting time, queue length and busy period for a queue with Erlang arrivals and general service.   Truslove \cite{Truslove_1975} considers this queue with finite waiting room.  %Niu and Cooper \cite{Niu_Cooper_1993} apply transform-free techniques 
Leonenko \cite{Leonenko_2009} studies the transient solution to the $M/E_k/1$ queue following an approach due to Parthasarathy, \cite{Parthasarathy_1987}.  A paper by Griffiths, Leonenko and Williams \cite{GriffithsLeonenkoWilliams_2006} also provides an exact solution to the transient distribution of the $M/E_k/1$ queue.  Arizono, et al \cite{Arizono} use generating functions for the number of minimal lattice paths to find the equilibrium distribution for the $E_k/E_m/1$ queue length distribution. 

 $E_k/E_m/1$ queues may also be analyzed using matrix analytic methods \cite{Latouche}.    Grassmann \cite{Grassmann_2011}, in his 2011 paper, derives an effective method for finding the characteristic roots of the $k+m$ degree polynomial related to the waiting time distribution that arises from these methods.  He builds on the approach due to Syski \cite{Syski} and Smith \cite{smith_1953}.  Ivo Adan and Yiqiang Zhao \cite{AdanZhao_1996} study a $GI/E_m/1$ system and show that for arbitrarily distributed inter-arrival times and Erlang service, the waiting time distribution can be expressed as the finite sum of exponentials which depend on the roots of an equation. They also develop a method for finding these roots. Luh and Liu \cite{LuhLiu} study the $E_k/E_m/1$ queue and show that the roots of the characteristic polynomial associated with the process are simple if the arrival and service rates are real.  They use this result to construct a general solution space of vectors for the stationary solution of the queue length distribution.  Poyntz and Jackson \cite{PoyntzJackson} find the steady-state solution for the $E_k/E_m/r$ queue, illustrating the method with the $E_k/E_m/2$ queue due to the ``tediousness of the algebra''.

In this paper, we are studying the $E_k/E_m/1$ queue when the transition rates vary with time.  For two fairly recent surveys of research on queueing systems with time-varying parameters, the reader is referred to the papers by Schwarz, et al \cite{schwarz_selinka_stolletz_2016} and Whitt \cite{Whitt}. 

The remainder of this paper is divided into several sections.  In section \ref{sec:review} we provide a brief review of results for QBDs with time-varying periodic transition rates.  Section \ref{sec:EkEm1} provides the set up and analysis of the $E_k/E_m/1$ queue.  We explore the singularities of the generating function for the level distribution in section \ref{sec:singularities} to find an exact formula for the level and phase probabilities.  In section \ref{sec:error}, we provide error bounds for the level probabilities and show that for smooth functions, the truncated series for our exact formulas can be made arbitrarily close.  Section \ref{sec:waiting} gives formulas for the waiting time distribution and section \ref{sec:busy} derives the busy period distribution as the solution of a Volterra equation of the second kind.

\section{Review of results for quasi-birth-death (QBDs) processes with time-varying periodic rates}
\label{sec:review}
We study the asymptotic periodic solution of ergodic queues with Erlang arrivals and service and time-varying periodic transition rates.  The solutions are expressed in terms of an integral over a single period.  The integrals involve the idle probabilities for the system.  These idle probabilities may be estimated using Tikhonov regularization.  We do this within the framework developed in \cite{margolius_2019} and \cite{margolius_2021}.  The solutions are exact, but involve an infinite series.  The solutions may be estimated to arbitrary precision using finitely many terms.

We begin by recapping one of the main results from \cite{margolius_2021}.    The infinitesimal generator for a (QBD) with time-varying periodic transition rates:
\begin{equation*}
{\bf Q}(t)=\left[
\begin{array}{ccccc}
{\bf B}(t)&{\bf A_1}(t)&&&\\
{\bf A}_{-1}(t)&{\bf A}_0(t)&{\bf A}_1(t)&&\\
&{\bf A}_{-1}(t)&{\bf A}_0(t)&{\bf A}_1(t)&\\
&&\ddots&\ddots\\
\end{array}
\right].
\end{equation*}
This leads to the system of differential equations:
\begin{eqnarray}\label{eq:QBDODE}
\dot{\bf p}_0(t)&=&{\bf p}_0(t){\bf B}(t)+{\bf p}_1(t){\bf A}_{-1}(t)\nonumber\\
\dot{\bf p}_n(t)&=&{\bf p}_{n-1}(t){\bf A}_1(t)+{\bf p}_n(t){\bf A}_{0}(t)+{\bf p}_{n+1}(t){\bf A}_{-1}(t), \;n\in\mathbb{N}\backslash\{0\},
\end{eqnarray}
where \({\bf p}_n(t)\) is a \(K\)-element row vector whose \(j\)th component reflects the probability of being in phase \(j\) and level \(n\) at time \(t\).  The \({\bf A}_m(t)\), \(m=-1,0,1\) are \(K\times K\) matrices reflecting transitions among phases and within the current level or to an adjacent level.

We can use the system of ordinary differential equations given in (\ref{eq:QBDODE}) to solve for the generating function for the asymptotic periodic distribution (see Breuer \cite{Breuer} for more details). The asymptotic periodic distribution is the limiting distribution at time \(t\) within the period as the number of periods tends to infinity.  Such a limit will exist if the process is ergodic.  To obtain equation (\ref{eq:QBDkey}), we have assumed that \({\bf p}_n(t)={\bf p}_n(t-1)\), so \({\bf P}(z,t)=\sum_{n=0}^\infty {\bf p}_n(t)z^n={\bf P}(z,t-1)\).  Note that \({\bf P}(z,t)\) is a row vector of generating functions.  The coefficient of \(z^n\) of the \(j\)th component gives the asymptotic periodic probability of being in level \(n\) and phase \(j\) at time \(t\) within the period.

The {\it key} equation for the generating function is given by:
\begin{multline}\label{eq:QBDkey}
{\bf P}(z,t)=\sum_{j=0}^\infty {\bf p}_j(t)z^j=\\
\int_{t-1}^t{\bf p}_0(u)\left({\bf B}(u)-{\bf A}_0(u)-z^{-1}{\bf A}_{-1}(u)\right)\Phi(z,u,t)du\\
\times
\left({\bf I}-\Phi(z,t-1,t)\right)^{-1},
\end{multline}
where \(\Phi(z,u,t)\) is the generating function for the unbounded process, that is, the process that permits negative levels.  \(\Phi(z,u,t)\) is an evolution operator that satisfies
\begin{equation}\label{eq:Phi1}
\frac{d}{dt}\Phi(z,u,t)=\Phi(z,u,t){\bf A}_z(t),
\end{equation}
\begin{equation}\label{eq:Phi2}
\frac{d}{du}\Phi(z,u,t)=-{\bf A}_z(u)\Phi(z,u,t),
\end{equation}
and
\begin{equation}\label{eq:Phi3}
\Phi(z,t,t)={\bf I},
\end{equation}
where
\begin{equation}\label{eq:Az}
{\bf A}_z(t)={\bf A}_{1}(t)z+{\bf A}_{0}(t)+{\bf A}_{-1}(t)z^{-1}.
\end{equation}
For further details, see \cite{margolius_2021}.

\section{Erlang arrivals and service, the $E_k/E_m/1$ queue}
\label{sec:EkEm1}

The queue with Erlang arrivals and service, the  $E_k/E_m/1$ queue, can be modeled with a three-dimensional state space $\{X(t),K(t),J(t)\}$ in which \(X(t)\) represents the level at time \(t\), $K(t)$ is the arrival phase, and $J(t)$ is the service phase.  Arrivals are \(k\)-Erlang with time-varying periodic transition rate  $\lambda(t)$ among arrival phases.  Service is \(m\)-Erlang with time-varying periodic transition rate $\mu(t)$ among service phases.  This process can be modeled as a QBD.  Define the following transition rate matrices:
\begin{equation*}
D_0(t)=
\left[
\begin{array}{cccc}
-\lambda(t)&\lambda(t)&&\\
&-\lambda(t)&\lambda(t)&\\
&\ddots&\ddots&\\
&&&-\lambda(t)
\end{array}
\right],
\end{equation*}  
\begin{equation*}
D_1(t)=
\left[
\begin{array}{cccc}
&&&\\
&&&\\
\lambda(t)&&&
\end{array}
\right],
\end{equation*}  
\begin{equation*}
C_0(t)=
\left[
\begin{array}{cccc}
-\mu(t)&\mu(t)&&\\
&-\mu(t)&\mu(t)&\\
&\ddots&\ddots&\\
&&&-\mu(t)
\end{array}
\right],
\end{equation*}  
\begin{equation*}
C_1(t)=
\left[
\begin{array}{cccc}
&&&\\
&&&\\
\mu(t)&&&
\end{array}
\right].
\end{equation*} 
The matrices $D_0(t)$ and $D_1(t)$ are $k\times k$, and the matrices $C_0(t)$ and $C_1(t)$ are $m\times m$ reflecting transitions among arrival and service phases, respectively.  Let \({\bf e}_1\) represent an appropriately dimensioned row vector with a one in the first position, and zeros elsewhere.  The inter-arrival arrival distribution in the constant rate case is given by 
\begin{equation*}
F_{T_a}(t)=1-{\bf e}_1{\rm exp}(D_1 t){\bf 1}_k=1-{\rm e}^{-\lambda t}\sum_{j=0}^{k-1}\frac{\lambda^j t^j}{j!}
\end{equation*}
where ${\bf 1}_k$ is a $k\time 1$ column vector of ones.  
In the time-varying case, we have
\begin{equation*}
F_{T_{a,u}}(t)=1-{\bf e}_1\Lambda(u,u+t){\bf 1}_k=1-{\rm e}^{-\int_u^{t+u}\lambda(\nu)d\nu}\sum_{j=0}^{k-1}\frac{\left({\rm e}^{\int_u^{t+u}\lambda(\nu)d\nu}\right)^j }{j!}
\end{equation*}
where $\Lambda(u,t)$ is an evolution operator satisfying
\begin{equation*}
\frac{d}{dt}\Lambda(u,t)=\Lambda(u,t)D_1(t),
\end{equation*}
\begin{equation*}
\frac{d}{du}\Lambda(u,t)=-D_1(u)\Lambda(u,t),
\end{equation*}
and
\begin{equation*}
\Lambda(t,t)={\bf I}.
\end{equation*}
An explicit formula for $\Lambda(u,u+t)$ is
\begin{equation*}
\Lambda(u,u+t)={\rm e}^{-\int_u^{t+u}\lambda(\nu)d\nu}\left[\begin{array}{cccc}
1&\int_u^{t+u}\lambda(\nu)d\nu&\cdots&\frac{\left(\int_u^{t+u}\lambda(\nu)d\nu\right)^{k-1}}{(k-1)!}\\
0&\ddots&&\frac{\left(\int_u^{t+u}\lambda(\nu)d\nu\right)^{k-2}}{(k-2)!}\\
\vdots&&\ddots&\vdots\\
0&\cdots&0&1\\
\end{array}\right].
\end{equation*}
The departure process is similarly defined with 
\begin{equation*}
F_{T_d}(t)=1-{\bf e}_1{\rm exp}(C_1 t){\bf 1}_k=1-{\rm e}^{-\mu t}\sum_{j=0}^{k-1}\frac{\mu^j t^j}{j!}
\end{equation*}
in the constant rate case and 
\begin{equation*}
F_{T_{d,u}}(t)=1-{\bf e}_1 M(u,u+t){\bf 1}_k=1-{\rm e}^{-\int_u^{t+u}\mu(\nu)d\nu}\sum_{j=0}^{m-1}\frac{\left({\rm e}^{\int_u^{t+u}\mu(\nu)d\nu}\right)^j }{j!}
\end{equation*}
when rates are time-varying.  The matrix function \(M(u,u+t)\) is given by
\begin{equation*}
M(u,u+t)={\rm e}^{-\int_u^{t+u}\mu(\nu)d\nu}\left[\begin{array}{cccc}
1&\int_u^{t+u}\mu(\nu)d\nu&\cdots&\frac{\left(\int_u^{t+u}\mu(\nu)d\nu\right)^{m-1}}{(m-1)!}\\
0&\ddots&&\frac{\left(\int_u^{t+u}\mu(\nu)d\nu\right)^{m-2}}{(m-2)!}\\
\vdots&&\ddots&\vdots\\
0&\cdots&0&1\\
\end{array}\right].
\end{equation*}

The infinitesimal generator for the $E_k/E_m/1$ queue, ${\bf Q}(t)$, is given below.  We arrange states in lexicographic order, i.e. $(0,1)$, $(0,2)$, $\dots$, $(0,k)$, $(1,1,1)$, $(1,1,2)$, $\dots$, $(1,1,m)$, $(1,2,1)$, $\dots$, $(1,2,m)$, $\dots$.  Then the infinitesimal generator for this process is given by
\begin{equation}\label{eq:Q}
{\bf Q}(t)=
\left[
\begin{array}{cccccc}
D_0&Q_{0,1}&&&&\\
Q_{1,0}&D_0\oplus C_0&D_1\otimes {\bf I}_m&&&\\
&{\bf I}_k\otimes C_1&D_0\oplus C_0&D_1\otimes {\bf I}_m&&\\
&&{\bf I}_k\otimes C_1&D_0\oplus C_0&D_1\otimes {\bf I}_m&\\
&&&\ddots&\ddots&\ddots\\
&&&&\ddots&\ddots
\end{array}
\right],
\end{equation}
where 
\[Q_{0,1}(t)=\left[
\begin{array}{cccc}
&&&\\
&&&\\
\lambda(t)&&&
\end{array}
\right]\]
is a \(k\times mk\) matrix, 
\[Q_{1,0}(t)={\bf I}_k\otimes\left[\begin{array}{c}0\\\vdots\\0\\\mu(t)\end{array}\right],\]
${\bf I}_k$ is a $k\times k$ identity matrix, $\otimes$ and $\oplus$ represent the Kronecker product and Kronecker sum, respectively.  For definitions of the Kronecker product and Kronecker sum see, for example,  the textbook by Alan Laub \cite{Laub}, or MathWorld \cite{KroneckerProduct}, \cite{KroneckerSum}.  
The dependence of ${\bf Q}(t)$ on $t$ is suppressed in the notation in equation (\ref{eq:Q}).

For the $E_k/E_m/1$ QBD, 
\begin{equation*}
{\bf A}_{-1}(t)={\bf I}_k\otimes C_1(t),
\end{equation*}
\begin{equation*}
{\bf A}_0(t)=D_0(t)\otimes{\bf I}_m+{\bf I}_k\otimes C_0(t)
\end{equation*} 
and
\begin{equation*}
{\bf A}_1(t)=D_1(t)\otimes {\bf I}_m.
\end{equation*}
With ${\bf A}(z,t)=z^{-1}{\bf A}_{-1}(t)+{\bf A}_{0}(t)+z{\bf A}_{1}(t)$ as in equation (\ref{eq:Az}), then the function $\Phi(z,u,t)$ is an evolution operator which satisfies equations (\ref{eq:Phi1}), (\ref{eq:Phi2}) and (\ref{eq:Phi3}).  For each of these $km\times km$ matrices, we reference the components of the matrix as $((a_1,s_1)(a_2,s_2))$ where $(a_1,s_1)$ refer to the arrival and service phases of the row and $(a_2,s_2)$ give the arrival and service phases of the column.

The function $\Phi(z,u,t)$ is a Laurent series in $z$ with $km\times km$ matrix coefficients where the $((a_1,s_1),(a_2,s_2))$ entry of the coefficient on $z^\ell$ represents the probability of a net change of $\ell$ levels during the time interval $[u,t)$ and a sequence of transitions that begin in arrival phase $a_1$ and service phase $s_1$ at time $u$ and end in arrival phase $a_2$ and service phase $s_2$ at time $t$.    

The {\it key} equation gives the generating function for this QBD in terms of an integral over a single time period.  See \cite{margolius_2021} for the general case for QBDs.  For the $E_k/E_m/1$ system, the {\it key} equation is given by
\begin{multline}\label{eq:EkEmkey}
{\bf P}(z,t)=\sum_{j=1}^\infty {\bf p}_j(t)z^j=\\
\int_{t-1}^t\left({\bf p}_0(u)z{\bf Q}_{0,1}(u)-{\bf p}_1(u){\bf A}_{-1}(u)\right)\Phi(z,u,t)du\\
\times
\left({\bf I}-\Phi(z,t-1,t)\right)^{-1}.
\end{multline}

We can write an explicit formula for the evolution operator $\Phi(z,u,t)$.  Note that ${\bf A}(z,t)$ may be expressed in terms of a Kronecker sum as
\begin{equation*}
{\bf A}(z,t)=\left(D_0(t)+zD_1(t)\right)\oplus\left(C_0(t)+z^{-1}C_1(t)\right).
\end{equation*}
%where $\oplus$ is the Kronecker sum defined by \[A_{k\times k}\oplus B_{m\times m}=A\otimes {\bf I}_m +{\bf I}_k\otimes B\] and $\otimes$ is the Kronecker product.  
The eigenvalues of ${\bf A}(z,t)$ are the sum of eigenvalues of the matrices $\left(D_0(t)+zD_1(t)\right)$ and $\left(C_0(t)+z^{-1}C_1(t)\right)$ and the eigenvectors are the Kronecker product of the corresponding eigenvectors.
 
Define
\begin{equation*}
\omega_{K}={\rm e}^{\frac{2\pi i}{K}}=\cos\!\left(\frac{2\pi}{K}\right)+i\sin\!\left(\frac{2\pi}{K}\right),
\end{equation*}
a \(K\)th primitive root of unity.

The matrix $D_0(t)+zD_1(t)$ has eigenvalues 
\begin{equation}\label{eq:aeigenvalue}
\xi_\ell(z,t)=\lambda(t)(\omega_k^\ell z^{1/k}-1),\;\;\ell=0,\dots,k-1,
\end{equation}
and corresponding eigenvectors
\begin{equation}\label{eq:aeigenvectors}
v_\ell=\frac{1}{\sqrt{k}}\left[
\begin{array}{c}
z^{(1-k)/k}\omega_k^0\\
z^{(2-k)/k}\omega_k^\ell\\
z^{(3-k)/k}\omega_k^{2\ell}\\
\vdots\\
z^{0}\omega_k^{(k-1)\ell}
\end{array}
\right].
\end{equation}
The matrix $C_0(t)+\frac{1}{z}C_1(t)$ has eigenvalues
\begin{equation}\label{eq:deigenvalue}
\epsilon_j(z,t)=\mu(t)(\omega_m^{j} z^{-1/m}-1),\;\;j=0,\dots,m-1,
\end{equation}
and corresponding eigenvectors
\begin{equation}\label{eq:deigenvectors}
u_j=\frac{1}{\sqrt{m}}\left[
\begin{array}{c}
z^{(m-1)/m}\omega_m^0\\
z^{(m-2)/m}\omega_m^{j}\\
z^{(m-3)/m}\omega_m^{2j}\\
\vdots\\
z^{0}\omega_m^{(m-1)j}
\end{array}
\right].
\end{equation}
We can now compute the eigenvalues and eigenvectors of ${\bf A}(z,t)=(D_0(t)+zD_1(t))\oplus(C_0(t)+\frac{1}{z}C_1(t))$.  The eigenvalues are
\begin{equation*}
\xi_\ell(z,t)+\epsilon_j(z,t),\;\;j=0,\dots,m-1,\; \ell=0,\dots,k-1.
\end{equation*}
The eigenvector for $(D_0(t)+zD_1(t))\oplus(C_0(t)+\frac{1}{z}C_1(t))$ corresponding to the eigenvalue $\xi_\ell(z,t)+\epsilon_j(z,t)$ is
\begin{equation}\label{eq:eigenvectors}
v_\ell\otimes u_j \;\;j=0,\dots,m-1,\; \ell=0,\dots,k-1.
\end{equation}
Note that while the eigenvalues depend on \(t\), the eigenvectors do not.

This enables us to easily compute the eigenvalues for the matrices:  $\Phi(z,u,t)$, $\left({\bf I}-\Phi(z,t-1,t)\right)^{-1}$ and $\Phi(z,u,t)\left({\bf I}-\Phi(z,t-1,t)\right)^{-1}$.  The eigenvectors for each of these matrices are those given in equation (\ref{eq:eigenvectors}), the same as for ${\bf A}(z,t)$.  Let 
\begin{equation*}
\bar\xi_\ell(z)=\int_0^1\lambda(u)(\omega_k^\ell z^{1/k}-1)du=\bar\lambda(\omega_k^\ell z^{1/k}-1)
\end{equation*}
and 
\begin{equation*}
\bar\epsilon_j(z)=\int_0^1\mu(u)(\omega_m^j z^{-1/m}-1)du=\bar\mu(\omega_m^j z^{-1/m}-1)
\end{equation*}
give the average value of the eigenvalues for the arrival and departure processes, respectively, over a single time-period.  We have defined $\bar\lambda=\int_{t-1}^t\lambda(u)du$,  the average value of $\lambda(t)$ over a single time period, and $\bar\mu=\int_{t-1}^t\mu(u)du$, the average value of $\mu(t)$ over a single time period.  Then the eigenvalues for the four matrices with common eigenvectors are as given in table \ref{tab:eigen}.

\begin{table}
\caption{Four matrix functions which share the eigenvectors \(v_\ell\otimes u_j\).}
\label{tab:eigen}   
\begin{center}   
\begin{tabular}{|c|c|}
\hline
Matrix&Eigenvalue\\
\hline
\hline
${\bf A}(z,t)$&$\xi_\ell(z,t)+\epsilon_j(z,t)$\\
\hline
$\Phi(z,u,t)$&${\rm exp}\{\int_u^t(\xi_\ell(z,\nu)+\epsilon_j(z,\nu))d\nu\}$\\
\hline
$({\bf I}-\Phi(z,t-1,t))^{-1}$&$(1-{\rm exp}\{(\bar\xi_\ell(z)+\bar\epsilon_j(z)\})^{-1}$\\
\hline
$\Phi(z,u,t)({\bf I}-\Phi(z,t-1,t))^{-1}$&${\rm exp}\{\int_u^t(\xi_\ell(z,\nu)+\epsilon_j(z,\nu))d\nu\}$\\
&$\times (1-{\rm exp}\{(\bar\xi_\ell(z)+\bar\epsilon_j(z)\})^{-1}$\\
\hline
\end{tabular}
\end{center}
\end{table}

Define the matrices 
\[\Omega_K = \left[
\begin{array}{cccc}
\omega_K^0&\omega_K^0&\cdots&\omega_K^0\\[.2 cm]
\omega_K^0&\omega_K^1&\cdots&\omega_K^{K-1}\\[.2 cm]
\vdots&\ddots&\ddots&\vdots\\[.2 cm]
\omega_K^0&\omega_K^{K-1}&\cdots&\omega_K^{(K-1)^2}
\end{array}
\right],\] 
and \(\overline\Omega_K\), its complex conjugate.  
Let 
\[H_C D_C H^{-1}_C=C_0+z^{-1}C_1\]
where $D_C$ is a diagonal matrix of the eigenvalues, $\epsilon_j(z,t)$, of $C_0+z^{-1}C_1$ and $H_C$ is a matrix whose columns are the eigenvectors, $u_j$, of $C_0+z^{-1}C_1$.
The matrix 
\[H_C=\frac{1}{\sqrt{m}}{\rm diag}[\begin{array}{cccc}z^{(m-1)/m}& z^{(m-2)/m}&\cdots& 1\end{array}]\Omega_m.\]
Similarly, let 
\[H_D D_D H^{-1}_D=D_0+zD_1\]
where $D_D$ is a diagonal matrix of the eigenvalues, $\xi_\ell(z,t)$, of $D_0+zD_1$ and $H_D$ is a matrix whose columns are the eigenvectors, $v_\ell$, of $D_0+zD_1$.  The matrix
\[H_D=\frac{1}{\sqrt{k}}{\rm diag}[\begin{array}{cccc}z^{(1-k)/k}& z^{(2-k)/k}&\cdots& 1\end{array}]\Omega_k.\]
Then 
\begin{multline*}
H_D\otimes H_C=\\
\frac{1}{\sqrt{km}}\left({\rm diag}[\begin{array}{cccc}z^{(1-k)/k}& z^{(2-k)/k}&\cdots& 1\end{array}]\otimes{\rm diag}[\begin{array}{cccc}z^{(m-1)/m}& z^{(m-2)/m}&\cdots& 1\end{array}]\right)\\
\times\left(\Omega_k\otimes \Omega_m\right).
\end{multline*}
Similarly,
\begin{multline*}
H_D^{-1}\otimes H_C^{-1}=
\frac{1}{\sqrt{km}}\left(\overline\Omega_k\otimes\overline\Omega_m\right)\\
\left({\rm diag}[\begin{array}{cccc}z^{(k-1)/k}& z^{(k-2)/k}&\cdots& 1\end{array}]\otimes{\rm diag}[\begin{array}{cccc}z^{(1-m)/m}& z^{(2-m)/m}&\cdots& 1\end{array}]\right).
\end{multline*}

Then 
\[\Xi_{(a_1,s_1)(a_2,s_2)}=\frac{1}{km}z^{\frac{s_2-s_1}{m}+\frac{a_1-a_2}{k}}\sum_{i=0}^{m-1}\omega_m^{i(s_1-s_2)}\sum_{\ell=0}^{k-1}d_{\ell,i}\omega_k^{\ell(a_1-a_2)}\]
where $\Xi$ is one of the matrix functions \({\bf A}(z,t)\), \(\Phi(z,u,t)\), \(\left({\bf I}-\Phi(z,t-1,t)\right)^{-1}\) or \linebreak \(\Phi(z,u,t)\left({\bf I}-\Phi(z,t-1,t)\right)^{-1}\), given in table \ref{tab:eigen} and the $d_{\ell,i}$ are the corresponding eigenvalues.

\paragraph{Generating functions and roots of unity.}
Recall that 
\begin{equation}\label{eq:pluck}
F^{(K,j)}(x)=\frac{1}{K}x^{-\frac{j}{K}}\sum_{\ell=0}^{K-1}F(\omega_K^\ell x^{\frac{1}{K}})\omega_K^{-j\ell}
\end{equation}
where $F(x)=\sum_{n=0}^\infty a_nx^n$ and $F^{(K,j)}(x)=\sum_{n=0}^\infty a_{Kn+j}x^n$.
See Herbert Wilf's text {\it generatingfunctionology}\cite{wilf} for more details on generating functions and the role of roots of unity.

We apply this formula twice to 
\begin{equation*}%\label{eq:main}
\frac{1}{km}z^{\frac{s_2-s_1}{m}+\frac{a_1-a_2}{k}}\sum_{i=0}^{m-1}\omega_m^{i(s_1-s_2)}\sum_{\ell=0}^{k-1}d_{\ell,i}\omega_k^{\ell(a_1-a_2)}
\end{equation*}
to obtain explicit formulas for  \([z^n]\Xi_{(a_1,s_1)(a_2,s_2)}\), the components of the coefficient matrices.

The generating function for a Poisson random variable appears several times in these expressions.  
Recall that the generating function for a Poisson random variable is
\begin{equation*}
{\rm e}^{\lambda(z-1)}={\rm e}^{-\lambda}\sum_{n=0}^\infty \frac{\lambda^n z^n}{n!}
\end{equation*}
and for steps to the left, the generating function is
\begin{equation*}
{\rm e}^{\mu(z^{-1}-1)}={\rm e}^{-\mu}\sum_{n=0}^\infty \frac{\mu^n z^{-n}}{n!}
\end{equation*}
with the product of these forming the generating function for a random walk:
\begin{equation*}
{\rm e}^{\lambda(z-1)+\mu(z^{-1}-1)}={\rm e}^{-\lambda-\mu}\sum_{n=-\infty}^\infty z^n
\sum_{\ell=0\vee(-n)}^\infty\frac{\mu^\ell \lambda^{n+\ell} }{\ell!(n+\ell)!}.
\end{equation*}
Note that these Poisson generating functions appear in three of our $\Xi$ matrices.  We work out in detail, the simplest of these. When $\Xi=\Phi(z,u,t)$ with $s=s_2-s_1$, $a=a_2-a_1$, $\lambda=\int_u^t\lambda(\nu)d\nu$ and $\mu=\int_u^t\mu(\nu)d\nu$, we have
\begin{align}\label{eq:Phicoeff}
&\left[\Phi(z,u,t)\right]_{(a_1,s_1),(a_2,s_2)}\\
&=
\frac{1}{km}z^{\frac{s}{m}-\frac{a}{k}}\sum_{i=0}^{m-1}\omega_m^{-is}\sum_{\ell=0}^{k-1}{\rm e}^{\lambda(\omega_k^\ell z^{1/k}-1)+\mu(\omega_m^i z^{-1/m}-1)}\omega_k^{-\ell a}\nonumber\\
&=\frac{z^{\frac{s}{m}}}{m}\sum_{i=0}^{m-1}{\rm e}^{\mu(\omega_m^i z^{-1/m}-1)}\omega_m^{-is}\frac{z^{-\frac{a}{k}}}{k}\sum_{\ell=0}^{k-1}{\rm e}^{\lambda(\omega_k^\ell z^{1/k}-1)}\omega_k^{-\ell a}&&\text{(factoring exponentials)} \nonumber\\
&=\frac{z^{\frac{s}{m}}}{m}\sum_{i=0}^{m-1}{\rm e}^{\mu(\omega_m^i z^{-1/m}-1)}\omega_m^{-is}
{\rm e}^{-\lambda}\sum_{n=0}^\infty\frac{\lambda^{nk+a}z^n}{(nk+a)!}&&\text{(applying equation (\ref{eq:pluck}))}\nonumber \\
&={\rm e}^{-\lambda-\mu}\sum_{\ell=0}^\infty \frac{\mu^{m\ell+s}z^{-\ell}}{(m\ell+s)!}\sum_{n=0}^\infty\frac{\lambda^{nk+a}z^n}{(nk+a)!}&&\text{(applying equation (\ref{eq:pluck}))}\nonumber \\
&={\rm e}^{-\lambda-\mu}\sum_{n=-\infty}^\infty z^n\sum_{\ell=0\vee(-n)}^\infty\frac{\mu^{\ell m+s}\lambda^{(n+\ell)k+a}}{(\ell m+s)!((n+\ell)k+a)!}&&\text{(coefficient on } z^n).\nonumber
\end{align}
The coefficient on $z^n$ in this Laurent series reflects the probability of $n$ more arrivals (which require completion of $k$ phases at rate $\lambda$) than service completions (which require completion of $m$ service phases at rate $\mu$) and a transition from arrival phase $a_1$ to arrival phase $a_2$ (net change $a=a_2-a_1$) and from service phase $s_1$ to $s_2$ (net change $s=s_2-s_1$) occurring during the time interval from $u$ to $t$.

The $((a_1,s_1),(a_2,s_2))$ component of the matrix coefficient on $z^n$ of the function  \(\Phi(z,u,t)\left({\bf I}-\Phi(z,t-1,t)\right)^{-1}\) gives the expected number of times $t$ within the period that the process has made a net change of $n$ levels and is in arrival phase $a_2$ and service phase $s_2$ having started at phases $(a_1,s_1)$ at time $u$ within an earlier period.  This coefficient does not count the expected number of visits, but rather the expected number of periods that the process is in a given state at time $t$ within the period.

The components of the matrix function \(\Phi(z,u,t)\left({\bf I}-\Phi(z,t-1,t)\right)^{-1}\) are linear combinations of generating functions of the form
\begin{multline*}
{\rm e}^{\int_u^t\left(\lambda(\nu)(z-1)+\mu(\nu)(z^{-1}-1)\right)d\nu}\left(1-{\rm e}^{\bar\lambda(z-1)+\bar\mu(z^{-1}-1)}\right)^{-1}\\
=\sum_{n=-\infty}^\infty z^n\sum_{\ell=0\vee(-n)}^\infty\sum_{j=0}^\infty{\rm e}^{-\int_u^{t+j}\left(\lambda(\nu)+\mu(\nu)\right)d\nu}\frac{\left(\int_u^{t+j}\lambda(\nu)d\nu\right)^{\ell+n}\left(\int_u^{t+j}\mu(\nu)d\nu\right)^{\ell}}{(\ell+n)!\ell!}
\end{multline*}
evaluated at $k$th and $m$th roots of the indeterminate $z$ times a root of unity.  An exact formula for the coefficient on $z^n$ of the $((a_1,s_1)(a_2,s_2))$ component is given by
\begin{multline*}
[z^n]\left[\Phi(z,u,t)\left({\bf I}-\Phi(z,t-1,t)\right)^{-1}\right]_{(a_1,s_1)(a_2,s_2)}=\\
\sum_{\ell=0\vee(-n)}^\infty\sum_{j=0}^\infty{\rm e}^{-\int_u^{t+j}\left(\lambda(\nu)+\mu(\nu)\right)d\nu}\frac{\left(\int_u^{t+j}\lambda(\nu)d\nu\right)^{(\ell+n)k+a}\left(\int_u^{t+j}\mu(\nu)d\nu\right)^{\ell m+s}}{((\ell+n)k+a)!(\ell m+s)!}.
\end{multline*}
For the $((a_1,s_1)(a_2,s_2))$ component of $\left({\bf I}-\Phi(z,t-1,t)\right)^{-1}$, we have
\begin{multline*}
[z^n]\left[\left({\bf I}-\Phi(z,t-1,t)\right)^{-1}\right]_{(a_1,s_1)(a_2,s_2)}=\\
\sum_{\ell=0\vee(-n)}^\infty\sum_{j=0}^\infty{\rm e}^{-(\bar\lambda+\bar\mu)j}\frac{(j\bar\lambda)^{(\ell+n)k+a}(j\bar\mu)^{\ell m+s}}{((\ell+n)k+a)!(\ell m+s)!}.
\end{multline*}
These formulas, while exact, are not conducive to computation.

\section{Singularity analysis}
\label{sec:singularities}
Following the approach of Sedgewick and Flajolet \cite{FLAJOLET}, we note that the singularities of the generating function are reflected in the coefficients.  In this section, we explore the zeros of the denominator of the generating function, \(P(z,t)\).
Note that the generating function $P(z,t)$ has singularities wherever 
\[1-{\rm exp}\left\{\bar\lambda\omega_k^\ell z^{1/k}-1)+\bar\mu(\omega_m^j z^{-1/m}-1)\right\}=0.\]  This occurs for $z$ such that  
\begin{equation}\label{eq:mypoly1}
\bar\lambda(\omega_k^\ell z^{1/k}-1)+\bar\mu(\omega_m^j z^{-1/m}-1)=2\pi i n,\;\;n\in\mathbb{Z}.
\end{equation}
Let $y=z^{\frac{1}{km}}$ when $\ell=j=0$, then equation (\ref{eq:mypoly1})  becomes
\begin{equation}\label{eq:mypoly2}
\bar\lambda y^{m+k}-(\bar\lambda+\bar\mu+2\pi i n)y^k+\bar\mu=0,\;\;n\in\mathbb{Z}.
\end{equation}

\begin{figure}[!tbp]
\begin{tabular}{cc}
\subfloat[Zeros of equation (\ref{eq:theroots}).  The $k=7$ petaled rose is inside the unit circle.]{\includegraphics[width=0.46\textwidth]{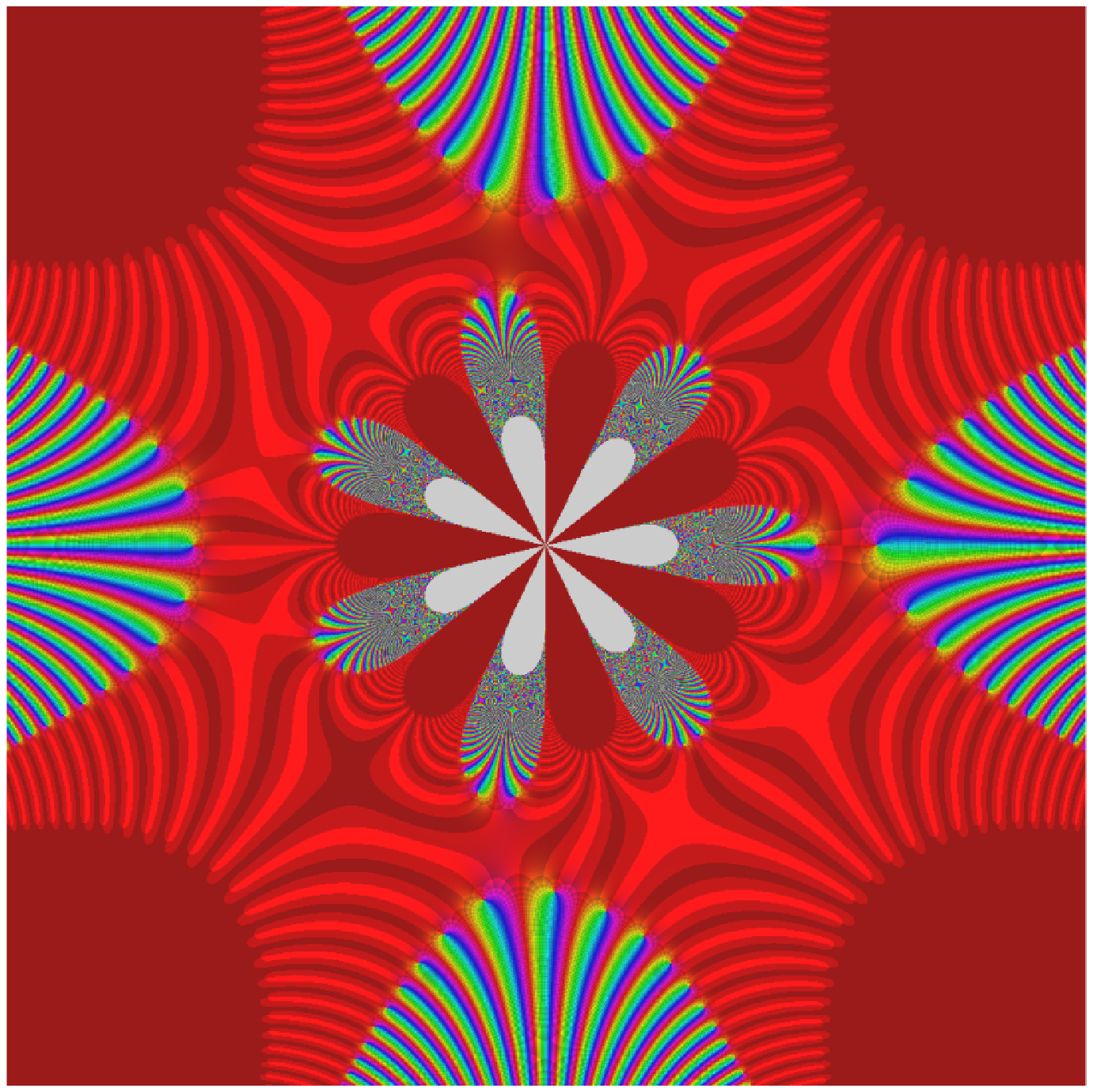}} &
\subfloat[Zeros of the real part of equation (\ref{eq:theroots}) are shown in white; zeros of the imaginary part for $n=3$ are shown in cyan.  The $m=4$ intersections of the cyan and white curves outside of the unit circle show the $m$ roots for $n=3$ outside the unit circle.  The intersection of the white and cyan seven petal roses show the $k=7$ roots corresponding to $n=3$ inside the unit circle.]{\includegraphics[width=0.46\textwidth]{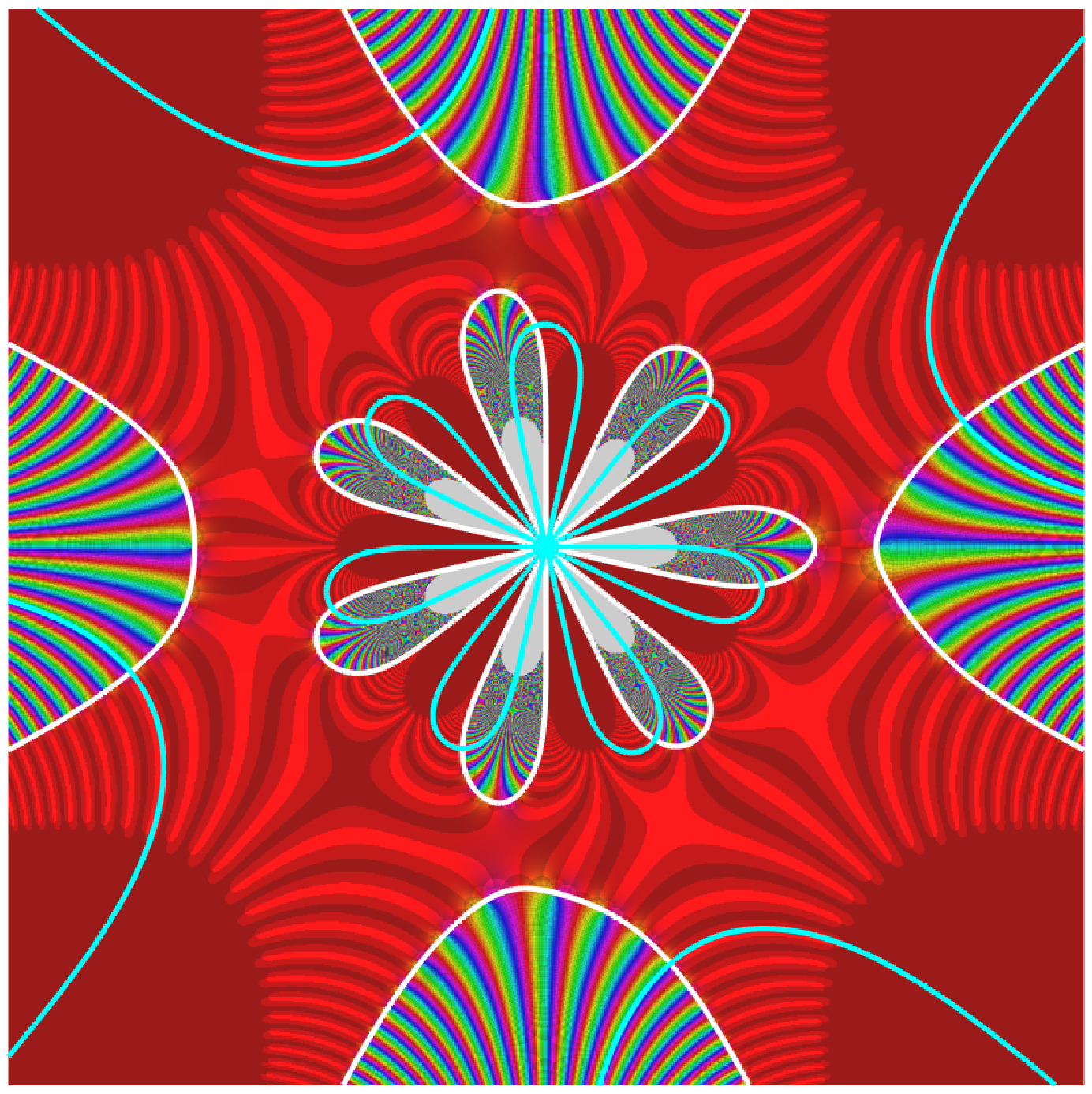}}
\end{tabular}
\caption{These figures were produced using Matlab code written by Elias Wegert \cite{Wegert}. Both plots show the region from $-2-2i$ to $2+2i$. The expression given in equation (\ref{eq:theroots}) is plotted in the complex plane.  Shading shows contour lines.  The colors represent the argument, so points where multiple colors come together are zeros of the function.  In this example, $k=7$, $m=4$, $\bar\lambda=3$ and $\bar\mu=5$.  Note that there are $k=7$ petals in the rose inside the unit circle and $m=4$ inverted petals outside the unit circle.}
\label{fig:zeros}
\end{figure}

Figure \ref{fig:zeros} shows the roots of 
\begin{equation}\label{eq:theroots}
1-{\rm e}^{\bar\lambda(y^m-1)+\bar\mu(y^{-k}-1)}.
\end{equation}
\begin{figure}
\begin{tabular}{cc}
\subfloat[Zeros of the real part of equation (\ref{eq:theroots}) are asymptotic to solutions to $r^k=\frac{\bar\mu}{\bar\lambda+\bar\mu}\cos(k\theta)$ (shown in blue inside the unit circle) as $r\to 0$ and asymptotic to $r^m=\frac{\bar\lambda+\bar\mu}{\bar\lambda}\sec(m\theta)$ (shown in red outside the unit circle) as $r\to\infty$.  Zeros of the real part of equation (\ref{eq:theroots}) are shown as dashed black lines.]
{\includegraphics[width=0.42\textwidth]{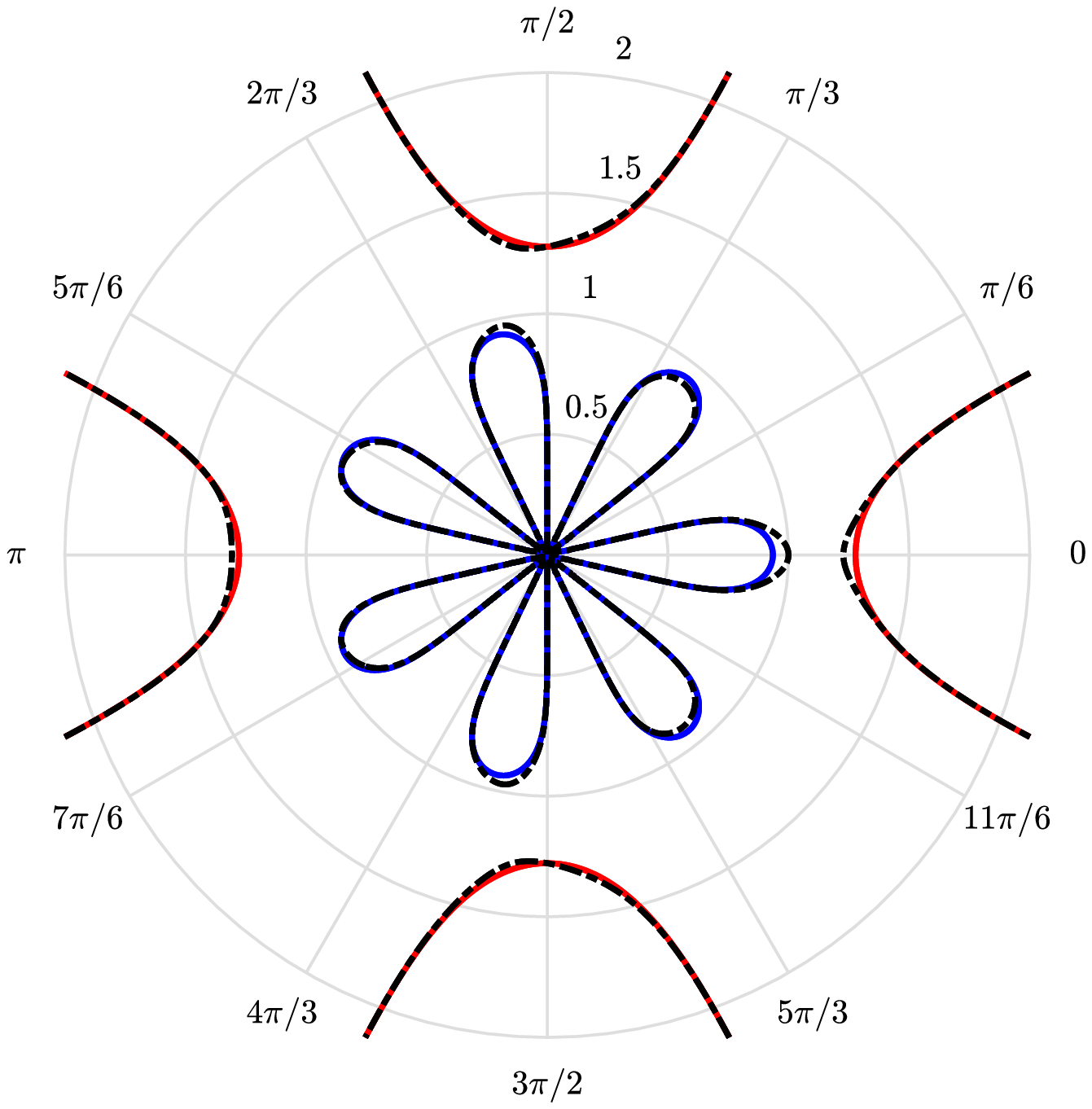}}&
\subfloat[Zeros of the real part of equation (\ref{eq:theroots}) are shown in white for an $E_3/E_5/1$ system; zeros of the imaginary part for $n=-7$ are shown in cyan.  The $m=5$ intersections of the cyan and white curves outside of the unit circle show the $m$ roots for $n=-7$ outside the unit circle.  The intersection of the white and cyan three petal roses show the $k=3$ roots corresponding to $n=-7$ inside the unit circle.]
{\includegraphics[width=0.48\textwidth]{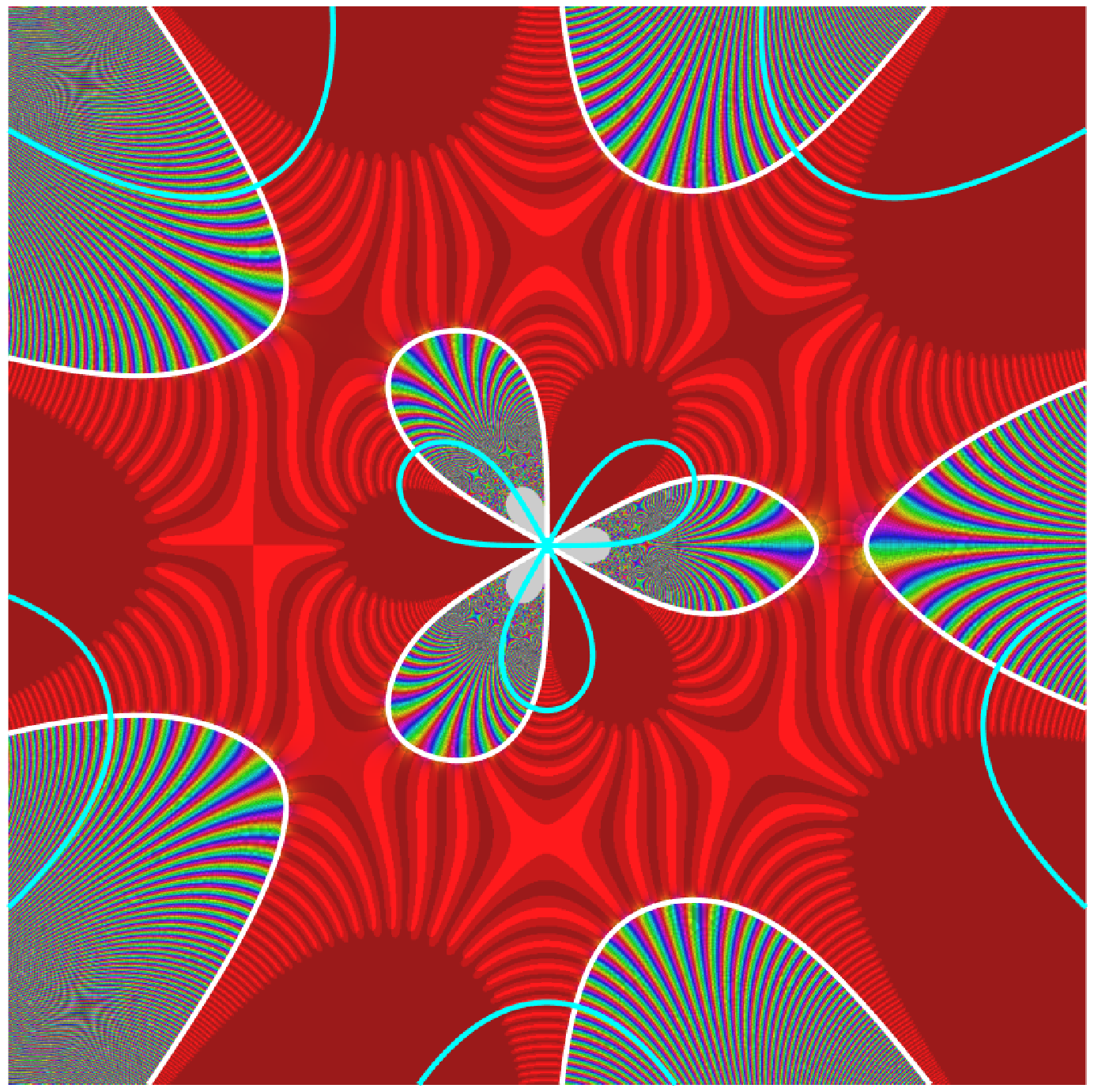}}
\end{tabular}
\caption{Zeros of the real part of the denominator of the generating function for an $E_7/E_4/1$ queueing system are shown on the left, and zeros of the denominator of the generating function for an $E_3/E_5/1$ queueing system are shown on the right.  Note that the petals inside the unit circle correspond to the number of arrival phases.}
\end{figure}

Using Rouch\`{e}'s theorem, we can show that the polynomial given in equation (\ref{eq:mypoly2}) has $k$ roots on or inside the unit circle and $m$ roots outside of the unit circle.  Provided that the $m$ solutions to $y^m=\frac {k \left( \bar\lambda+\bar\mu+2\pi i n \right) }{\bar\lambda\,\left( k+m \right) }$ are not solutions to (\ref{eq:mypoly2}), the roots are distinct.  We can substitute 
\begin{equation*}%\label{eq:mroot}
y^*=\left(\frac{k(\bar\lambda+\bar\mu+2\pi i n)}{\bar\lambda(k+m)}\right)^{\frac{1}{m}}
\end{equation*}
into equation (\ref{eq:mypoly2}) to show that if $y^*$ solves (\ref{eq:mypoly2}), then
\begin{equation*}
\left( {\frac {\bar\lambda+\bar\mu+n\pi i n}{k+m}} \right) ^{k+m}=
\left( {\frac {\bar\lambda}{k}} \right) ^{k} \left( {\frac {\bar\mu}{m}}
 \right) ^{m}.
\end{equation*}
When $n=0$, 
\begin{equation*}
\left( {\frac {\bar\lambda+\bar\mu}{k+m}} \right) \ge
\left( {\frac {\bar\lambda}{k}} \right)^{\frac{k}{k+m}} \left( {\frac {\bar\mu}{m}}
 \right) ^{\frac{m}{k+m}}
\end{equation*}
by Young's inequality.  Equality holds only if $\frac{\bar\lambda}{k}=\frac{\bar\mu}{m}$, or if $k$ or $m$ is zero.  We have assumed ergodicity, so $\frac{\bar\lambda}{k}<\frac{\bar\mu}{m}$; that is, the mean arrival rate must be less than the mean service rate.  neither $k$ nor $m$ equals zero since the arrival and service processes must have at least one phase.  When $n\neq 0$, it is clear that the two sides of the equation are not equal because the real and imaginary parts are not equal.

Therefore, the roots of equation (\ref{eq:mypoly2}) are distinct.  In fact, we can use the following contraction mappings to find the $k+m$ roots for each fixed $n$.  To find the $m$ roots outside the unit circle, we may use the iteration:
\begin{equation*}
y_0^{(q)}={\rm e}^{2\pi i q/m}\left(\frac{1}{\bar\lambda}\left(2\pi i n+\bar\lambda+\bar\mu\right)\right)^{1/m}, \;q=0,\dots,m-1
\end{equation*}
with
\begin{equation*}
y_{n+1}^{(q)}={\rm e}^{2\pi i q/m}\left(\frac{1}{\bar\lambda}\left(2\pi i n+\bar\lambda+\bar\mu\right(1-(y_n^{(q)})^{-k}))\right)^{1/m}, \;q=0,\dots,m-1,
\end{equation*}
though the roots command from Matlab, for example, works perfectly well.  To find the $k$ roots on or in the unit circle, we may use the iteration:
\begin{equation*}
y_0^{(q)}={\rm e}^{2\pi i q/k}\left(\frac{\bar\mu}{2\pi i n+\bar\lambda+\bar\mu}\right)^{1/k}, \;q=0,\dots,k-1
\end{equation*}
with
\begin{equation*}
y_{n+1}^{(q)}={\rm e}^{2\pi i q/k}\left(\frac{\bar\mu}{2\pi i n+\bar\lambda(1-(y_n^{(q)})^{m}))+\bar\mu}\right)^{1/k}, \;q=0,\dots,k-1.
\end{equation*}

Because $P(z,t)$ is a generating function for an ergodic process, it must converge for all complex $|z|<1$.  This means that zeros of the denominator inside the unit circle are also zeros of the numerator of the generating function.  We focus our attention on the $m$ roots of equation (\ref{eq:mypoly2}) outside of the unit circle.  We label these roots, $\chi_{j,n}$, $j=1,\dots,m$ and $n\in\mathbb{Z}$.

We consider examples where $k$ and $m$ are relatively prime. 
Suppose $\alpha$ is a root of the polynomial (\ref{eq:mypoly2}), 
then $\omega_{km}^x \alpha$ is a root of
\[\bar\lambda\omega_k^\ell y^{k+m}-(\bar\lambda+\bar\mu+2\pi i n)y^k+\bar\mu\omega_m^j=0,\;\;n\in\mathbb{Z},\]
where $x$ is the minimum non-negative integer such that \(mj+mx\equiv 0 \;{\rm mod}\; mk\) and \(kq-kx\equiv 0\; {\rm mod}\; mk\).  Note that $\alpha$ is the $km$th root of $\chi_{j,n}$ for some $j\in\{1,\dots,m\}$.

More generally, we have 
\begin{equation}\label{eq:lim}
\lim_{y\to \omega_{km}^x \alpha}\frac{\left(1-\frac{y}{\omega_{km}^x \alpha}\right){\rm e}^{\int_u^t\epsilon_q(y^{km},\nu)+\xi_j(y^{km},\nu)d\nu}}{1-{\rm e}^{\bar\epsilon_q(y^{km})+\bar\xi_j(y^{km})}}
=\frac{{\rm e}^{\int_u^t(\lambda(\nu)(\alpha^m-1)+\mu(\nu)(\alpha^{-k}-1))d\nu}}
{m\bar\lambda \alpha^m-k\bar\mu \alpha^{-k}}
\end{equation}
independent of the indices \(j\) and \(q\). (Note that if $k$ and $m$ are not relatively prime, the  approach in this paper can still be used, but the limit given in equation (\ref{eq:lim}) would not be independent of $j$ and $q$.  We would need to find the roots of more than one equation for each \(n\).) So, we approximate $\frac{{\rm e}^{\int_u^t\epsilon_q(y^{km},\nu)+\xi_j(y^{km},\nu)d\nu}}{1-{\rm e}^{\bar\epsilon_q(y^{km})+\bar\xi_j(y^{km})}}$ with the series
\begin{equation*}
\frac{{\rm e}^{\int_u^t(\lambda(\nu)(\alpha^{m}-1)+\mu(\nu)(\alpha^{-k}-1))d\nu}}
{m\bar\lambda \alpha^m-k\bar\mu \alpha^{-k}}
\sum_{n=0}^\infty \frac{y^{n}}{\omega^{nx}_{km}\alpha^{n}}.
\end{equation*}

Define 
\begin{equation*}
D_\chi=\left[\begin{array}{cccc}1&\chi^{-1/k}&\cdots&\chi^{(1-k)/k}\\\chi^{1/k}&1&\cdots&\chi^{(2-k)/k}\\\vdots&\ddots&\ddots&\vdots\\\chi^{(k-1)/k}&\chi^{(k-2)/k}&\cdots&1\end{array}\right]
\end{equation*}
 and 
\begin{equation*}
C_\chi=\left[\begin{array}{cccc}1&\chi^{1/m}&\cdots&\chi^{(m-1)/m}\\\chi^{-1/m}&1&\cdots&\chi^{(m-2)/m}\\\vdots&\ddots&\ddots&\vdots\\\chi^{(1-m)/m}&\chi^{(2-m)/m}&\cdots&1\end{array}\right].
\end{equation*}
Then we may express the generating function given in equation (\ref{eq:EkEmkey}), as
\begin{multline}\label{eq:gf}
{\bf P}(z,t)=\sum_{j=1}^\infty {\bf p}_j(t)z^j=\\
\sum_{j=1}^\infty\int_{t-1}^t\sum_{n=-\infty}^\infty\sum_{\ell=1}^m
\frac{{\rm e}^{\int_u^t(\lambda(\nu)(\chi_{\ell,n}^{1/k}-1)+\mu(\nu)(\chi_{\ell,n}^{-1/m}-1))d\nu}}
{m\bar\lambda \chi_{\ell,n}^{1/k}-k\bar\mu \chi_{\ell,n}^{-1/m}}\\
\times
\left({\bf p}_0(u)\chi_{\ell,n}{\bf Q}_{0,1}(u)-{\bf p}_1(u){\bf A}_{-1}(u)\right)du\\
\times
D_{\chi_{\ell,n}}\otimes C_{\chi_{\ell,n}}\left(\frac{z^j}{\chi_{\ell,n}^j}\right)\\
=\sum_{j=1}^\infty\int_{t-1}^t\sum_{n=-\infty}^\infty\sum_{\ell=1}^m
\frac{{\rm e}^{\int_u^t(\lambda(\nu)(\chi_{\ell,n}^{1/k}-1)+\mu(\nu)(\chi_{\ell,n}^{-1/m}-1))d\nu}}{m\bar\lambda \chi_{\ell,n}^{1/k}-k\bar\mu \chi_{\ell,n}^{-1/m}}\\
\times
\left(p_{0,k-1}(u)\chi_{\ell,n}\lambda(u)-\mu(u)\sum_{q=0}^{k-1}p_{1,q,m-1}(u)\chi_{\ell,n}^{q/k}\right)du\\
\times\left[\begin{array}{cccc}1&\chi_{\ell,n}^{-1/k}&\cdots&\chi_{\ell,n}^{(1-k)/k}\end{array}\right]\otimes\left[\begin{array}{cccc}1&\chi_{\ell,n}^{1/m}&\cdots&\chi_{\ell,n}^{(m-1)/m}\end{array}\right]\left(\frac{z^j}{\chi_{\ell,n}^j}\right)
\end{multline}
Define
\begin{multline}\label{eq:fxt}
f(x,t)=\int_{t-1}^t\frac{{\rm e}^{\int_u^t(\lambda(\nu)(x^{1/k}-1)+\mu(\nu)(x^{-1/m}-1))d\nu}}{m\bar\lambda x^{1/k}-k\bar\mu x^{-1/m}}\\
\times
\left(p_{0,k-1}(u)x\lambda(u)-\mu(u)\sum_{q=0}^{k-1}p_{1,q,m-1}(u)x^{q/k}\right)du, 
\end{multline}
then the probability vector for level $j$ is $[z^n]{\bf P}(z,t)$ from equation (\ref{eq:gf}),
\begin{equation}\label{eq:pjt}
{\bf p}_j(t)=\sum_{n=-\infty}^\infty\sum_{\ell=1}^m f(\chi_{\ell,n},t)\chi_{\ell,n}^{-j}\left[\begin{array}{cccc}1&\chi_{\ell,n}^{-1/k}&\cdots&\chi_{\ell,n}^{(1-k)/k}\end{array}\right]\otimes\left[\begin{array}{cccc}1&\chi_{\ell,n}^{1/m}&\cdots&\chi_{\ell,n}^{(m-1)/m}\end{array}\right].
\end{equation}
This expression is exact.  See \cite{margolius_2021} for more details.  

To illustrate the method, we consider an example of an $E_7/E_4/1$ queue with 
\begin{equation*}
\lambda(t)=3-2\sin(2\pi t)
\end{equation*}
and
\begin{equation*}
\mu(t)=5+4\sin(2\pi t).
\end{equation*}
We approximate the distribution with 
\begin{equation}\label{eq:pjest}
{\bf p}^{(q)}_j(t)=\sum_{n=-q}^q\sum_{\ell=1}^m f(\chi_{\ell,n},t)\chi_{\ell,n}^{-j}\left[\begin{array}{cccc}1&\chi_{\ell,n}^{-1/k}&\cdots&\chi_{\ell,n}^{(1-k)/k}\end{array}\right]\otimes\left[\begin{array}{cccc}1&\chi_{\ell,n}^{1/m}&\cdots&\chi_{\ell,n}^{(m-1)/m}\end{array}\right].
\end{equation}

%In figure \ref{fig:level1q1}, we have taken $q=1$ in equation (\ref{eq:pjest}).
\begin{figure}[!tbp]
\begin{tabular}{cccc}
\subfloat[$p_{1,0,0}$]{\includegraphics[width=0.21\textwidth]{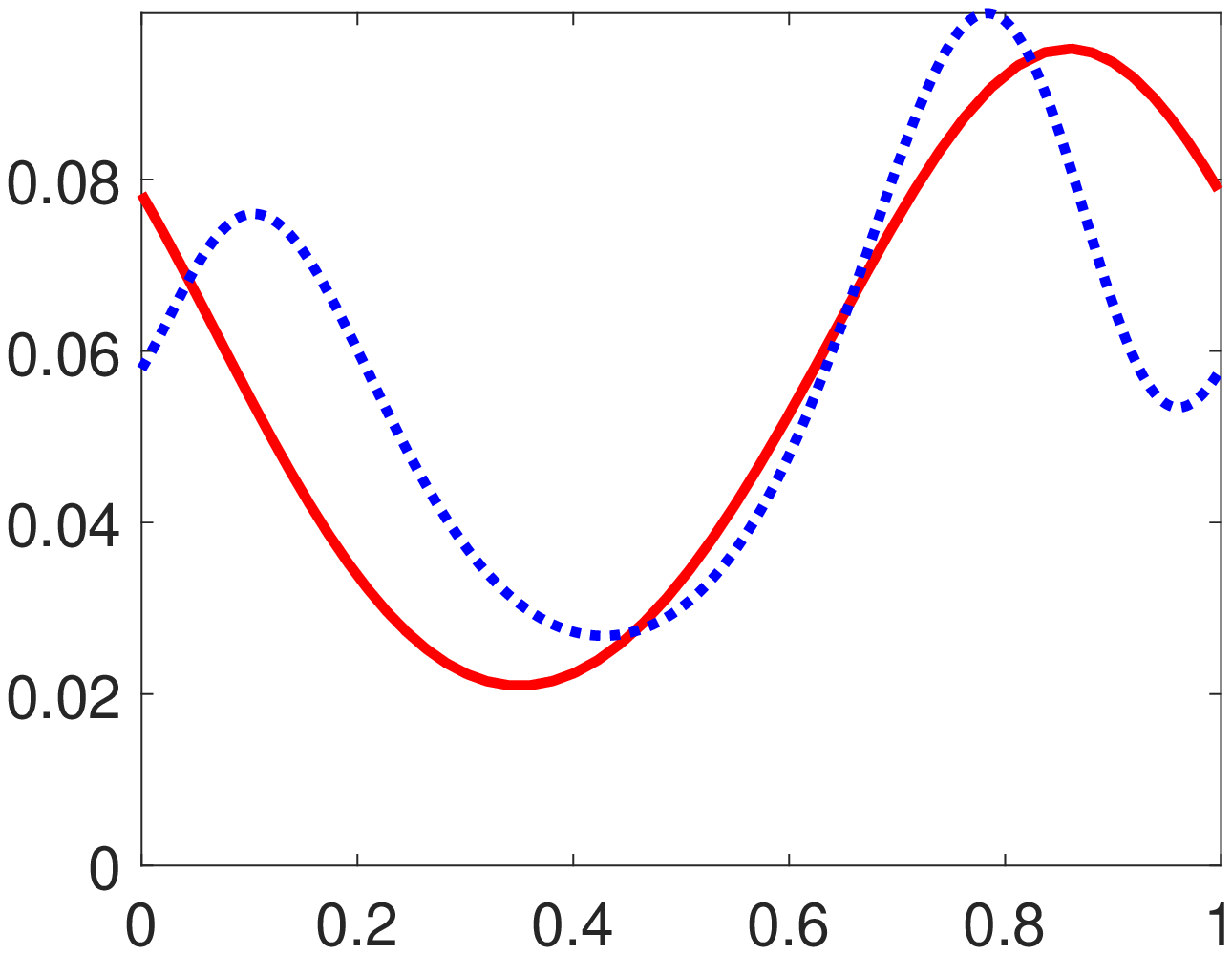}} &
\subfloat[$p_{1,0,1}$]{\includegraphics[width=0.21\textwidth]{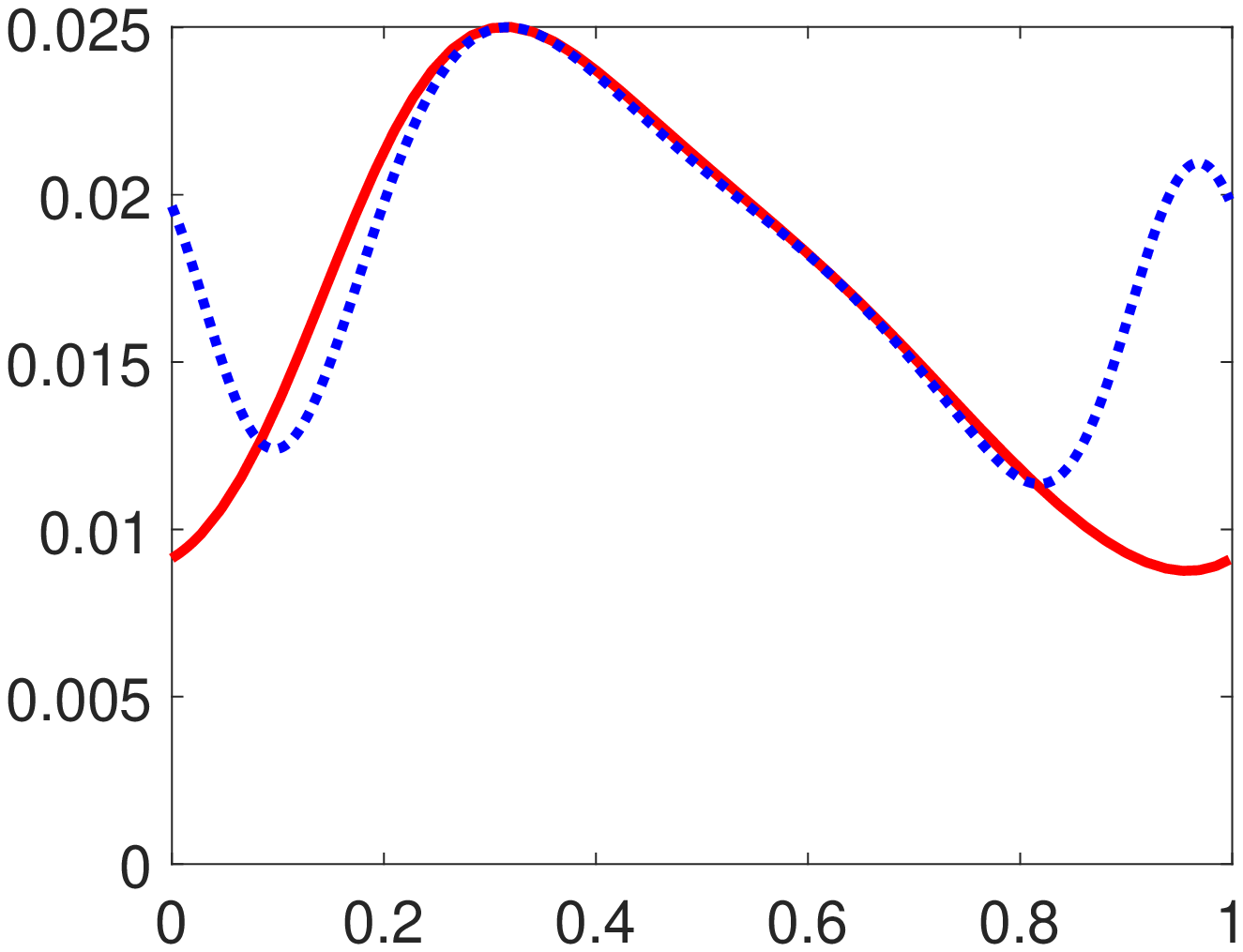}} &
\subfloat[$p_{1,0,2}$]{\includegraphics[width=0.21\textwidth]{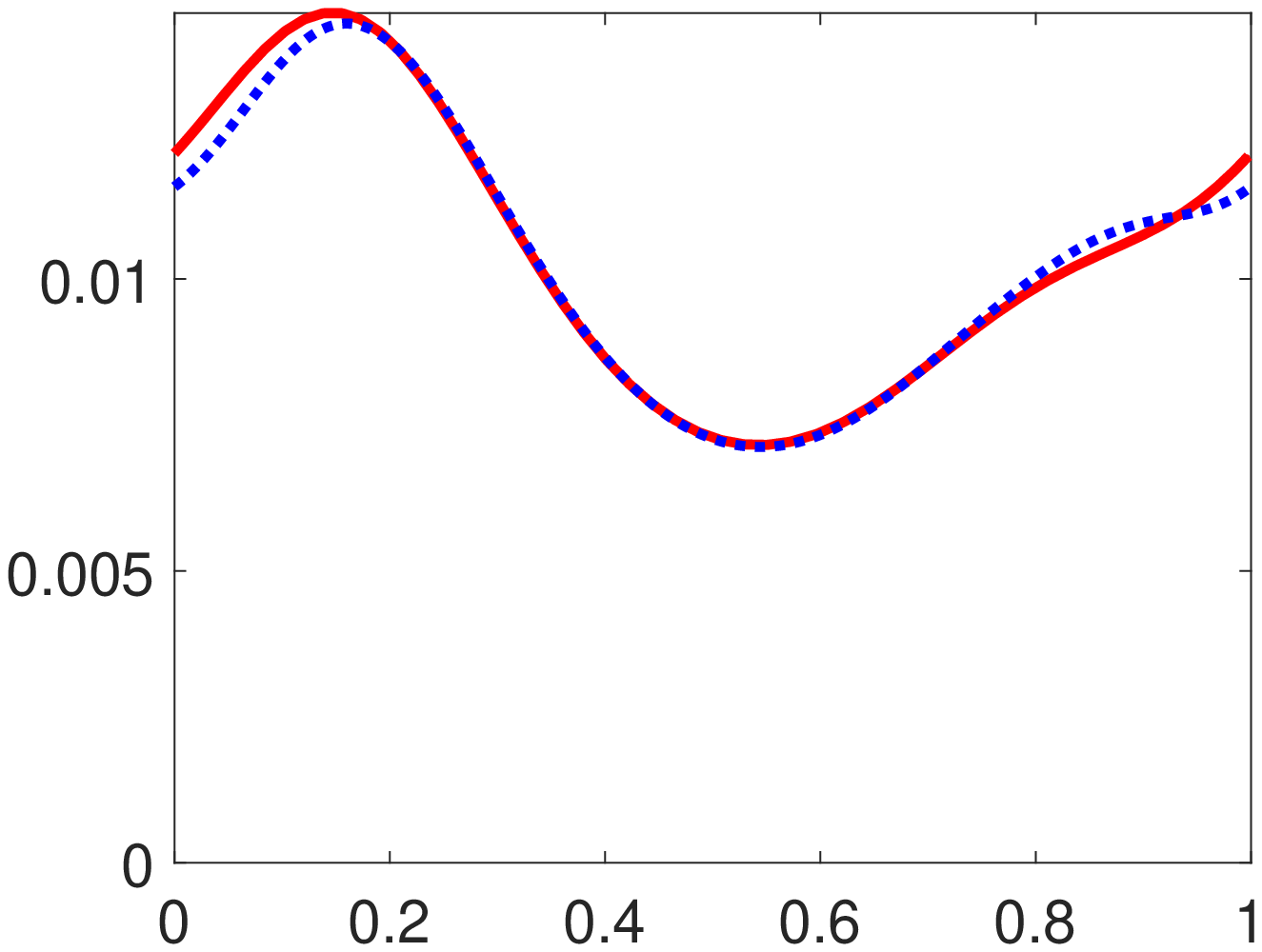}} &
\subfloat[$p_{1,0,3}$]{\includegraphics[width=0.21\textwidth]{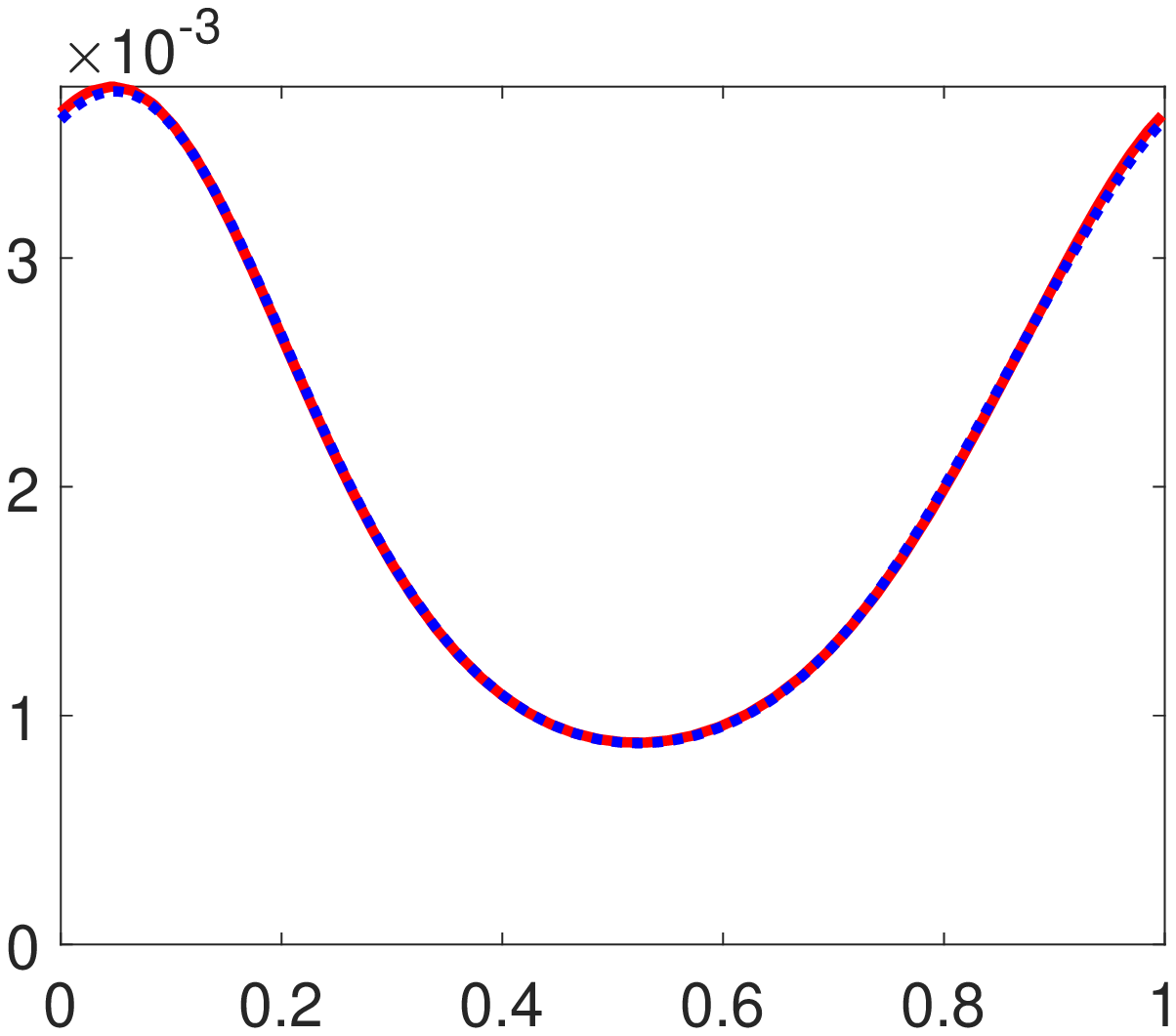}}\\
\subfloat[$p_{1,1,0}$]{\includegraphics[width=0.21\textwidth]{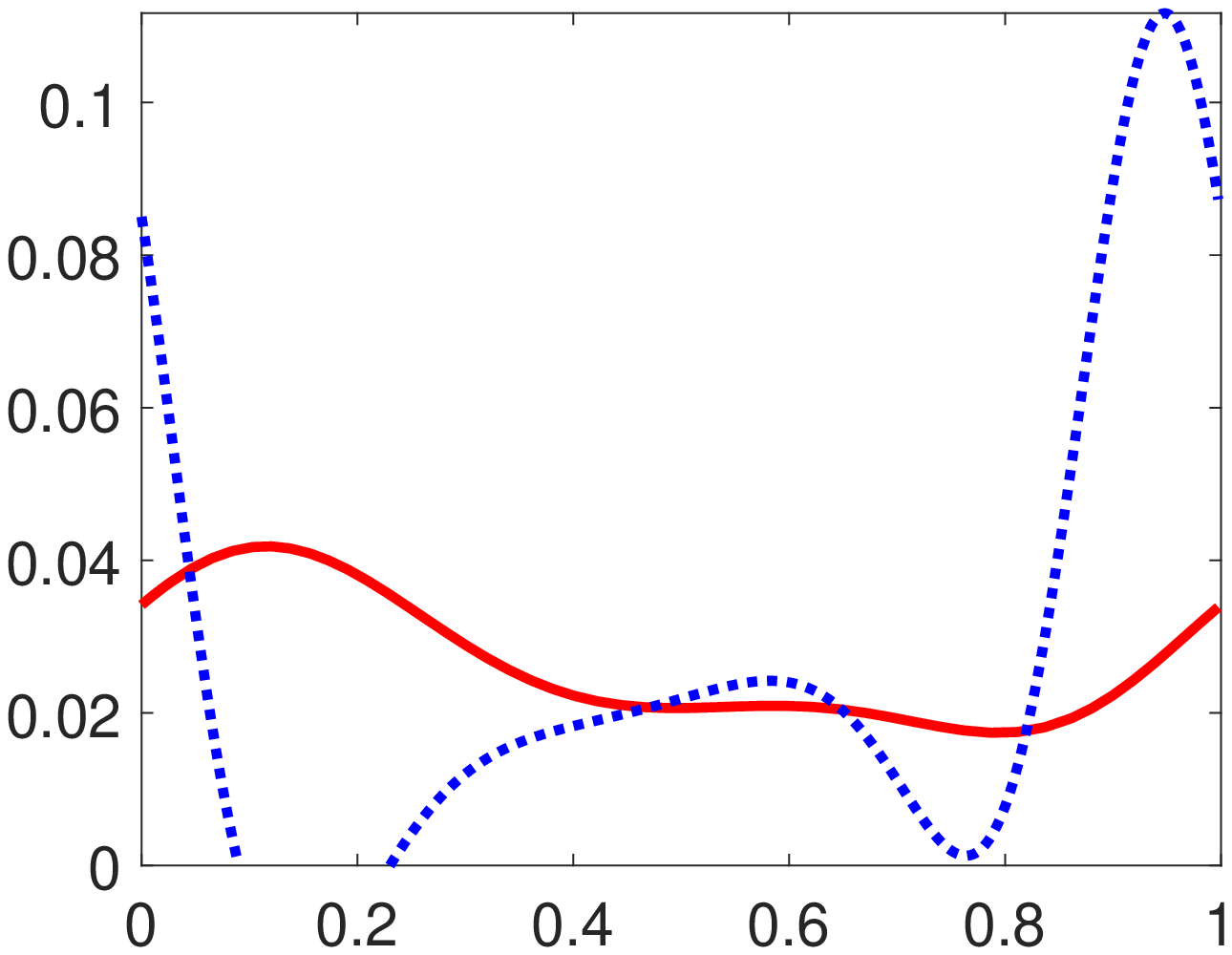}} &
\subfloat[$p_{1,1,1}$]{\includegraphics[width=0.21\textwidth]{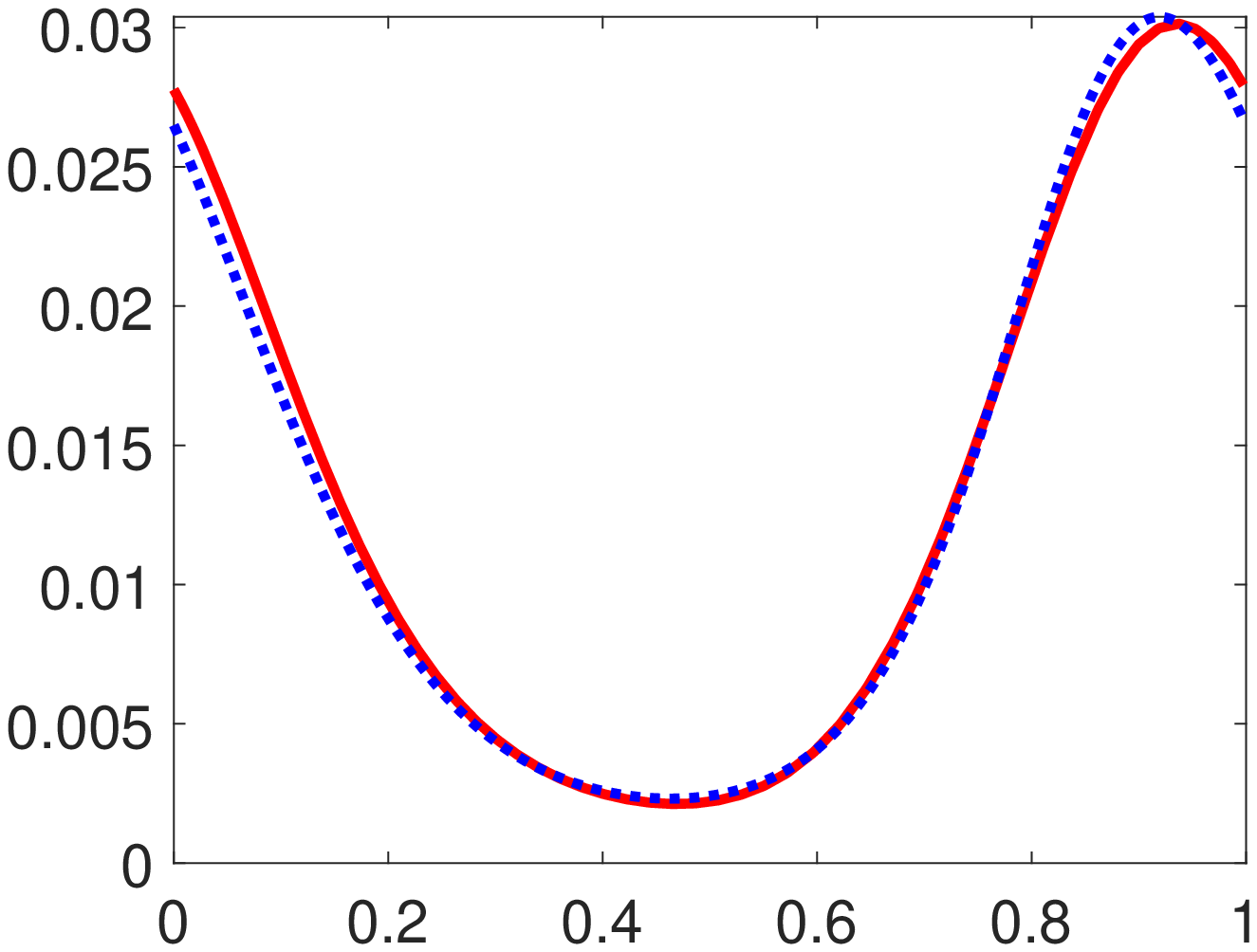}} &
\subfloat[$p_{1,1,2}$]{\includegraphics[width=0.21\textwidth]{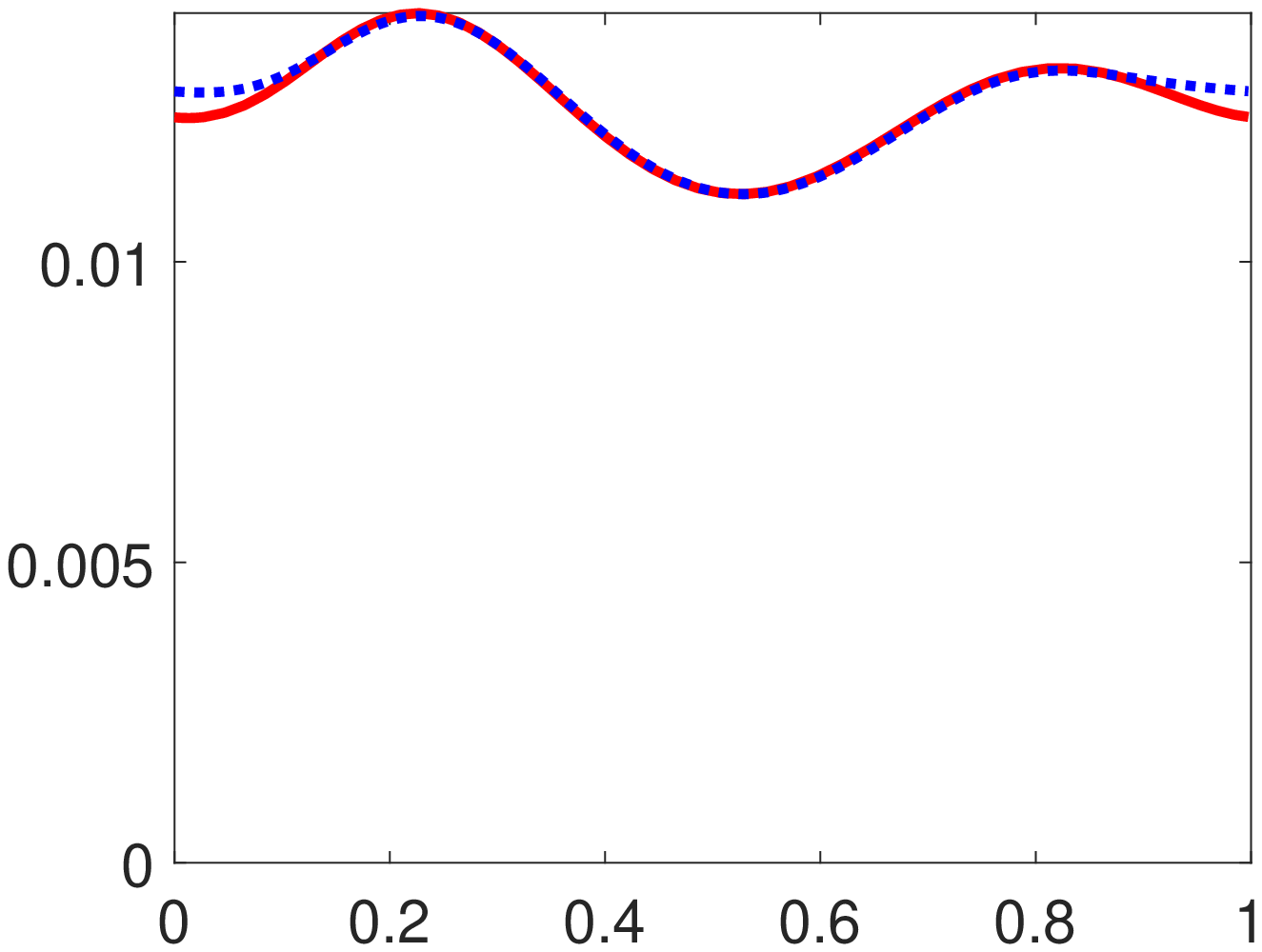}} &
\subfloat[$p_{1,1,3}$]{\includegraphics[width=0.21\textwidth]{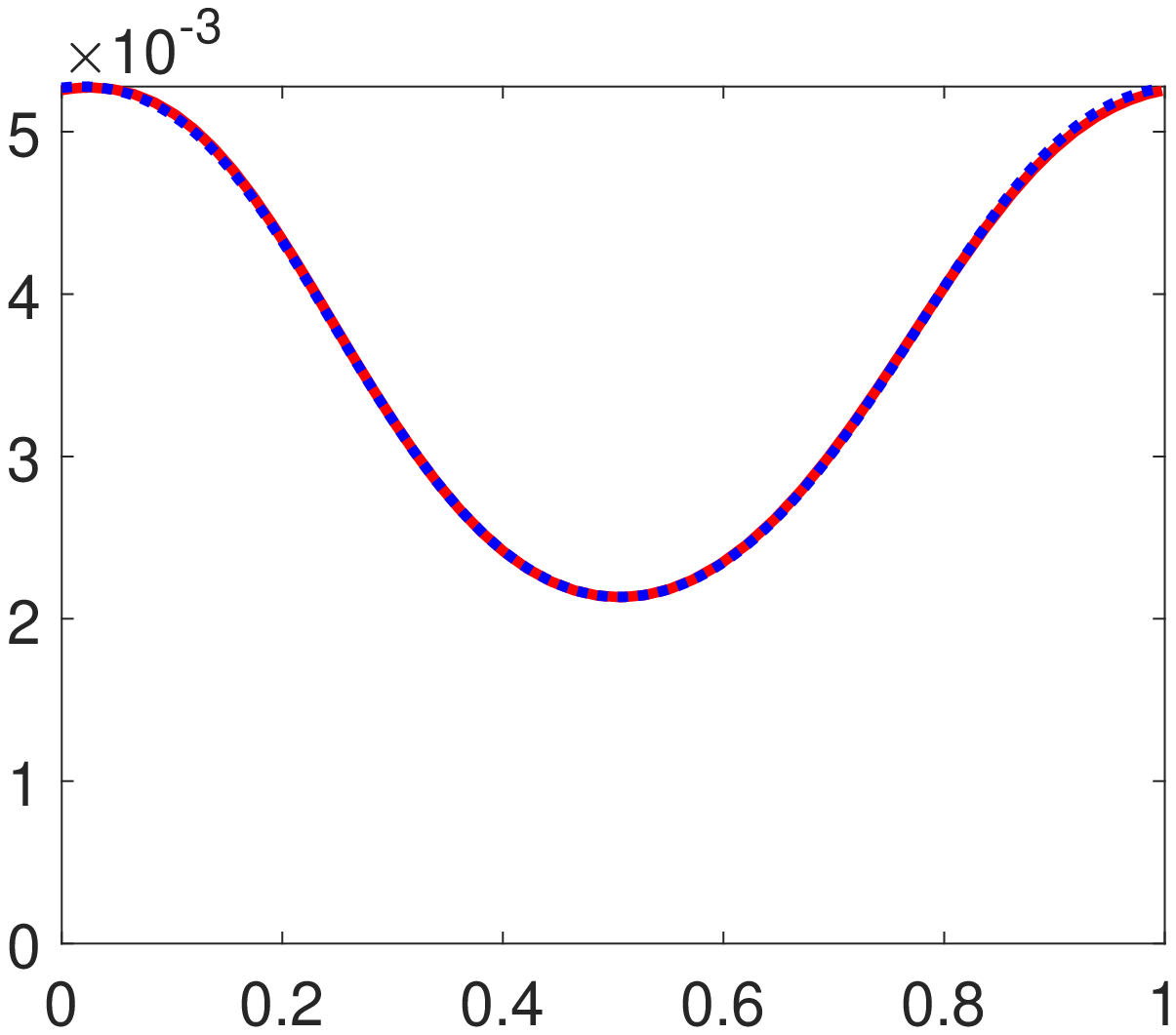}}\\
\subfloat[$p_{1,2,0}$]{\includegraphics[width=0.21\textwidth]{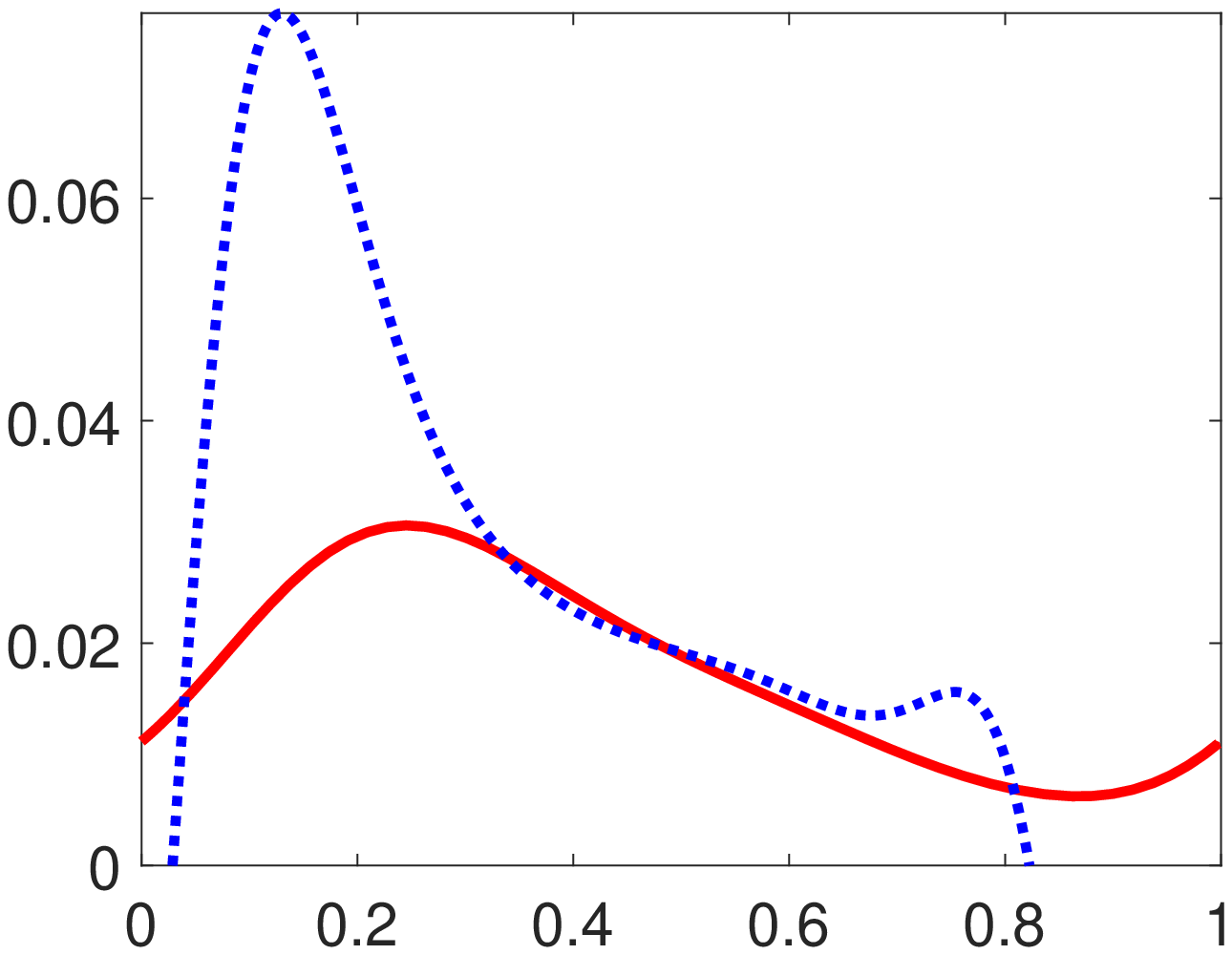}} &
\subfloat[$p_{1,2,1}$]{\includegraphics[width=0.21\textwidth]{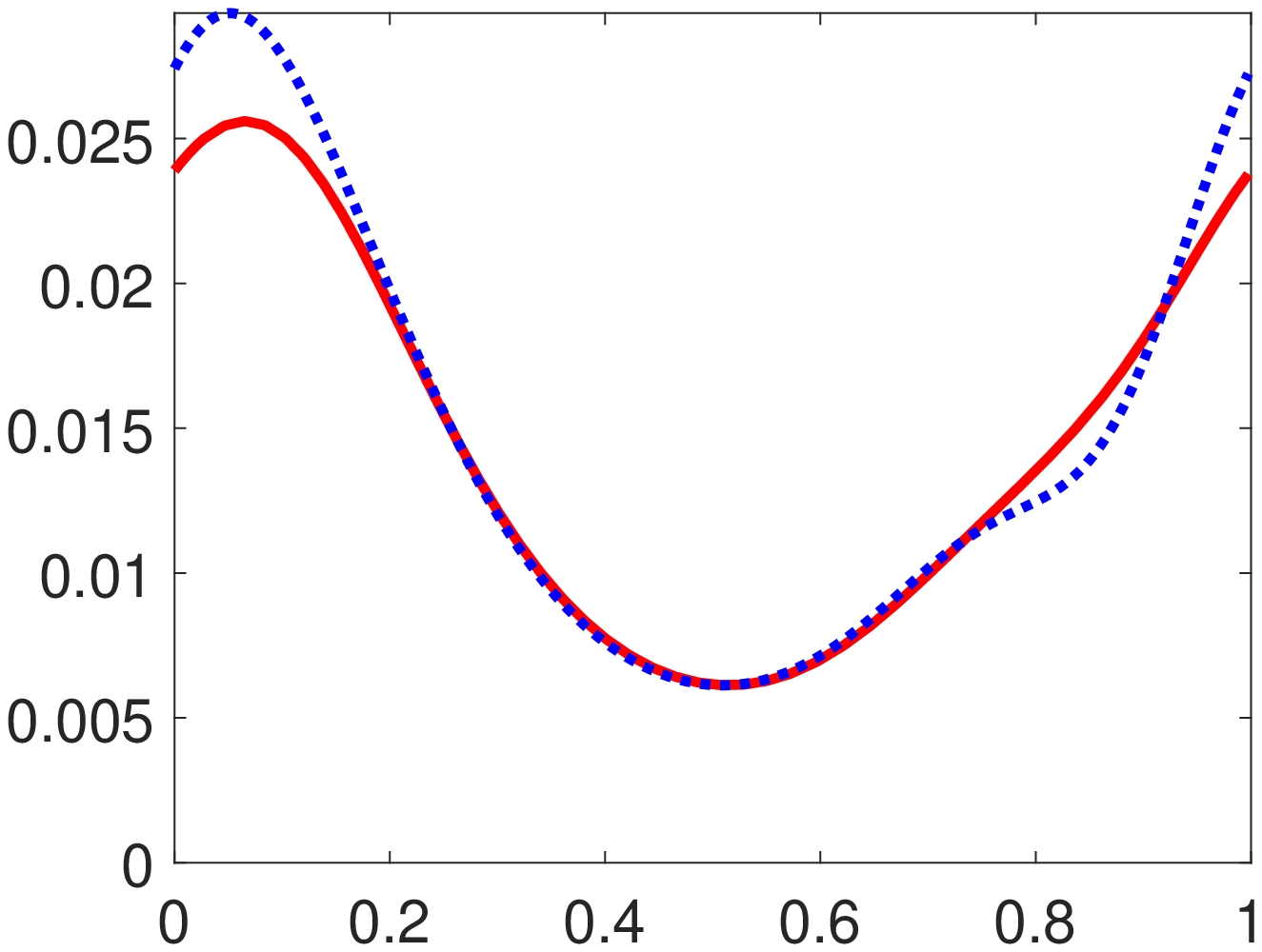}} &
\subfloat[$p_{1,2,2}$]{\includegraphics[width=0.21\textwidth]{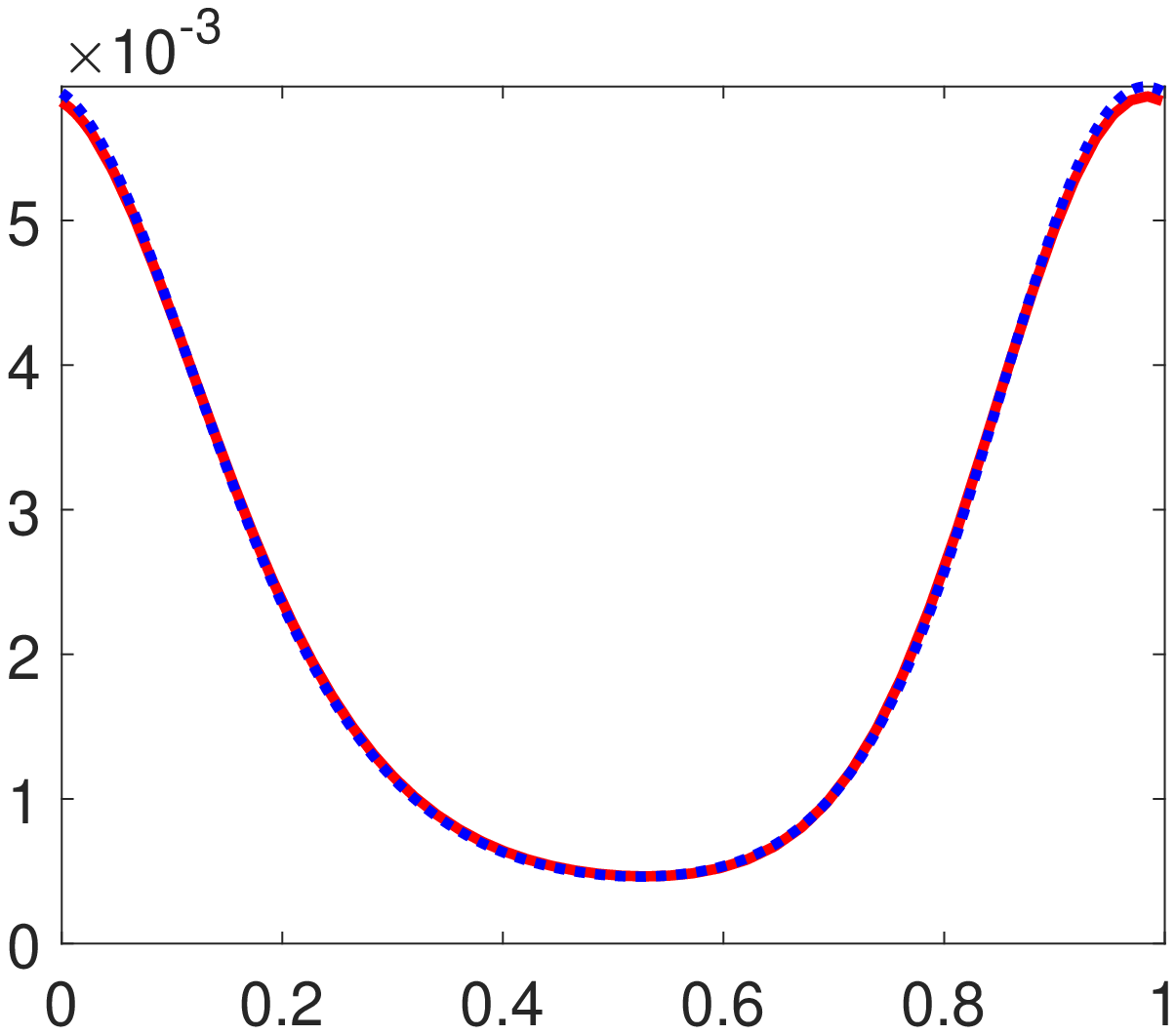}} &
\subfloat[$p_{1,2,3}$]{\includegraphics[width=0.21\textwidth]{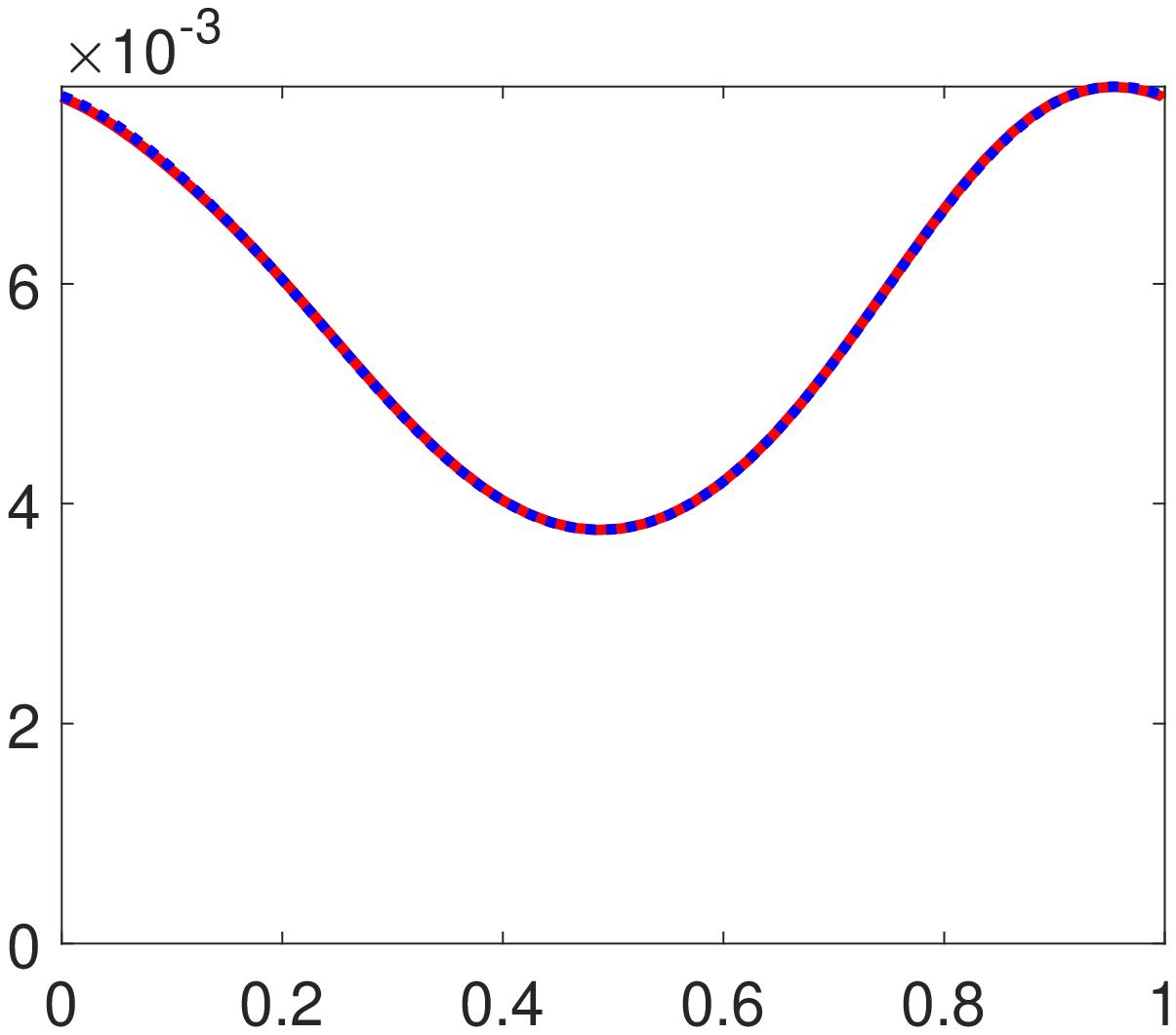}}\\
\subfloat[$p_{1,3,0}$]{\includegraphics[width=0.21\textwidth]{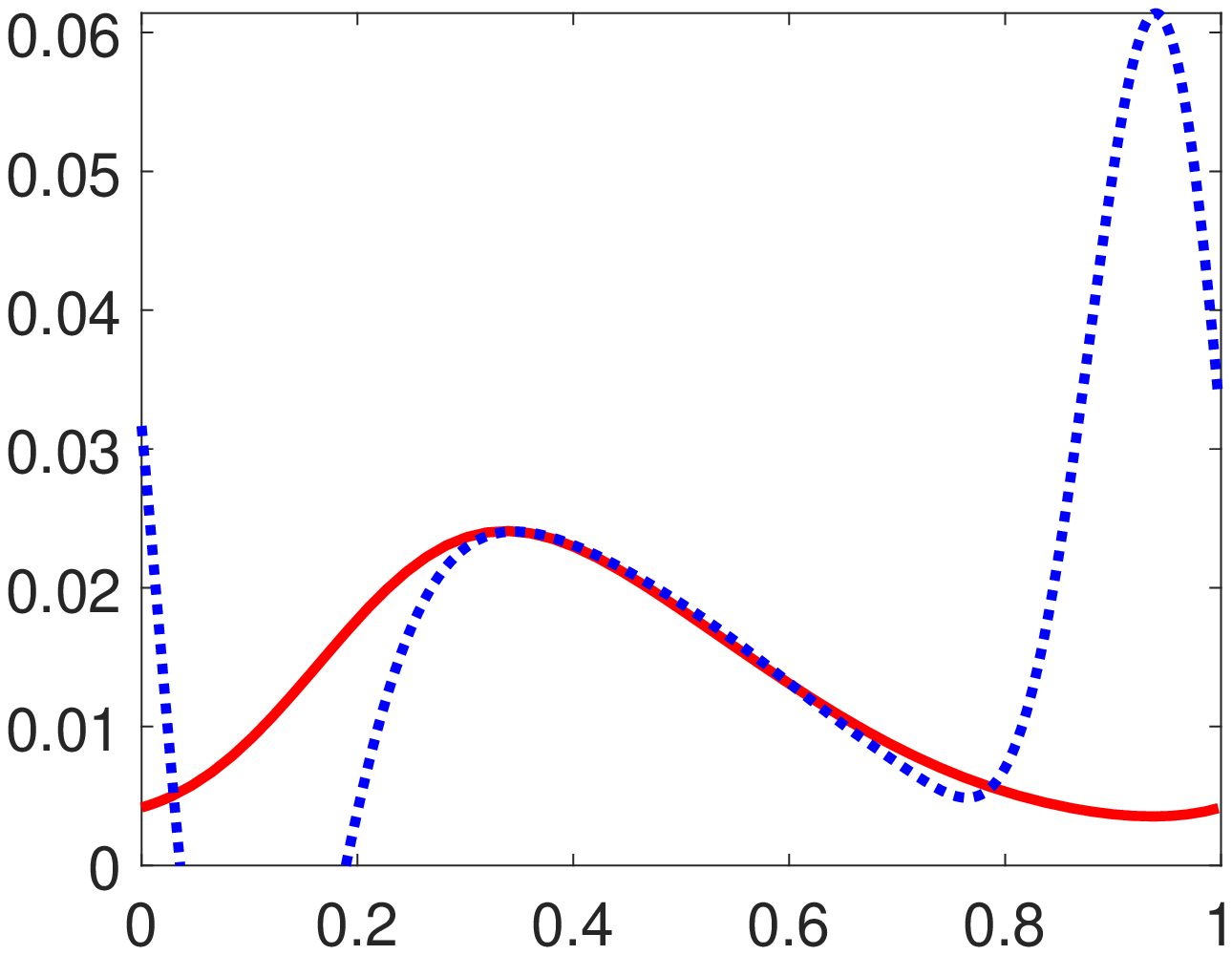}} &
\subfloat[$p_{1,3,1}$]{\includegraphics[width=0.21\textwidth]{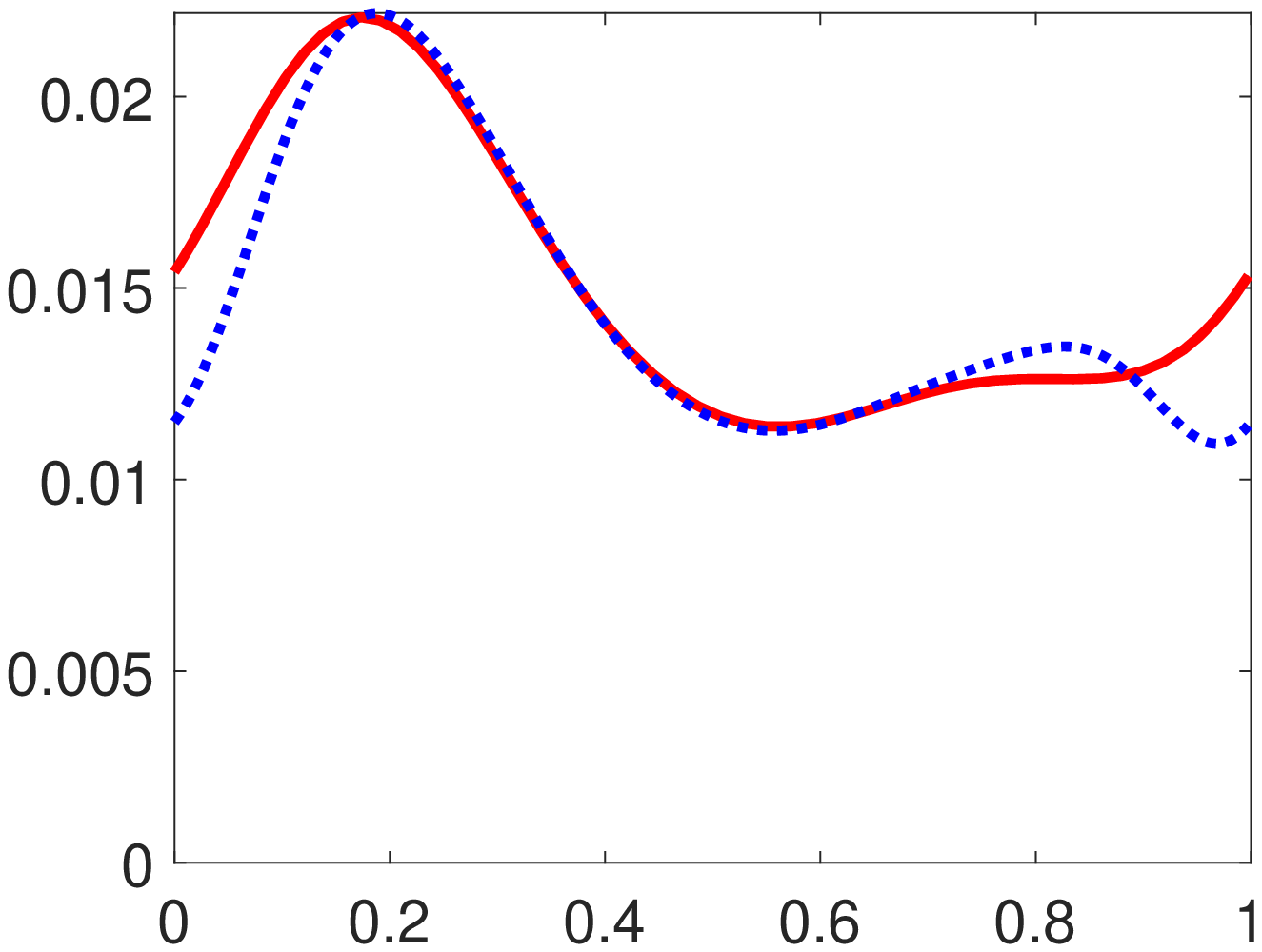}} &
\subfloat[$p_{1,3,2}$]{\includegraphics[width=0.21\textwidth]{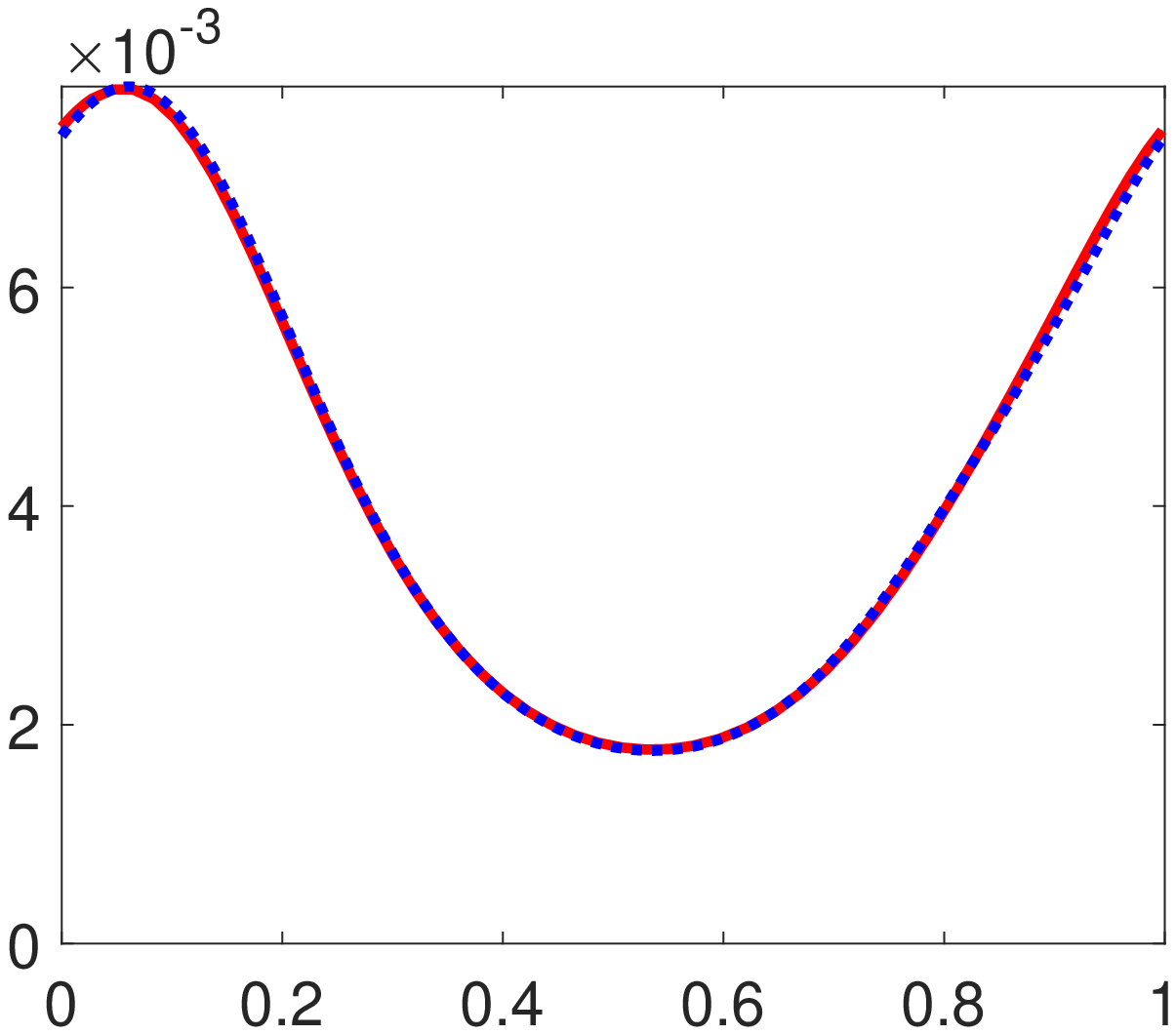}} &
\subfloat[$p_{1,3,3}$]{\includegraphics[width=0.21\textwidth]{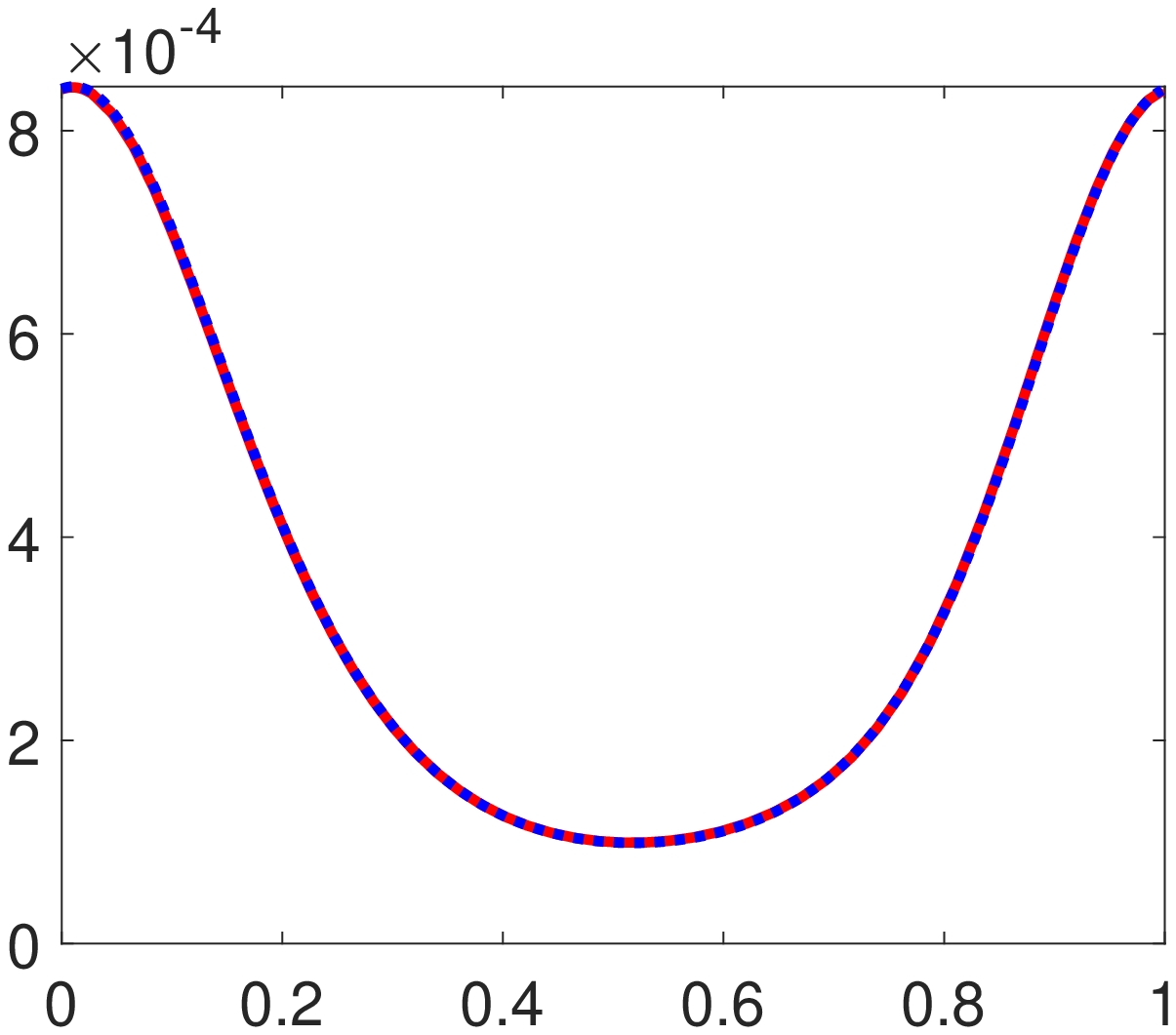}}\\
\subfloat[$p_{1,4,0}$]{\includegraphics[width=0.21\textwidth]{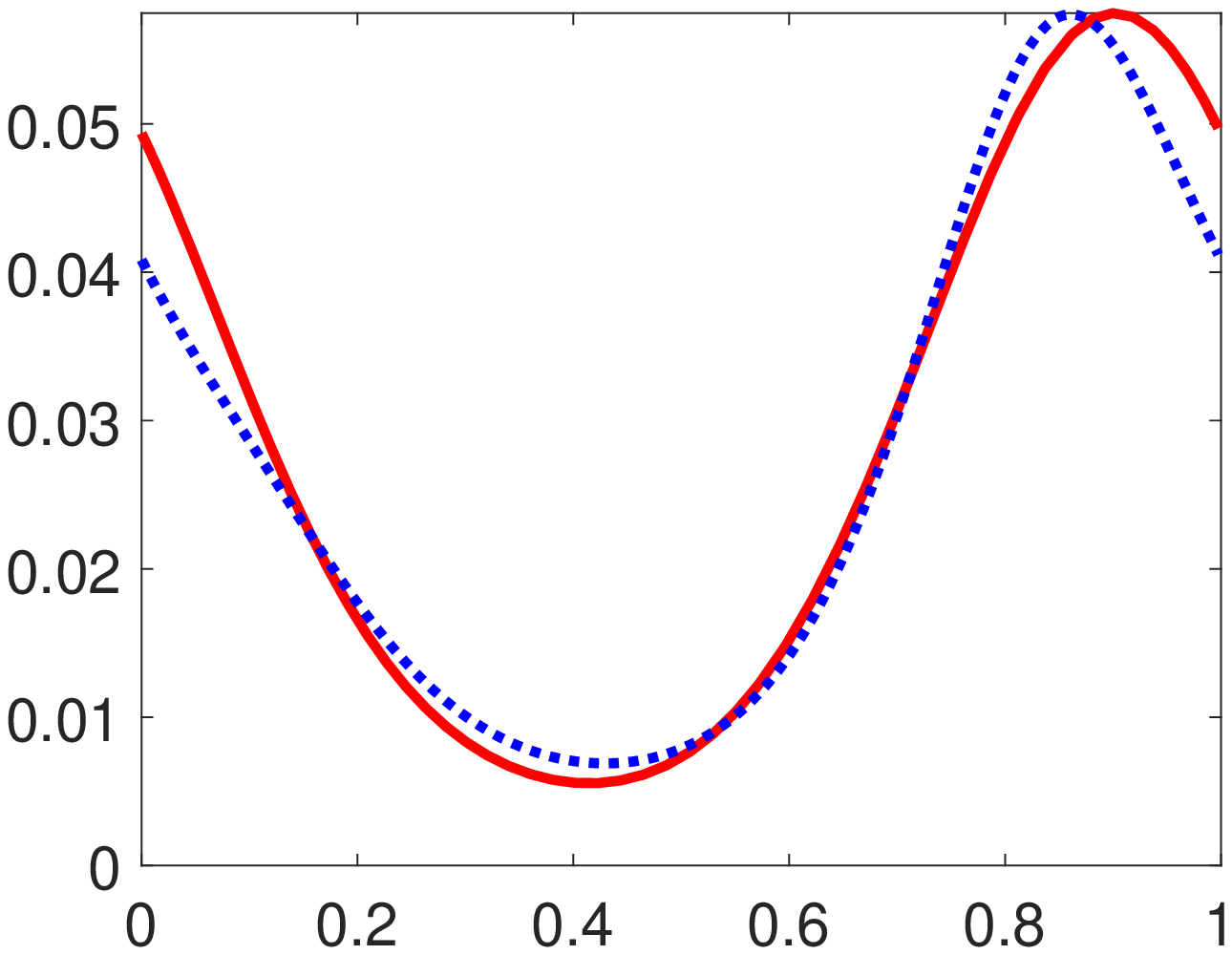}} &
\subfloat[$p_{1,4,1}$]{\includegraphics[width=0.21\textwidth]{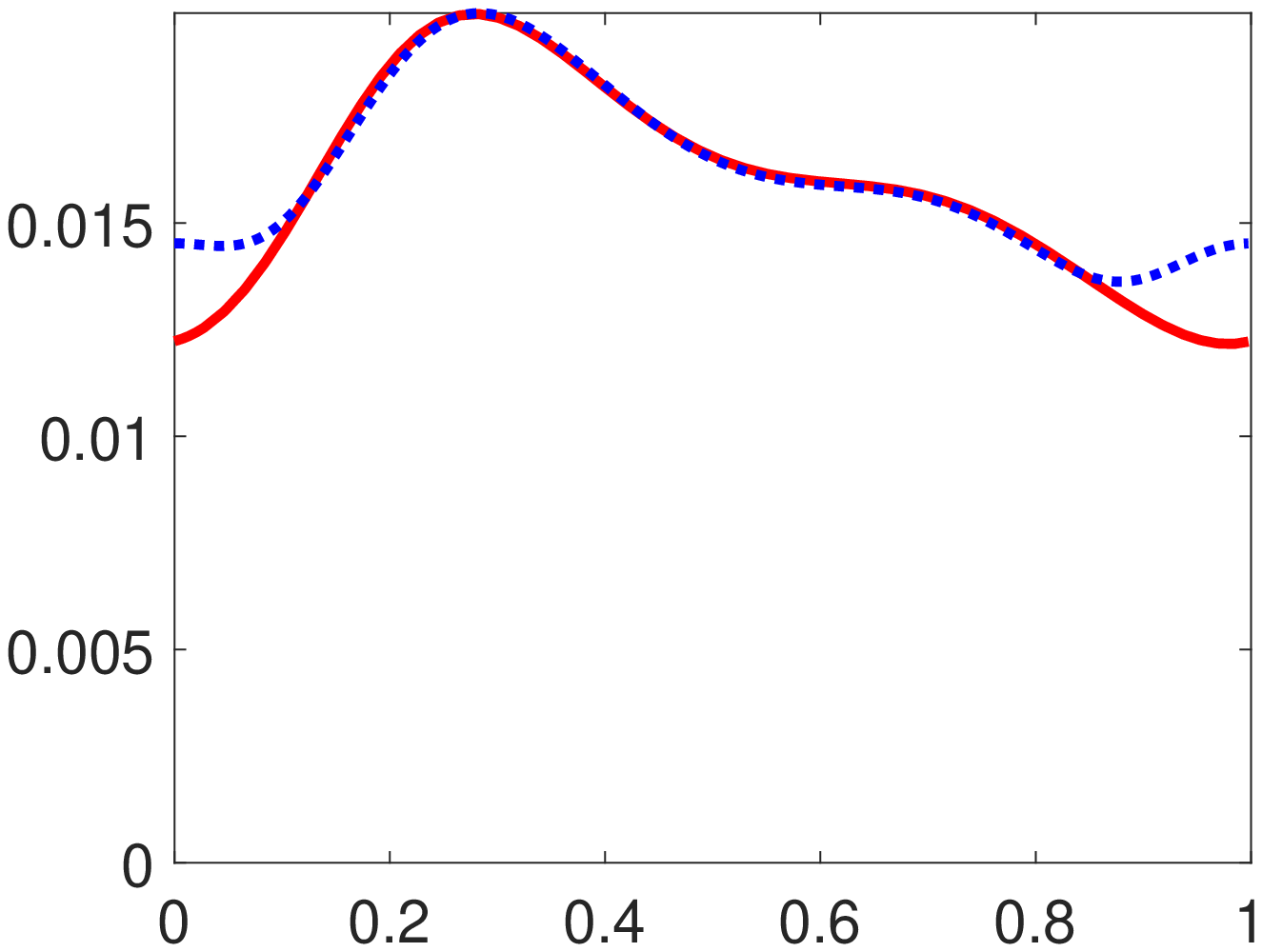}} &
\subfloat[$p_{1,4,2}$]{\includegraphics[width=0.21\textwidth]{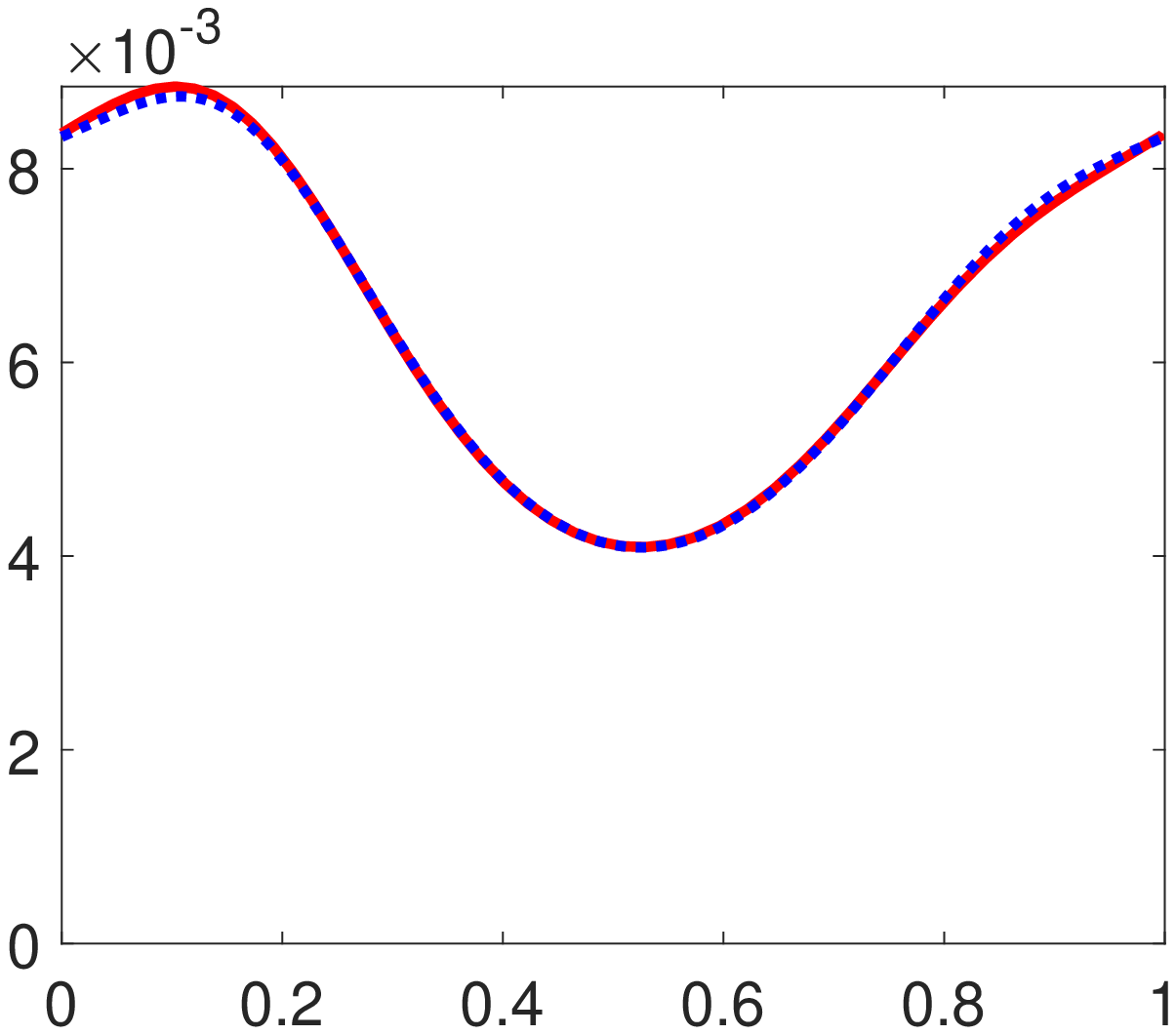}} &
\subfloat[$p_{1,4,3}$]{\includegraphics[width=0.21\textwidth]{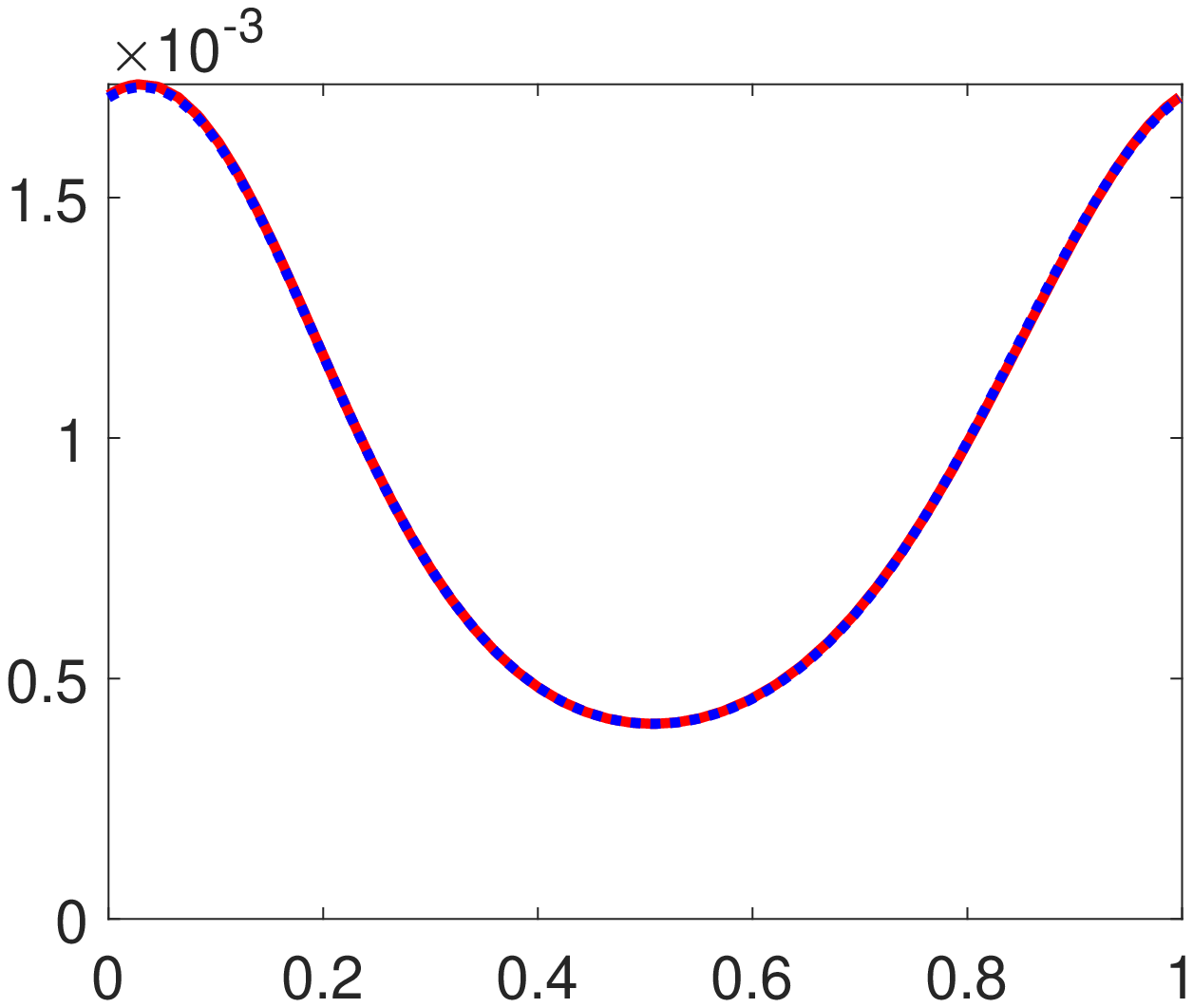}}\\
\subfloat[$p_{1,5,0}$]{\includegraphics[width=0.21\textwidth]{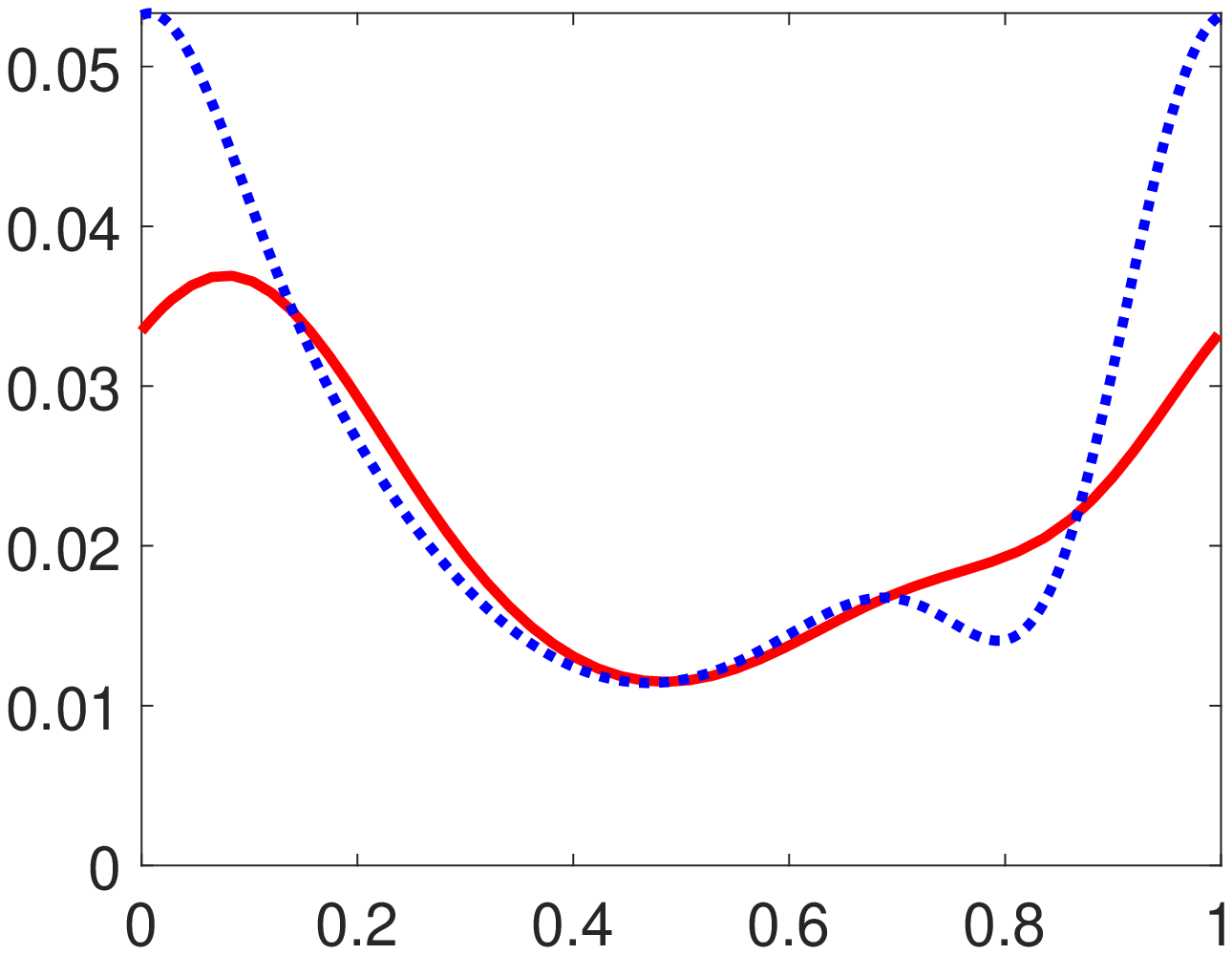}} &
\subfloat[$p_{1,5,1}$]{\includegraphics[width=0.21\textwidth]{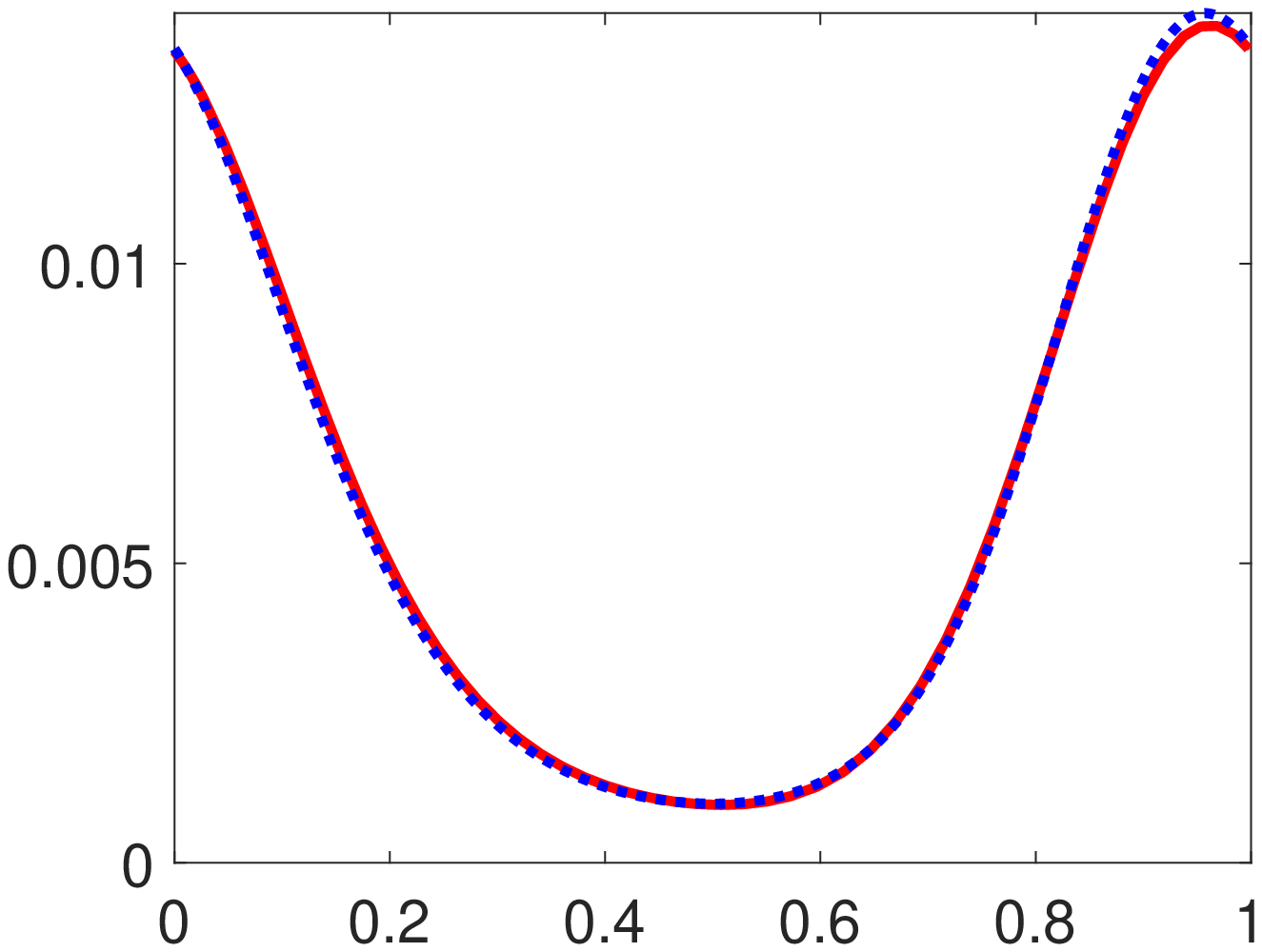}} &
\subfloat[$p_{1,5,2}$]{\includegraphics[width=0.21\textwidth]{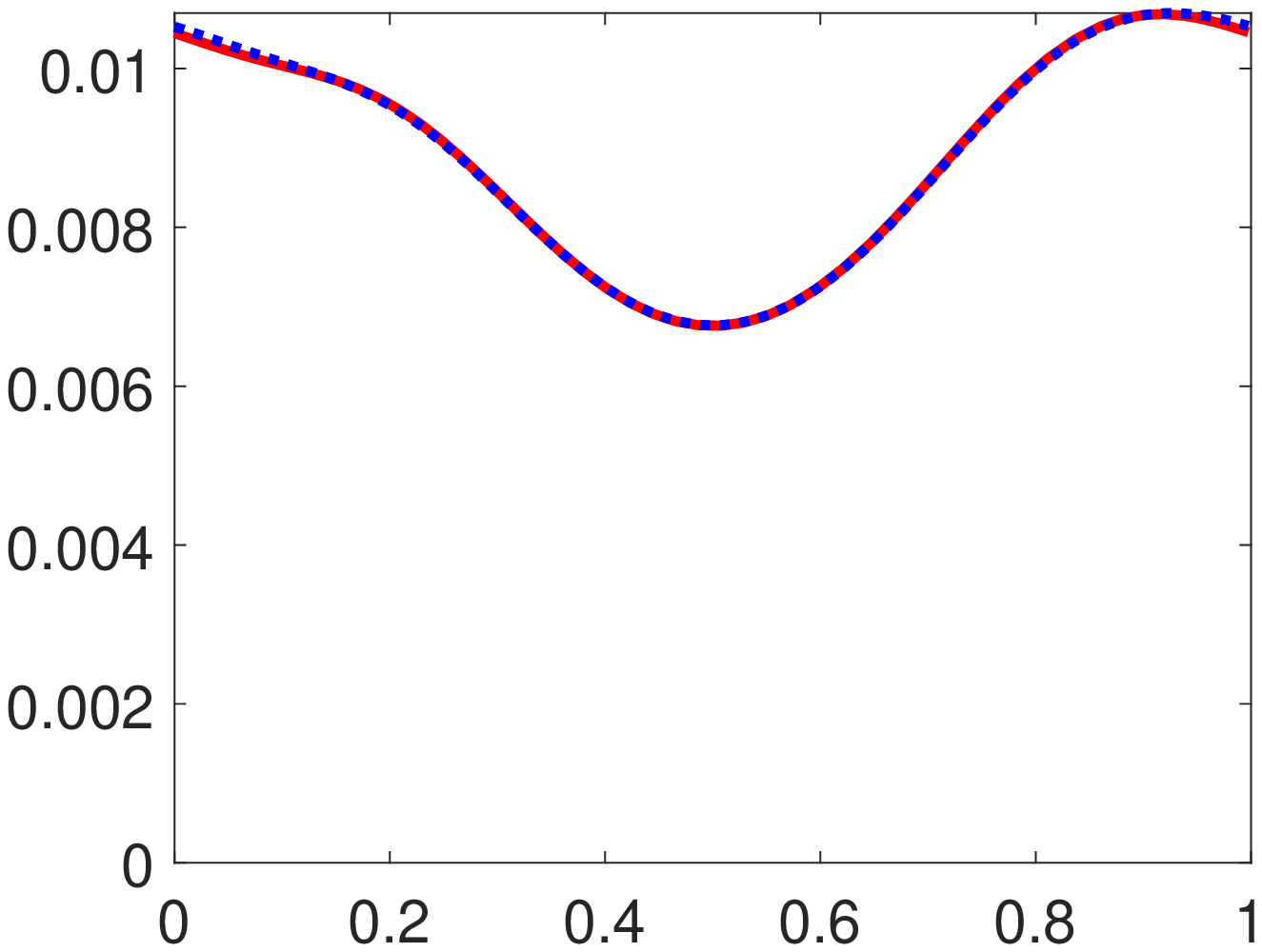}} &
\subfloat[$p_{1,5,3}$]{\includegraphics[width=0.21\textwidth]{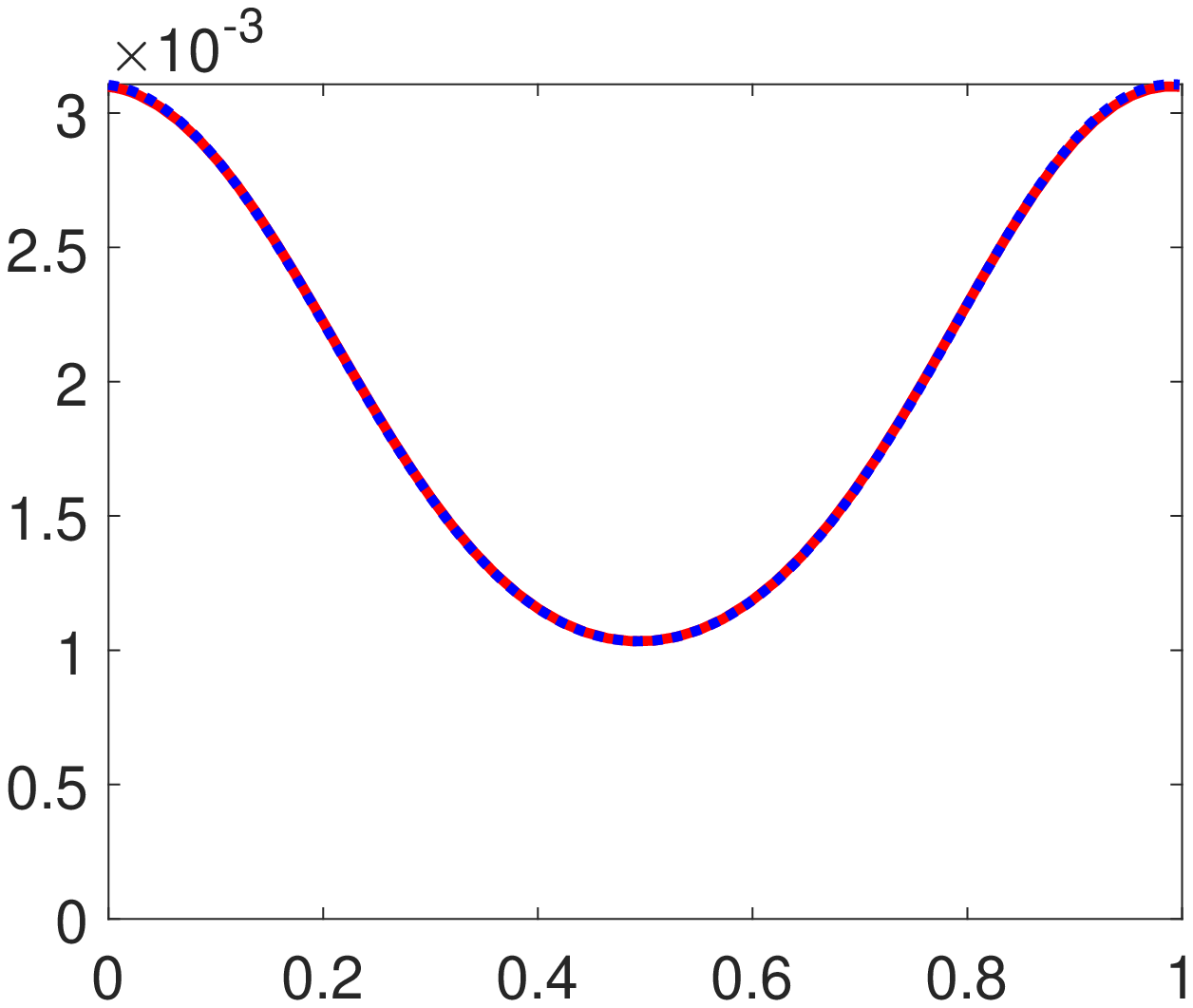}}\\
\subfloat[$p_{1,6,0}$]{\includegraphics[width=0.21\textwidth]{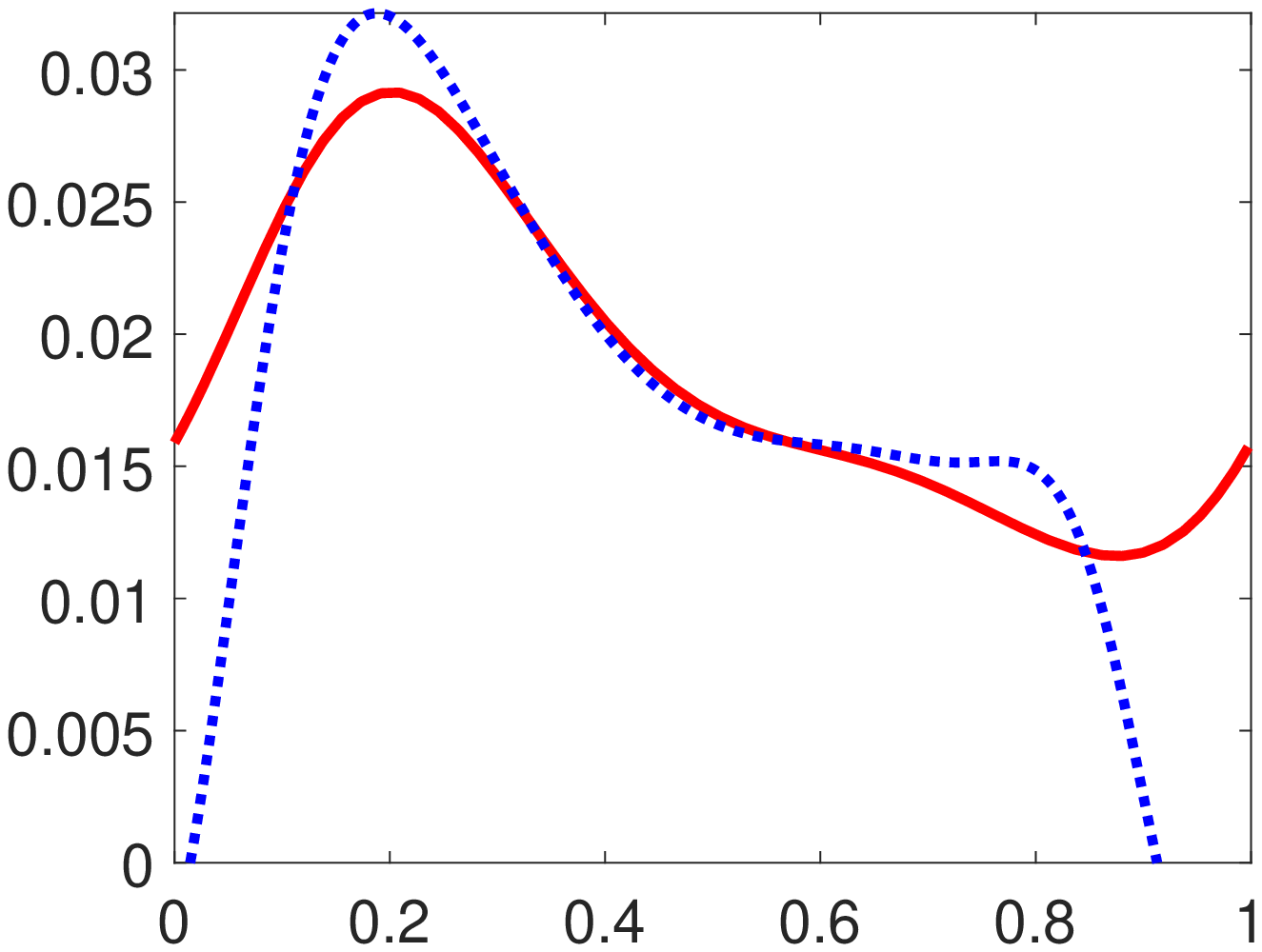}} &
\subfloat[$p_{1,6,1}$]{\includegraphics[width=0.21\textwidth]{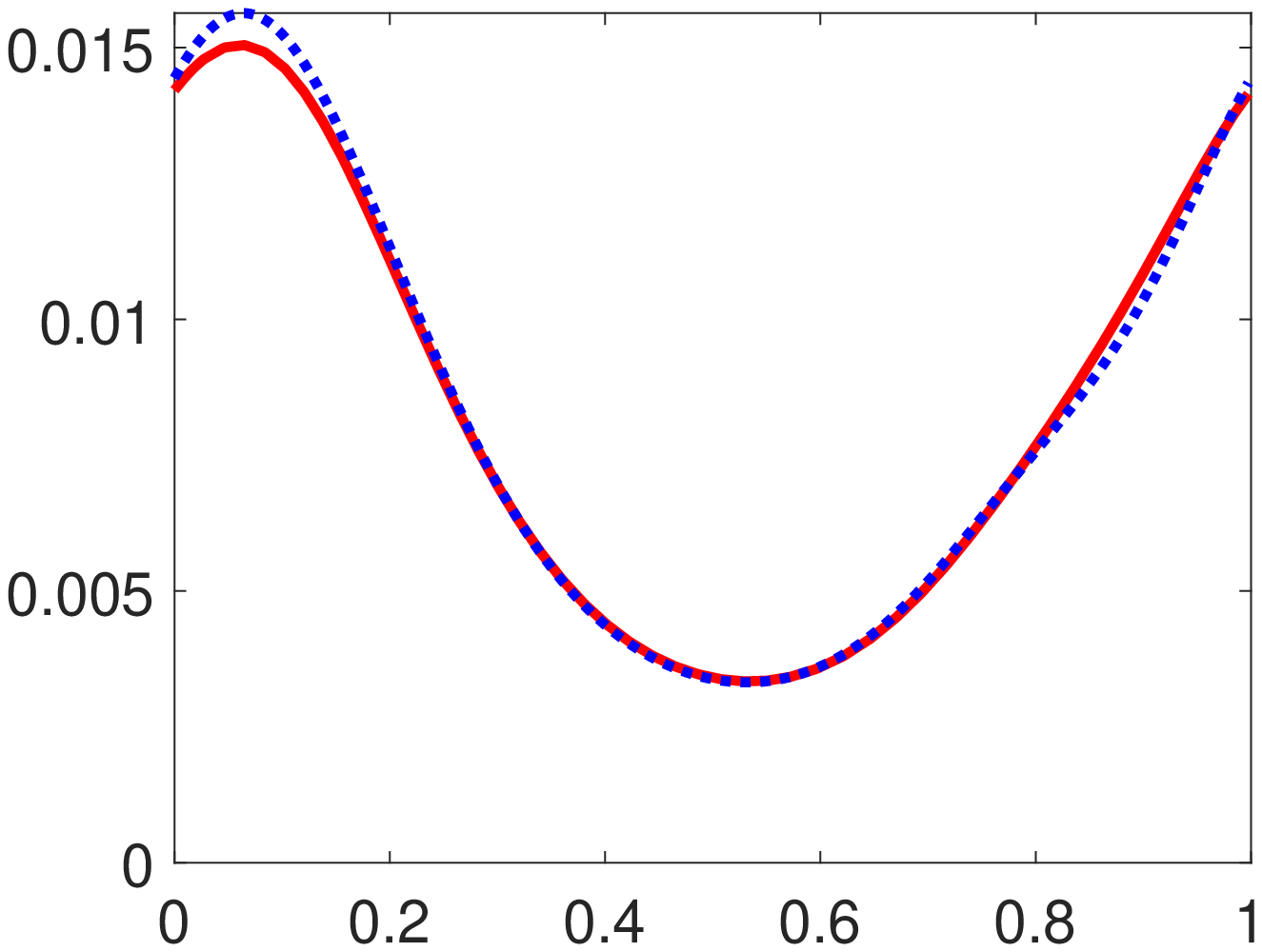}} &
\subfloat[$p_{1,6,2}$]{\includegraphics[width=0.21\textwidth]{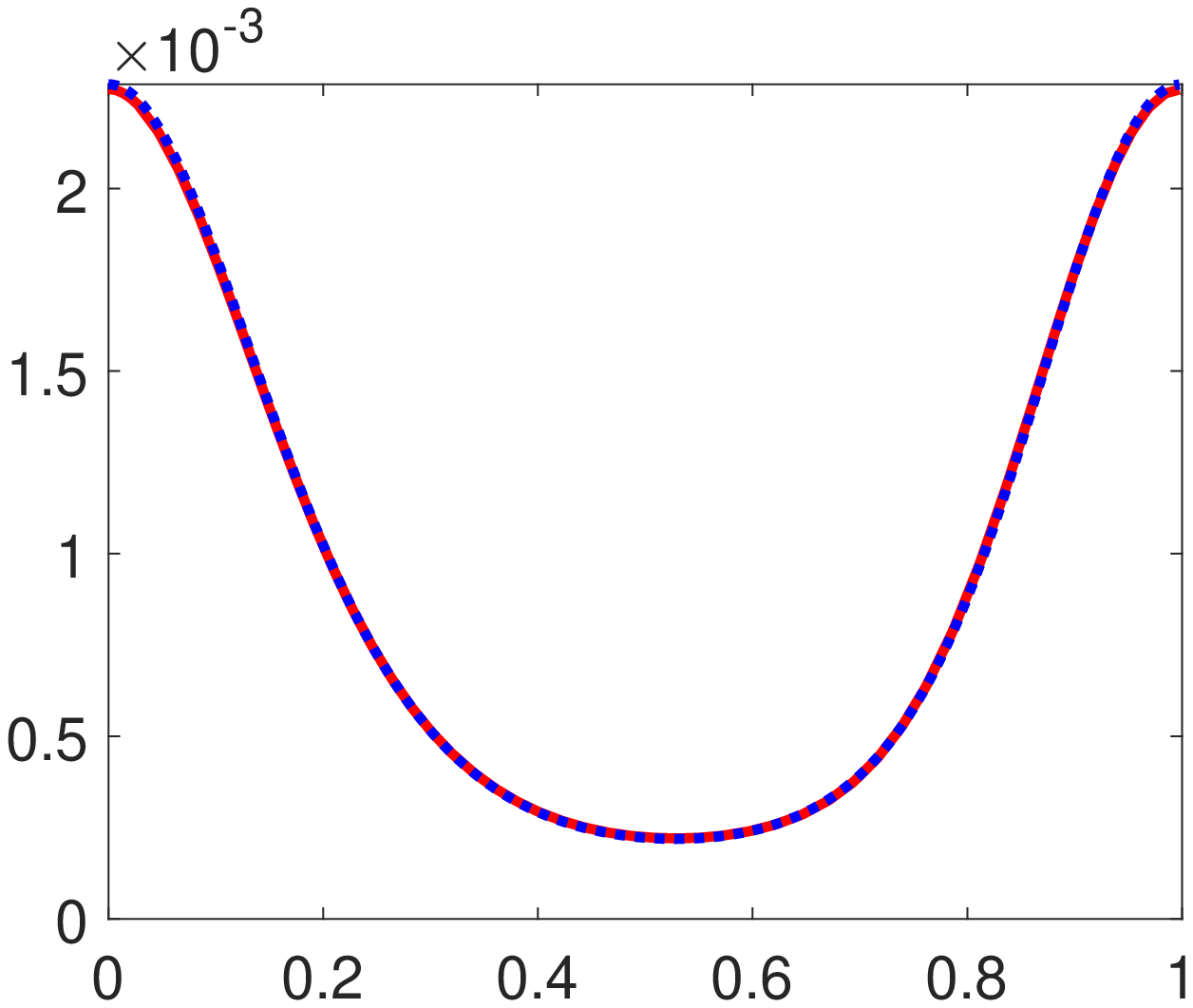}} &
\subfloat[$p_{1,6,3}$]{\includegraphics[width=0.21\textwidth]{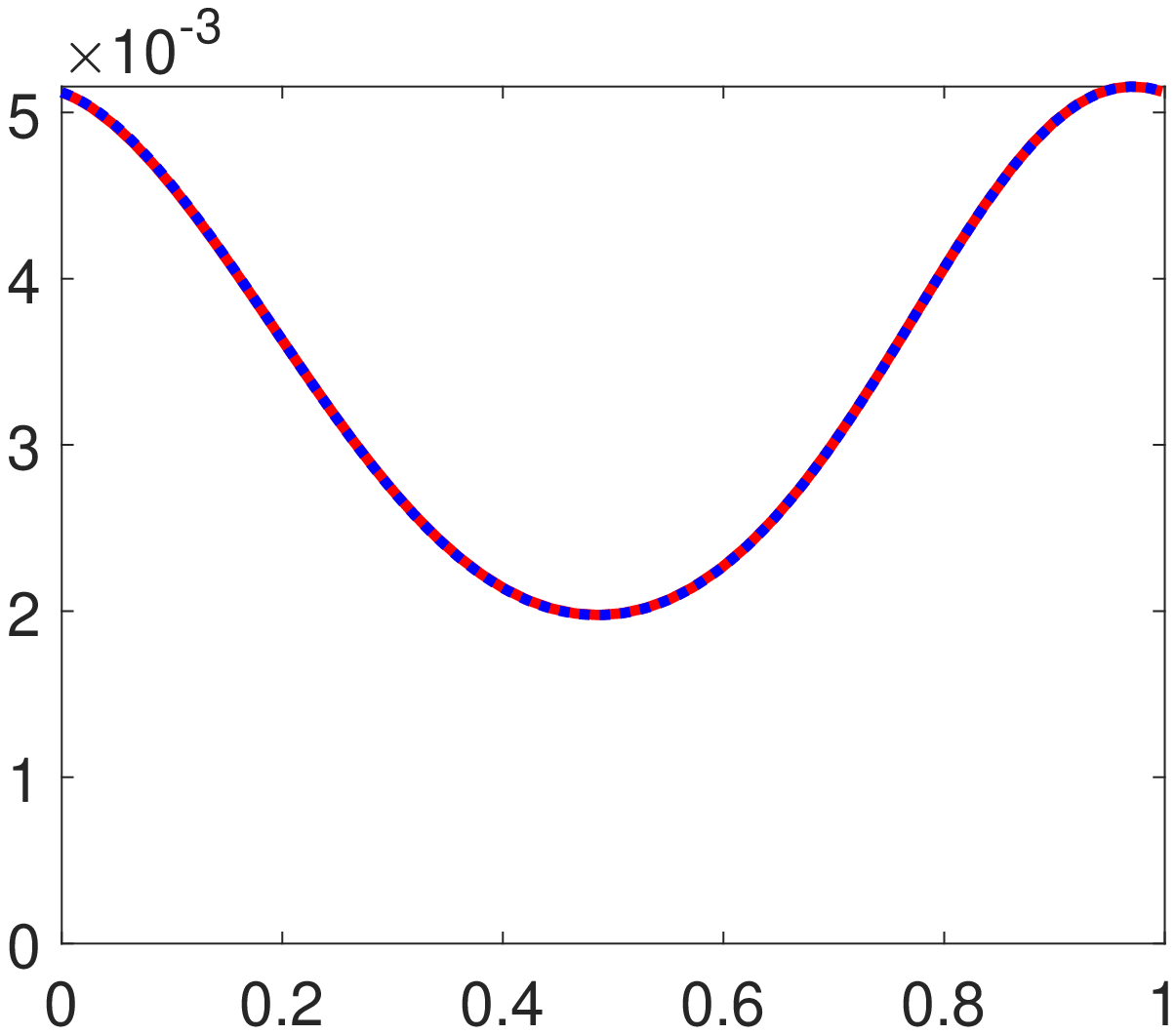}}\\
\multicolumn{2}{c}{\adjustbox{trim={.0\width} {.81\height} {0.0\width} {.0\height},clip}%
  {\includegraphics[width=0.45\textwidth]{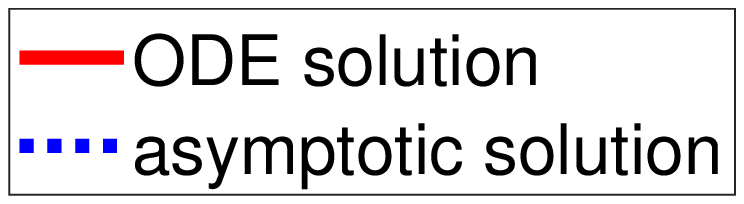}}}
&&
\end{tabular}
\caption{These graphs compare the asymptotic estimate given by ${\bf p}_1^{(1)}(t)$ for the probability of being in level 1 and the specified arrival and service phases for the $E_7/E_4/1$ system. See equation (\ref{eq:pjest}).}
\label{fig:level1q1}
\end{figure}

\begin{figure}[!tbp]
\begin{tabular}{ccc}
\subfloat[$p_{1,1,0}$,\;$q=1$]{\includegraphics[width=0.28\textwidth]{pi1_0_1qeq1}} &
\subfloat[$p_{1,1,1}$,\;$q=5$]{\includegraphics[width=0.28\textwidth]{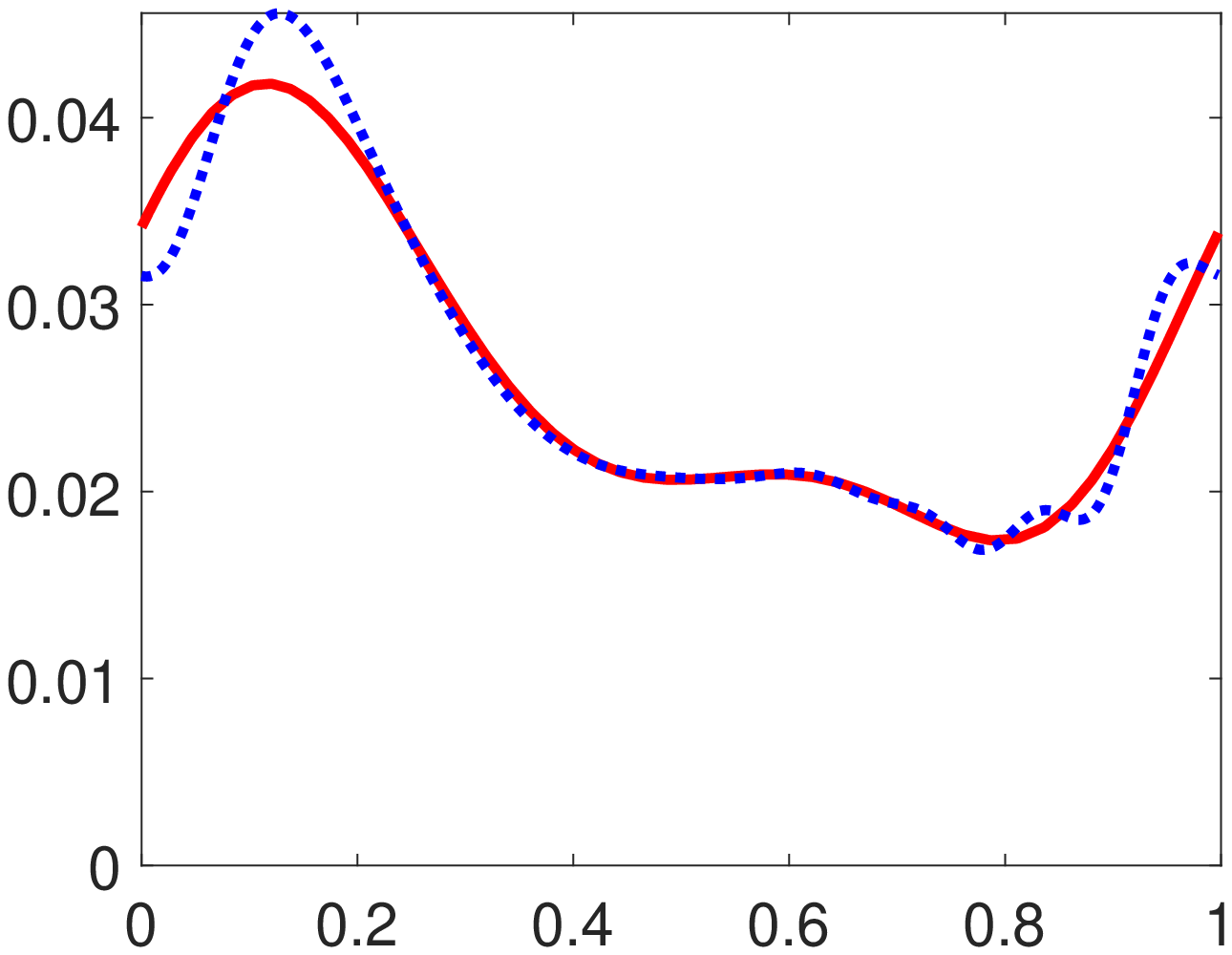}} &
\subfloat[$p_{1,1,2}$,\;$q=10$]{\includegraphics[width=0.28\textwidth]{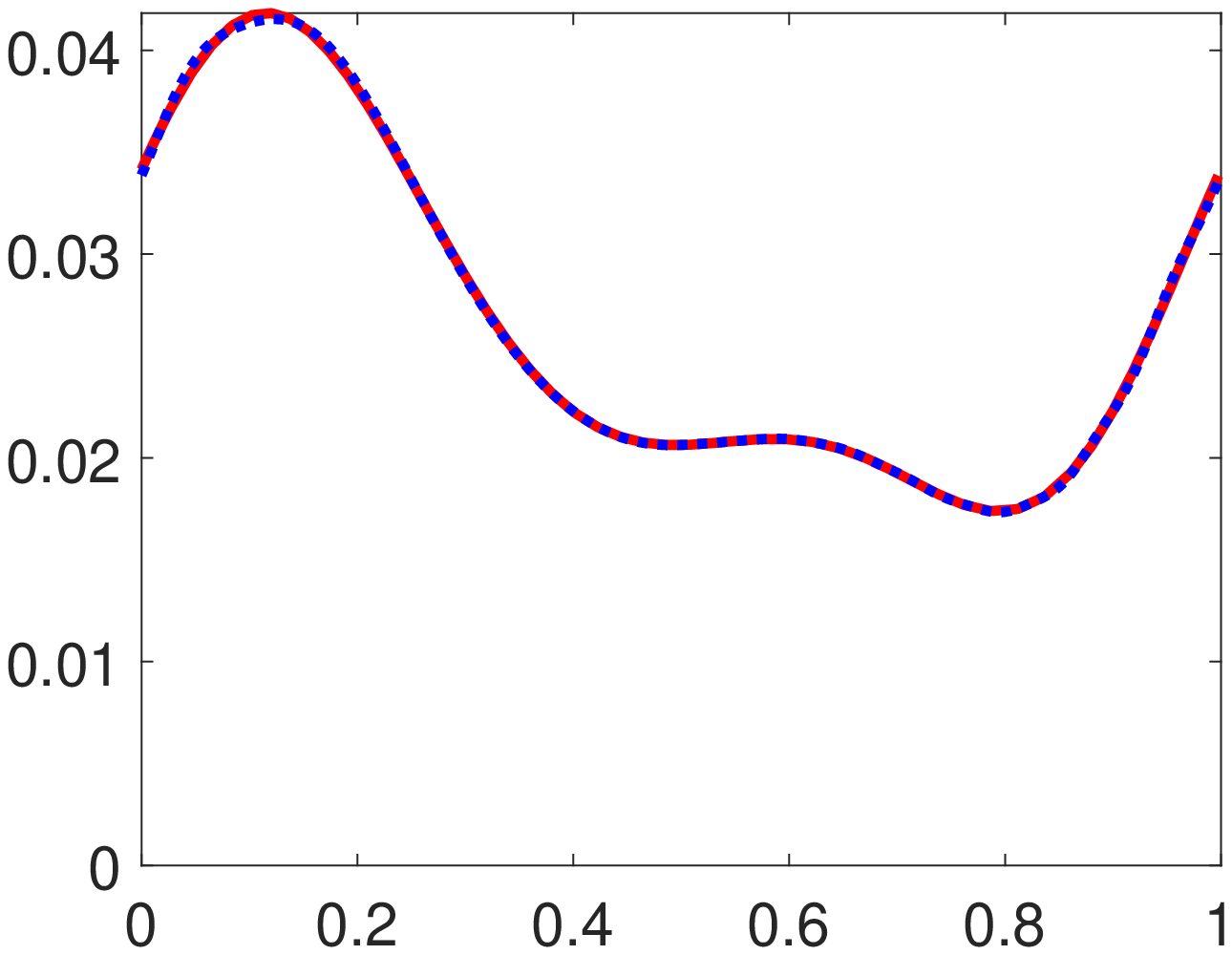}} \\
\multicolumn{2}{c}{\adjustbox{trim={.0\width} {.81\height} {0.0\width} {.0\height},clip}%
  {\includegraphics[width=0.45\textwidth]{legend}}}
&
\end{tabular}
\caption{These graphs compare the asymptotic estimate given by $p_{1,1,0}(t)$ for the probability of being in level 1, arrival phase 1 and service phase 0 for the $E_7/E_4/1$ system.  The asymptotic estimate is shown with the blue dashed line and the solution from a system of ordinary differential equations, truncated at 50 levels is shown in red.  The graphs are for three different values of $q$, with estimate given by $p_{1,1,0}^{(q)}(t)$ as defined in equation (\ref{eq:pjest}).}
\label{fig:level1q1-10}
\end{figure}

The plots in figures \ref{fig:level1q1} and \ref{fig:level1q1-10} show convergence of the asymptotic estimates \({\bf p}_j^{(q)}(t)\) to the level probabilities \({\bf p}_j(t)\) as the number of terms in the estimate increases.  In section \ref{sec:error} we provide error bounds for the asymptotic estimates of the level probabilities.

\section{Error Bound}\label{sec:error}

Our goal is to estimate the error
\begin{align*}
&\left\lVert{\bf p}_j(t)-{\bf p}^{(q)}_j(t)\right\rVert_\infty
\end{align*}
where ${\bf p}^{(q)}_j(t)$ is defined in equation (\ref{eq:pjest}).  Our first bound applies for $j\ge 3$.  We do this by finding bounds for
\begin{enumerate}[label=(\alph*)]
\item the modulus of the roots $|\chi_{\ell,n}|$,
\item $f(\chi_{\ell,n},t)$, defined in equation (\ref{eq:fxt}),  and on
\item $\left\lVert\left[\begin{array}{cccc}1&\chi_{\ell,n}^{-1/k}&\cdots&\chi_{\ell,n}^{(1-k)/k}\end{array}\right]\otimes\left[\begin{array}{cccc}1&\chi_{\ell,n}^{1/m}&\cdots&\chi_{\ell,n}^{(m-1)/m}\end{array}\right]\right\rVert_\infty$.
\end{enumerate}

Our asymptotic estimates for the \({\bf p}_j(t)\) are governed by the singularities of the generating function $P(z,t)$ and the function $f(x,t)$ given in equation (\ref{eq:fxt}).  The \(km\)th roots of these are the zeros of 
\begin{equation*}
1-{\rm e}^{\bar\lambda(y^m-1)+\bar\mu(y^{-k}-1)}.
\end{equation*}
We examine the asymptotic behavior of the roots which are outside the unit circle.    Return again to equation (\ref{eq:mypoly1}).  Write $z=r{\rm e}^{i\theta}$ in polar form and consider the limit of $\frac{z^{1/k}}{n}$ as $n\to\infty$.  Assume $r>1$ and that $r^{1/k}$ is also positive.  From equation (\ref{eq:mypoly1}), the roots of the singularities of the generating function satisfy
\begin{equation}\label{eq:useroot}
\bar\lambda z^{1/k}\omega_k^\ell=\bar\lambda-\bar\mu(\omega_m^jz^{-1/m}-1)+2\pi i n,\;\;n\in\mathbb{Z}.
\end{equation}
Dividing both sides by \(\bar\lambda\omega_k^\ell{\rm e}^{i\theta/k}n\),
\begin{multline*}
\lim_{n\to\pm\infty}\frac{r^{1/k}}{n}=\lim_{n\to\pm\infty}\left(\frac{\bar\lambda}{\bar\lambda n}-\frac{\bar\mu(\omega_m^j r^{-1/m}{\rm e}^{-i\theta/m}-1)}{\bar\lambda n}+\frac{2\pi i n}{\bar\lambda n}\right){\rm e}^{-i\theta/k}\omega_k^{-\ell}\\
=\lim_{n\pm\to\infty}\left(\frac{2\pi i n}{\bar\lambda n}\right){\rm e}^{-i\theta/k}\omega_k^{-\ell}
=\frac{2\pi }{\bar\lambda }
\end{multline*}
where the last equality follows from the fact that $r^{1/k}$ is real and positive.  This, in turn, implies that the limiting value of $\theta$, $\theta^*$, as $n\to\pm\infty$ is such that
\[{\rm e}^{-i\theta^*/k-2\pi i\ell/k+\pi i/2 }=1\]
or
\[\theta^*=\frac{k\pi}{2}.\]
Hence, for $r>1$,
\[r\sim\left(\frac{2\pi n}{\bar\lambda}\right)^k.\]
Similarly, if $r<1$, then
\[\lim_{n\to\pm\infty}\frac{r^{-1/m}}{n}=\frac{2\pi}{\bar\mu}\]
and the limiting value of $\theta$ is \(-\frac{m\pi}{2}\).

From the preceding analysis, we see that the modulus of the $k$th root of the singularity is bounded by
\begin{equation}\label{eq:rootbound}
\frac{2\pi |n|}{\bar\lambda}-\frac{\bar\lambda+2\bar\mu}{\sqrt{2}\bar\lambda}<\left\lvert \chi_{\ell,n}^{1/k}\right\rvert<\frac{2\pi |n|}{\bar\lambda}+\frac{\bar\lambda+2\bar\mu}{\sqrt{2}\bar\lambda}.
\end{equation}
These bounds are independent of $\ell=0,\dots,m-1$.

Next we consider the function $f(x,t)$ given in equation (\ref{eq:fxt}).  We consider three expressions separately:
\begin{equation*}
f_1(x,u,t)={\rm e}^{\int_u^t(\lambda(\nu)(x^{1/k}-1)+\mu(\nu)(x^{-1/m}-1))d\nu},
\end{equation*}
\begin{equation*}
f_2(x)={m\bar\lambda x^{1/k}-k\bar\mu x^{-1/m}}
\end{equation*}
and
\begin{equation*}
f_3(x,u)=\left(p_{0,k-1}(u)x\lambda(u)-\mu(u)\sum_{q=0}^{k-1}p_{1,q,m-1}(u)x^{q/k}\right)
\end{equation*}
so that
\begin{equation*}
f(x,t)=\int_{t-1}^t\frac{f_1(x,u,t)}
{f_2(x)}f_3(x,u)du.
\end{equation*}
We can compute the following bound for  $\lvert f_3(\chi,u)\rvert$, with \(\left\lvert \chi\right\rvert\ge 1\):
\begin{multline*}
\left\lvert
p_{0,k-1}(u)\chi\lambda(u)-\mu(u)\sum_{q=0}^{k-1}p_{1,q,m-1}(u)\chi^{q/k}\right\rvert\\
\le 
p_{0,k-1}(u)\left\lvert\chi\right\rvert\lambda(u)+\mu(u)\sum_{q=0}^{k-1}p_{1,q,m-1}(u)\left\lvert\chi^{q/k}\right\rvert\\
\le
\lambda(u)\left\lvert
\chi\right\rvert+\mu(u)\sum_{q=0}^{k-1}p_{1,q,m-1}(u)\left\lvert\chi^{q/k}\right\rvert\\
\le(\lambda(u)+\mu(u))\left\lvert\chi\right\rvert,
\end{multline*}
where we have used the fact that the phase transition rates are real and non-negative, as are probabilities.

We can find a lower bound for $\lvert f_2(\chi)\rvert$, $\lvert \chi\rvert>1$.  $\chi_{\ell,n}$ is a root of equation (\ref{eq:useroot}).  Because of the limit (\ref{eq:lim}), we may take the exponents $j$ and $\ell$ equal to zero, so
\begin{equation*}
\chi_{\ell,n}^{1/k}=\frac{1}{\bar\lambda}\left(\bar\lambda+\bar\mu-\bar\mu\chi_{\ell,n}^{-1/m}+2\pi i n\right).
\end{equation*}
Then
\begin{align*}
&\left\lvert f_2(\chi_{\ell,n})\right\rvert=\left\lvert m\bar\lambda \chi_{\ell,n}^{1/k}-k\bar\mu \chi_{\ell,n}^{-1/m}\right\rvert\\
&=\left\lvert m\left(\bar\lambda+\bar\mu-\bar\mu\chi_{\ell,n}^{-1/m}+2\pi i n\right)-k\bar\mu \chi_{\ell,n}^{-1/m}\right\rvert&&\text{substituting for $\chi_{\ell,n}^{1/k}$}\\
&=\left\lvert m\left(\bar\lambda+\bar\mu+2\pi i n\right)-(k+m)\bar\mu \chi_{\ell,n}^{-1/m}\right\rvert&&\text{collect terms}\\
&\ge
\left\lvert 
m\left\lvert \left(\bar\lambda+\bar\mu+2\pi i n\right)\right\rvert-(k+m)\bar\mu \left\lvert\chi_{\ell,n}^{-1/m}\right\rvert\right\rvert&&\parbox{.25\textwidth}{$k$, $m$ and $\bar\lambda$ are positive and $|a-b|\ge||a|-|b||$}\\
&\ge
\left\lvert 
m\left\lvert \left(\bar\lambda+\bar\mu+2\pi i n\right)\right\rvert-(k+m)\bar\mu \right\rvert&&\left\lvert\chi_{\ell,n}^{-1/m}\right\rvert<1\\
&=m\sqrt{(\bar\lambda+\bar\mu)^2+4\pi^2n^2}-(k+m)\bar\mu&&m\left\lvert \left(\bar\lambda+\bar\mu+2\pi i n\right)\right\rvert\ge(k+m)\bar\mu 
\end{align*}
for $\frac{m}{\bar\mu}>\frac{k}{\bar\lambda}$ (our ergodicity condition).

Now consider $\lvert f_1(\chi_{\ell,n},u,t)\rvert$.  
\begin{align*}
&\lvert f_1(\chi_{\ell,n},u,t)\rvert=\left\lvert{\rm e}^{\int_u^t(\lambda(\nu)(\chi_{\ell,n}^{1/k}-1)+\mu(\nu)(\chi_{\ell,n}^{-1/m}-1))d\nu}\right\rvert\\
&=\left\lvert{\rm e}^{\int_u^t\left(\lambda(\nu)\left(\frac{1}{\bar\lambda}\left[\bar\lambda+\bar\mu\left(1-\chi_{\ell,n}^{-1/m}\right)+2\pi i n\right]-1\right)+\mu(\nu)(\chi_{\ell,n}^{-1/m}-1)\right)d\nu}\right\rvert&&\text{substitution for $\chi_{\ell,n}^{1/k}$}\\
&=\left\lvert{\rm e}^{\int_u^t\left(\frac{\lambda(\nu)}{\bar\lambda}\left[\bar\mu+2\pi i n\right]+\left(\mu(\nu)-\frac{\lambda(\nu)\bar\mu}{\bar\lambda}\right)\chi_{\ell,n}^{-1/m}-\mu(\nu)\right)d\nu}\right\rvert&&\text{simplification}\\
&=\left\lvert{\rm e}^{\int_u^t\left(\frac{\lambda(\nu)\bar\mu}{\bar\lambda}+\left(\mu(\nu)-\frac{\lambda(\nu)\bar\mu}{\bar\lambda}\right)\chi_{\ell,n}^{-1/m}-\mu(\nu)\right)d\nu}\right\rvert&&\left\lvert{\rm e}^{\int_u^t\frac{\lambda(\nu)}{\bar\lambda} 2\pi i nd\nu}\right\rvert=1\\
&={\rm e}^{\int_u^t\frac{\lambda(\nu)\bar\mu}{\bar\lambda}-\mu(\nu)+\mu(\nu)\Re\{\chi_{\ell,n}^{-1/m}\} d\nu}\le {\rm e}^{\int_u^t\frac{\lambda(\nu)\bar\mu}{\bar\lambda}d\nu}.&&|\chi_{\ell,n}^{-1/m}|<1
\end{align*}
Then, putting these inequalities all together,
\begin{multline*}
\lvert f(\chi_{\ell,n},t)\rvert\le\\
\frac{|\chi_{\ell,n}|}{m\sqrt{(\bar\lambda+\bar\mu)^2+4\pi^2n^2}-(k+m)\bar\mu}\int_{t-1}^t(\lambda(u)+\mu(u)){\rm e}^{\frac{\bar\mu}{\bar\lambda}\int_u^t \lambda(\nu)d\nu}
du\\
=|\chi_{\ell,n}|C_n,
\end{multline*}
where
\begin{equation*}
C_n=\frac{\int_{t-1}^t(\lambda(u)+\mu(u)){\rm e}^{\frac{\bar\mu}{\bar\lambda}\int_u^t \lambda(\nu)d\nu}
du}{m\sqrt{(\bar\lambda+\bar\mu)^2+4\pi^2n^2}-(k+m)\bar\mu}.
\end{equation*}
The $C_n$ form a decreasing sequence. 

A bound on $\left\lVert\left[\begin{array}{cccc}1&\chi_{\ell,n}^{-1/k}&\cdots&\chi_{\ell,n}^{(1-k)/k}\end{array}\right]\otimes\left[\begin{array}{cccc}1&\chi_{\ell,n}^{1/m}&\cdots&\chi_{\ell,n}^{(m-1)/m}\end{array}\right]\right\rVert_\infty$ is $\lvert\chi_{\ell,n}\rvert$.

Applying the lower bound (because the exponent is negative) for  $\left|\chi_{\ell,n}^{1/k}\right|$ given in inequality (\ref{eq:rootbound}), we have
\begin{align*}
&\left\lVert{\bf p}_j(t)-{\bf p}^{(q)}_j(t)\right\rVert_\infty\\
&\le 2mC_q\sum_{n=q+1}^\infty \left(\frac{2\pi n}{\bar\lambda}-\frac{\bar\lambda+2\bar\mu}{\bar\lambda\sqrt{2}}\right)^{-kj+2k}&&\parbox{.3\textwidth}{The leading coefficient $2m$ is for $m$ roots for each fixed $n$, and two tails of the sum over $n$.  We also employ a bound on $|\chi_{\ell,n}|$ in this step.}\\
&\le 2mC_q\int_q^\infty \left(\frac{2\pi x}{\bar\lambda}-\frac{\bar\lambda+2\bar\mu}{\bar\lambda\sqrt{2}}\right)^{-kj+2k}dx&&\parbox{.3\textwidth}{For a monotone decreasing function, the given integral is greater than the sum.}\\
&=\frac{mC_q}{\pi}\frac{\left(2\pi q-\frac{1}{\sqrt{2}}(\bar\lambda+2\bar\mu)\right)^{k(2-j)+1}}{(k(j-2)-1)\bar\lambda^{k(2-j)}}&&\parbox{.3\textwidth}{This bound goes to zero as $q\to\infty$ for $j\ge 3$.}
\end{align*}

The plots in figures \ref{fig:level1q1} and \ref{fig:level1q1-10} show rapid convergence even for level one.  We explore why this is so in subsection \ref{subs:RLlemma}.

\subsection{A Riemann-Lebesgue type lemma}\label{subs:RLlemma}

\begin{figure}
\includegraphics[width=\textwidth]{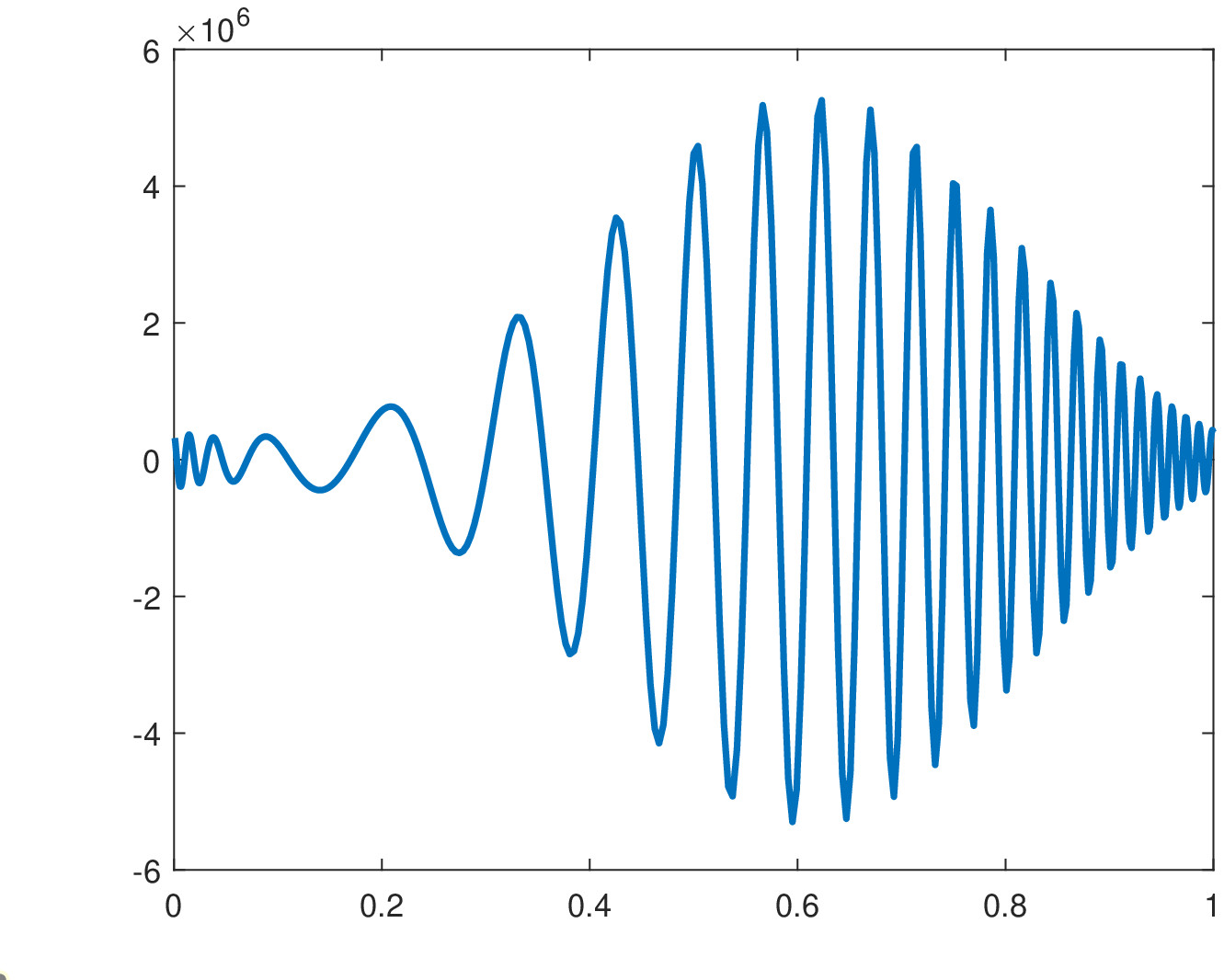}
\caption{Graph of ${\rm e}^{2\pi i n x}h_{\ell,25}(x)$ for $t=0.25$ for $E_7/E_4/1$ example}
\label{fig:osc}
\end{figure}

The functions $f(\chi_{\ell,n},t)$ defined in equation (\ref{eq:fxt}), for fixed $t$, are not Fourier coefficients, but they behave somewhat similarly.   We have
\begin{multline*}
f(\chi_{\ell,n},t)
=\int_{t-1}^t
\frac{{\rm e}
^{\int_u^t\left(\frac{\lambda(\nu)}{\bar\lambda}\left[\bar\mu+2\pi i n\right]+\left(\mu(\nu)-\frac{\lambda(\nu)\bar\mu}{\bar\lambda}\right)\chi_{\ell,n}^{-1/m}-\mu(\nu)\right)d\nu}}
{m(\bar\lambda+\bar\mu+2\pi i n)-(k+m)\bar\mu \chi_{\ell,n}^{-1/m}}
\\
\times
\left(p_{0,k-1}(u)\chi_{\ell,n}\lambda(u)-\mu(u)\sum_{q=0}^{k-1}p_{1,q,m-1}(u)\chi_{\ell,n}^{q/k}\right)
du.
\end{multline*}
We perform a change of variables.  For fixed $t$, let 
\begin{equation*}
x=\int_u^t\frac{\lambda(\nu)}{\bar\lambda}d\nu=g(u).
\end{equation*}
$g(u)$ is a decreasing function. Hence it has an inverse.  The differential 
\begin{equation*}
dx=-\frac{\lambda(u)}{\bar\lambda}du,
\end{equation*}
so
\begin{equation*}
du=-\frac{\bar\lambda}{\lambda(g^{-1}(x))}dx.
\end{equation*}
\begin{multline*}
f(\chi_{\ell,n},t)
=\int_{0}^1
\frac{{\rm e}
^{2\pi i n x+\bar\mu x-\bar\mu x\chi_{\ell,n}^{-1/m}+\int_{g^{-1}(x)}^t\mu(\nu)\left(\chi_{\ell,n}^{-1/m}-1\right)d\nu}}
{m(\bar\lambda+\bar\mu+2\pi i n)-(k+m)\bar\mu \chi_{\ell,n}^{-1/m}}
\\
\times
\left(p_{0,k-1}(g^{-1}(x))\chi_{\ell,n}\lambda(g^{-1}(x))-\mu(g^{-1}(x))\sum_{q=0}^{k-1}p_{1,q,m-1}(g^{-1}(x))\chi_{\ell,n}^{q/k}\right)\\
\times
\frac{\bar\lambda}{\lambda(g^{-1}(x))}dx
\end{multline*}
As $|n|\to\infty$,
\begin{equation*}
\chi_{\ell,n}^{-1/m}\approx\bar\lambda^{k/m}\left(\bar\lambda+\bar\mu+2\pi i n\right)^{-k/m}.
\end{equation*} 
Define 
\begin{multline*}
h_{\ell,n}(x)=
\frac{{\rm e}
^{\bar\mu x-\bar\mu x\chi_{\ell,n}^{-1/m}+\int_{g^{-1}(x)}^t\mu(\nu)\left(\chi_{\ell,n}^{-1/m}-1\right)d\nu}}
{m(\bar\lambda+\bar\mu+2\pi i n)-(k+m)\bar\mu \chi_{\ell,n}^{-1/m}}\\
\times
\left(p_{0,k-1}(g^{-1}(x))\chi_{\ell,n}\lambda(g^{-1}(x))-\mu(g^{-1}(x))\sum_{q=0}^{k-1}p_{1,q,m-1}(g^{-1}(x))\chi_{\ell,n}^{q/k}\right)\\
\times
\frac{\bar\lambda}{\lambda(g^{-1}(x))}
\end{multline*}
so that
\begin{equation*}
f(\chi_{\ell,n},t)=\int_0^1{\rm e}^{2\pi i n x}h_{\ell,n}(x)dx.
\end{equation*}
Let $\mathbb{T}=[0,1)$.
If $h_{\ell,n}(x)$ is a continuous $N$ times differentiable periodic function ($h_n(x)\in C^N(\mathbb{T})$ with $h_{\ell,n}^{(k)}(0)=h_{\ell,n}^{(k)}(1)$ for $0\le k\le N$), then repeated applications of integration by parts will yield
\begin{equation*}
f(\chi_{\ell,n},t)=\left(\frac{-1}{2\pi i n}\right)^N\int_0^1 h_{\ell,n}^{(N)}(x){\rm e}^{2\pi i n x}dx
\end{equation*}
where $h_{\ell,n}^{(N)}(x)$ is the $N$th derivative of $h_{\ell,n}(x)$.
Then
\begin{equation*}
\left\lvert f(\chi_{\ell,n},t)\right\rvert\le \left(\frac{1}{2\pi |n|}\right)^N \int_0^1\left|h_{\ell,n}^{(N)}(x)\right|dx.
\end{equation*}
Note that 
\begin{equation*}
\chi_{\ell,n}=\left(\frac{1}{\bar\lambda}\left[\bar\lambda+\bar\mu\left(1-\chi_{\ell,n}^{-1/m}\right)+2\pi i n\right]\right)^k
\end{equation*}
so that $h_{\ell,n}(x)\sim C(x) n^{k-1}$ for some function $C(x)$ that does not depend on $n$. The contribution from $\chi_{\ell,n}^{-1/m}\to 0$ as $|n|$ increases.  However, so long as $h_{\ell,n}(x)$ is sufficiently smooth, the integral 
\begin{equation*}
f(\chi_{\ell,n},t)=\int_0^1{\rm e}^{2\pi i n x}h_{\ell,n}(x)dx\to 0
\end{equation*}
as $n\to\infty$.  This happens because the rapid oscillations introduced by the factor ${\rm e}^{2\pi i n x}$ cause the integral to go to zero.  Figure \ref{fig:osc} shows a graph of ${\rm e}^{2\pi i n x}h_{\ell,25}(x)$ for $t=0.25$ to illustrate this idea.  See Loukas Grafakos text {\it Classical Fourier Analysis} \cite{Grafakos}, theorem 3.3.9, p. 196 for a similar result for Fourier coefficients.

% Grafakos p. 196 thm 3.3.9. 2014 edition

\section{Waiting time distribution}\label{sec:waiting}

If a customer enters the system when there are already $j$ customers ahead of him and the customer being served is in service phase $s$, then at least $m-s+m(j-1)$ additional service phases must be completed before he begins service and $m-s+mj$ must be completed before his service is finished and he leaves the queue.  Let $W_q(u)$ represent the waiting time until a customer arriving at time $u$ reaches the front of the queue  and $W(u)$ represent the waiting time including service for that customer.  The waiting time distributions, given that \(X(u)=j\), \(J(u)=s\) and \(j\ge 1\) is 
\begin{equation*}
P\{W_q(u)\le t|X(u)=j,J(u)=s\}=\sum_{q=mj-s}^\infty \frac{\left(\int_u^{u+t}\mu(\nu)d\nu\right)^q}{q!}{\rm e}^{-\int_u^{u+t}\mu(\nu)d\nu}
\end{equation*}
and
\begin{equation*}
P\{W(u)\le t|X(u)=j,J(u)=s\}=\sum_{q=mj-s+m}^\infty \frac{\left(\int_u^{u+t}\mu(\nu)d\nu\right)^q}{q!}{\rm e}^{-\int_u^{u+t}\mu(\nu)d\nu}.
\end{equation*}
From equation (\ref{eq:pjt}),
\begin{equation*}
P\{X(u)=j,J(u)=s\}=\sum_{n=-\infty}^\infty\sum_{\ell=1}^m f(\chi_{\ell,n},u)\left(\sum_{a=0}^{k-1}\chi_{\ell,n}^{-a/k}\right)\chi_{\ell,n}^{-j}\chi_{\ell,n}^{s/m}.
\end{equation*}
so the waiting time distribution for a customer entering at time \(u\) is given by
\begin{multline*}%\label{eq:Wq}
P\{W_q(u)\le t\}=\sum_{n=-\infty}^\infty \sum_{\ell=1}^m  f(\chi_{\ell,n},u)\left(\sum_{a=0}^{k-1}\chi_{\ell,n}^{-a/k}\right)\frac{\chi_{\ell,n}^{-1/m}}{1-\chi_{\ell,n}^{-1/m}}\times\\
\left(1-{\rm e}^{\int_u^{u+t}\mu(\nu)(\chi_{\ell,n}^{-1/m}-1)d\nu}\right)
\end{multline*}
and
\begin{multline*}
P\{W(u)\le t\}=
\sum_{n=-\infty}^\infty \sum_{\ell=1}^m  f(\chi_{\ell,n},u)\left(\sum_{a=0}^{k-1}\chi_{\ell,n}^{-a/k}\right)\frac{\chi_{\ell,n}^{-1/m}}{1-\chi_{\ell,n}^{-1/m}}\times\\
\left(1-{\rm e}^{-\int_u^{u+t}\mu(\nu)d\nu}\sum_{q=0}^m\frac{\left(\int_u^{u+t}\mu(\nu)d\nu\right)^q}{q!}\right.\\
\left.-\chi_{\ell,n}{\rm e}^{-\int_u^{u+t}\mu(\nu)d\nu}\left({\rm e}^{\chi_{\ell,n}^{-1/m}\int_u^{u+t}\mu(\nu)d\nu}
-\sum_{q=0}^m\frac{\left(\chi_{\ell,n}^{-1/m}\int_u^{u+t}\mu(\nu)d\nu\right)^q}{q!}
\right)
\right).
\end{multline*}

\begin{figure}[!tbp]
\begin{tabular}{cc}
\subfloat[$P\{W_q(.2)\le t\}$, $q=0$]{\includegraphics[width=0.42\textwidth]{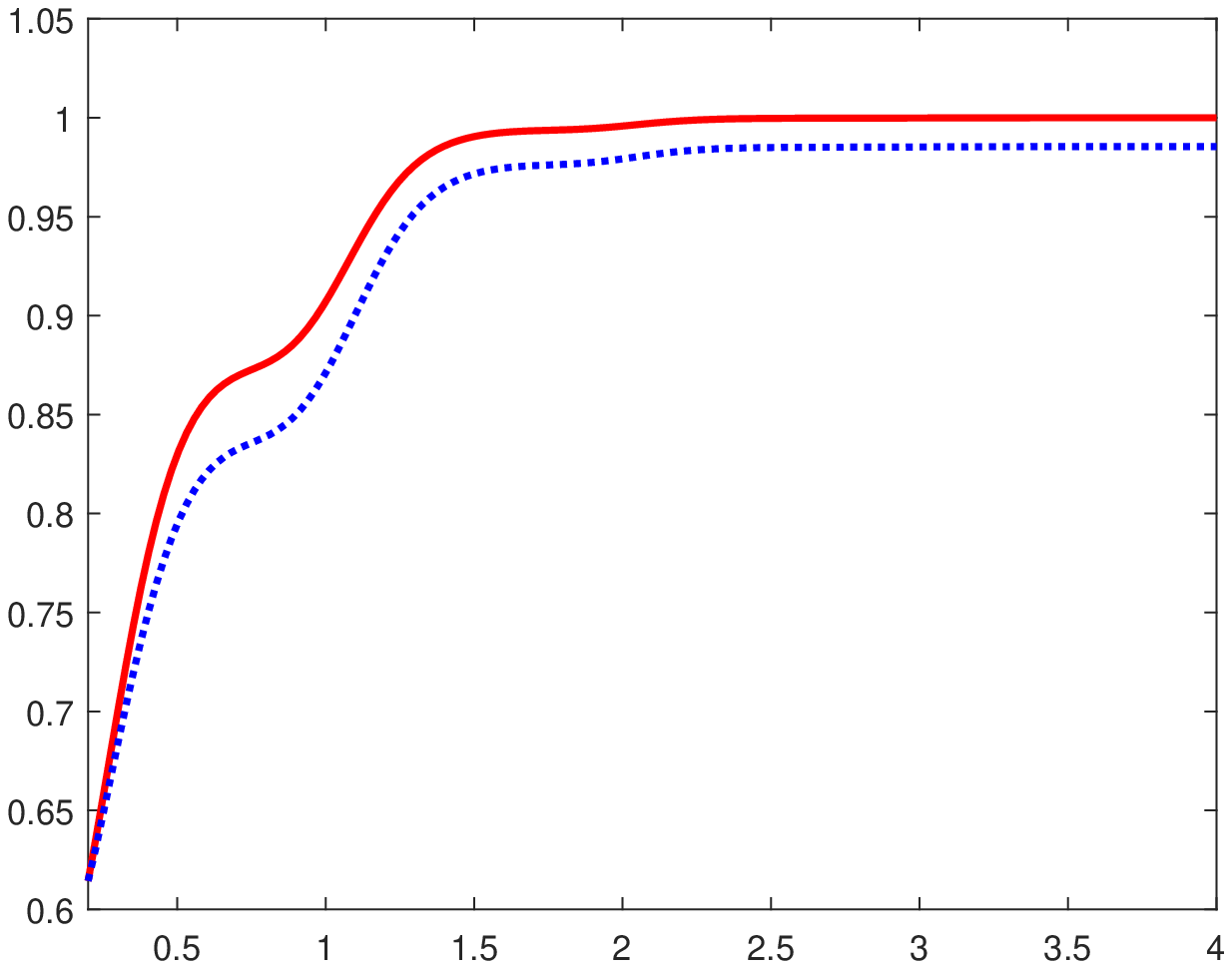}} &
\subfloat[$P\{W_q(.2)\le t\}$, $q=2$]{\includegraphics[width=0.42\textwidth]{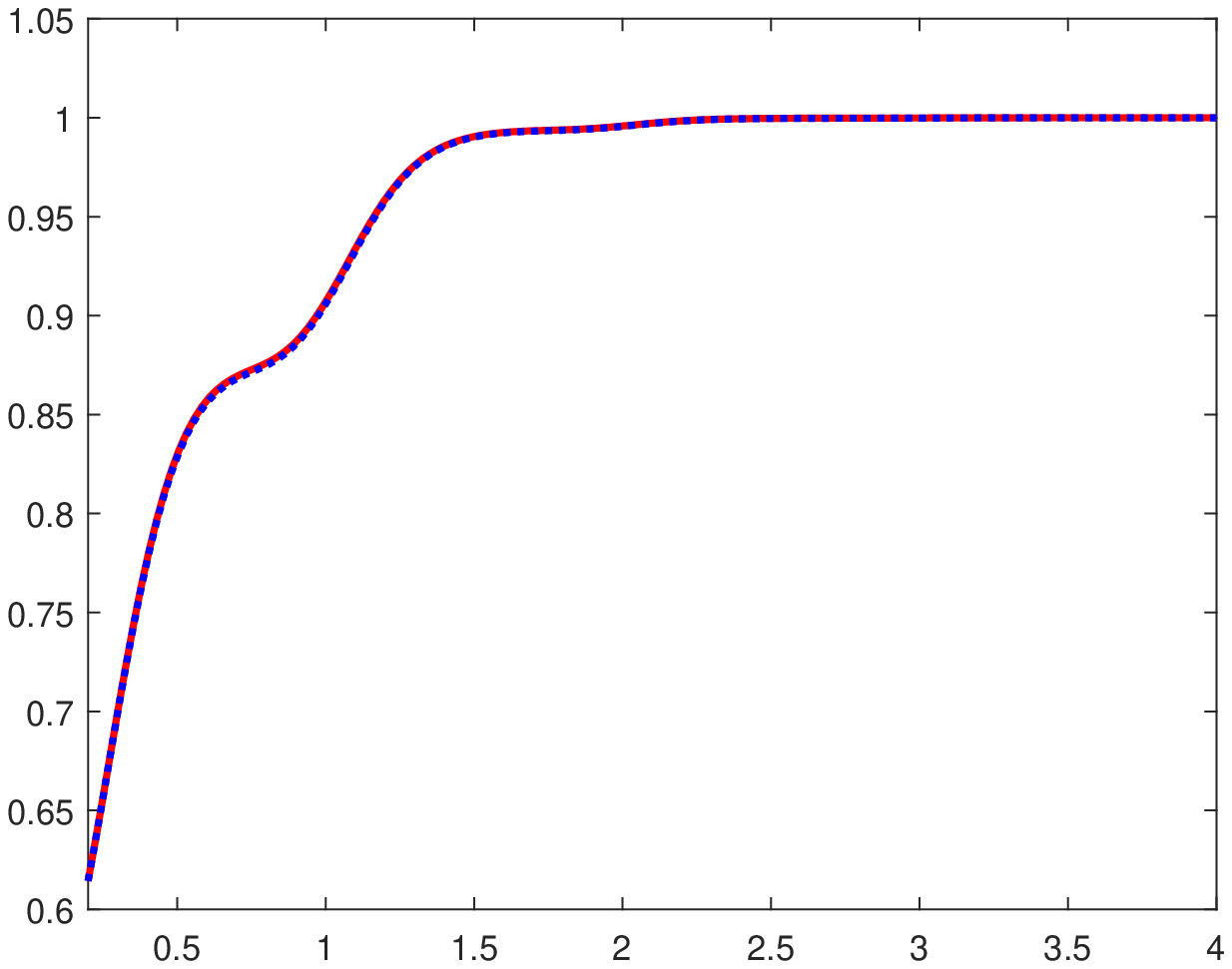}}\\
\subfloat[$P\{W_q(.7)\le t\}$, $q=0$]{\includegraphics[width=0.42\textwidth]{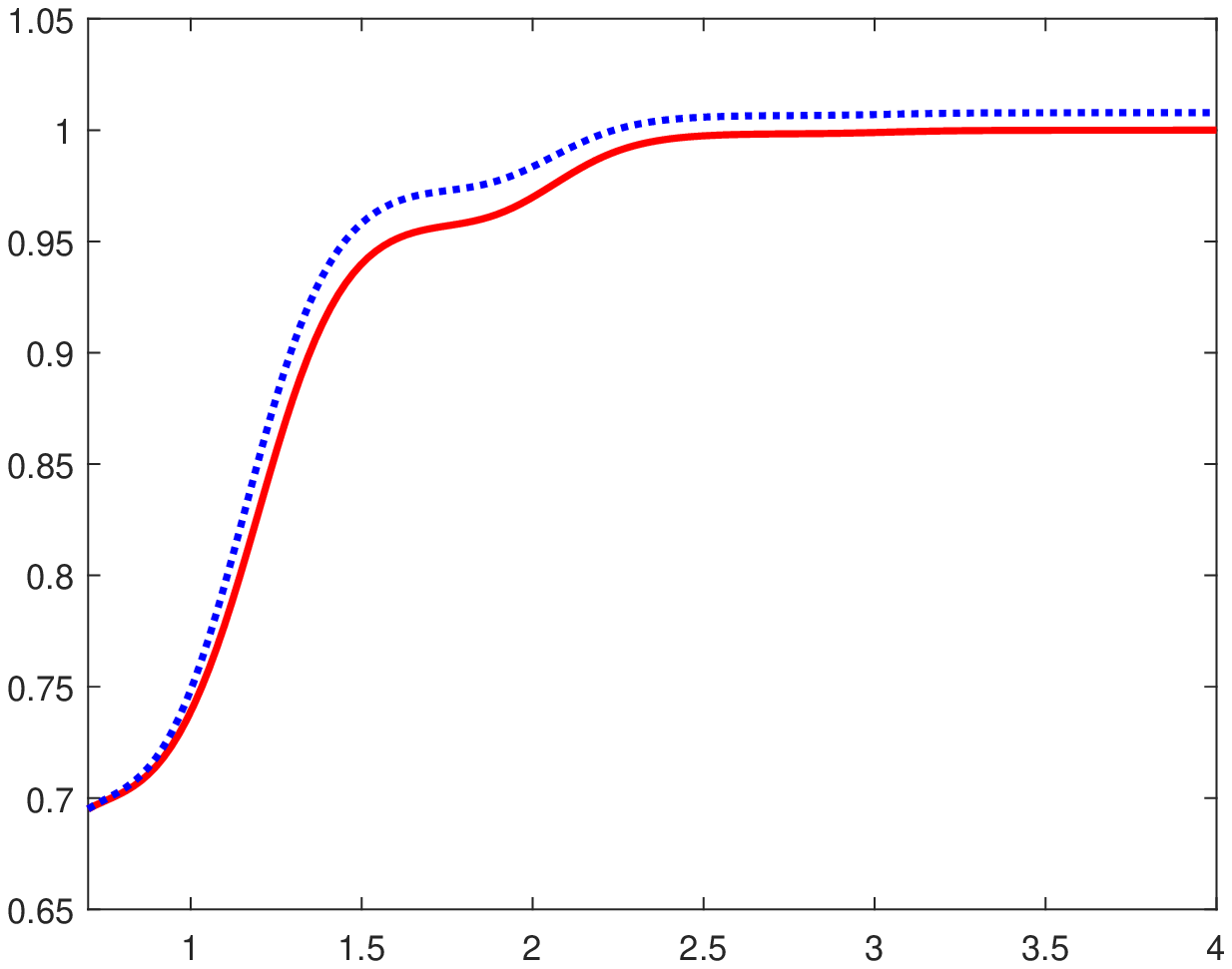}} &
\subfloat[$P\{W_q(.7)\le t\}$, $q=2$]{\includegraphics[width=0.42\textwidth]{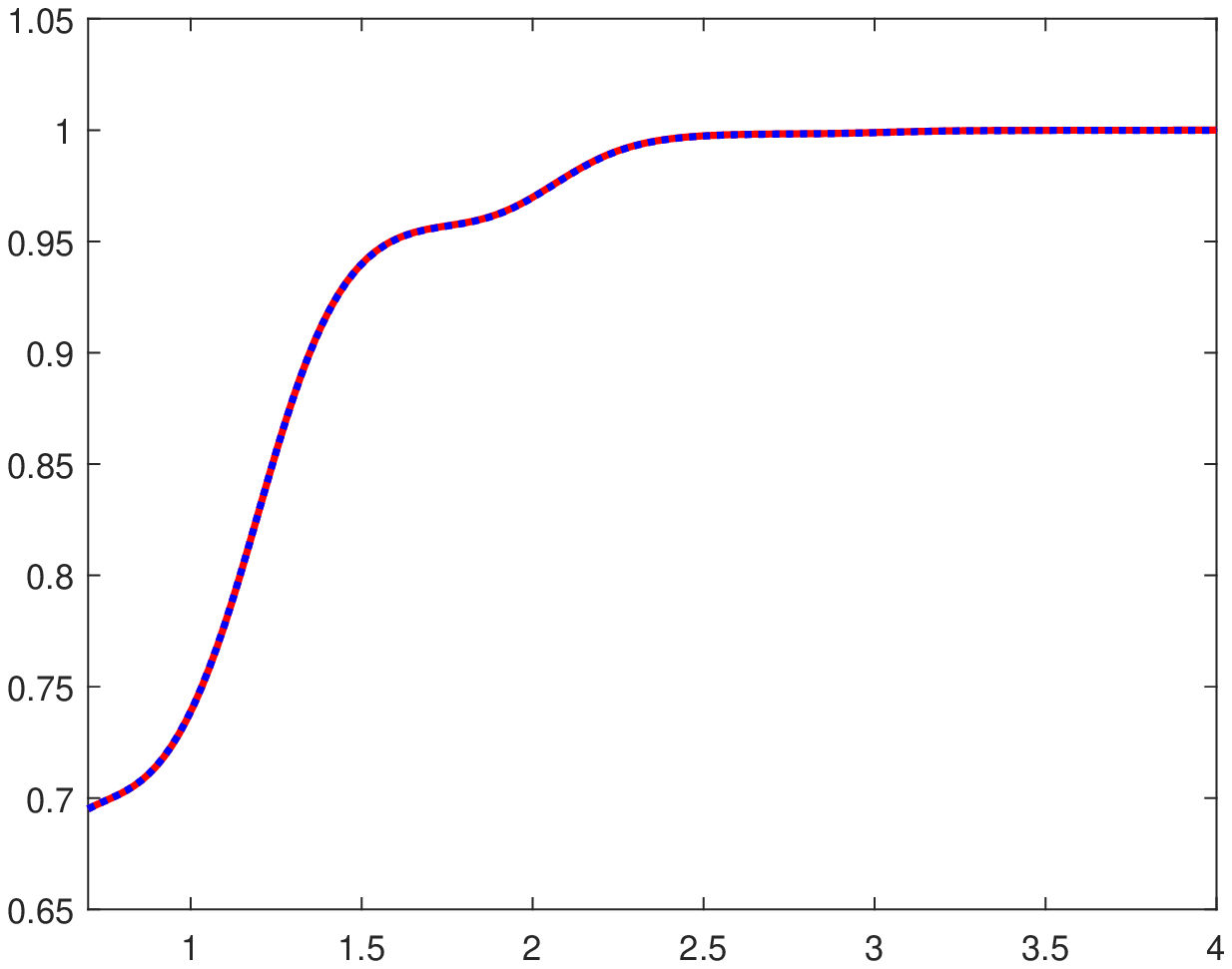}}\\
\multicolumn{2}{c}{\adjustbox{trim={.0\width} {.81\height} {0.0\width} {.0\height},clip}%
  {\includegraphics[width=0.45\textwidth]{legend}}}
\end{tabular}
\caption{These graphs compare the asymptotic estimate for the waiting time distribution to the ODE estimate for the distribution for the specified arrival and service phases for the $E_7/E_4/1$ system. See equation (\ref{eq:Wqest}).}
\label{fig:waitq0q2}
\end{figure}

We estimate the waiting time distribution with the expression
\begin{multline}\label{eq:Wqest}
P\{W_q^{(q)}(u)\le t\}=\sum_{n=-q}^q \sum_{\ell=1}^m  f(\chi_{\ell,n},u)\left(\sum_{a=0}^{k-1}\chi_{\ell,n}^{-a/k}\right)\frac{\chi_{\ell,n}^{-1/m}}{1-\chi_{\ell,n}^{-1/m}}\times\\
\left(1-{\rm e}^{\int_u^{u+t}\mu(\nu)(\chi_{\ell,n}^{-1/m}-1)d\nu}\right).
\end{multline}
An example appears in figure \ref{fig:waitq0q2}.

\section{Busy period distribution}\label{sec:busy}

In this section, we follow the approach of Baek, Moon and Lee \cite{BaekMoonLee} and apply it to the case of time-varying periodic parameters to find the busy period in terms of a Volterra equation of the second kind.
Let us define the first passage time 
\begin{equation*}
\tau_j=\inf\{t>u,N(t)=0|N(u)=j,J(u)=q\}.
\end{equation*}
We note that \(\tau_j\) is the length of a busy period that starts with \(j\) customers in the system and with an arriving customer in phase \(q\).

Let us define the following probabilities:
\begin{equation*}
{\bf Q}_n^{(j)}(t)=P\{N(t)=n,\tau_j>t|N(u)=j,J(u)=q\},
\end{equation*}

\begin{equation*}
{\bf Q}_0^{(j)}(t)=P\{\tau_j<t|N(u)=j,J(u)=q\}.
\end{equation*}
We find the busy time distribution.
\(\{N(t),J(t)\}\) is a continuous time Markov chain with absorbing boundary at \(N(\tau_j)=0\).  We set up the following system of ordinary differential equations:
\begin{eqnarray}\label{eq:busyODE}
\frac{d}{dt}{\bf Q}_0^{(j)}(t)&=&{\bf Q}_1^{(j)}(t){\bf A}_{-1}(t)\nonumber\\
\frac{d}{dt}{\bf Q}_1^{(j)}(t)&=&{\bf Q}_1^{(j)}(t){\bf A}_{0}(t)+{\bf Q}_2^{(j)}(t){\bf A}_{-1}(t)\\
\frac{d}{dt}{\bf Q}_n^{(j)}(t)&=&{\bf Q}_{n-1}^{(j)}(t){\bf A}_{1}(t)+{\bf Q}_{n}^{(j)}(t){\bf A}_{0}(t)+{\bf Q}_{n+1}^{(j)}(t){\bf A}_{-1}(t)\nonumber
\end{eqnarray}
To solve the system of differential equations (\ref{eq:busyODE}), we define the generating function
\begin{equation*}
G^{(j)}(z,u,t)=\sum_{n=0}^\infty {\bf Q}^{(j)}_n(t)z^n.
\end{equation*}
The differential equation for the generating function is
\begin{equation*}
\frac{d}{dt}G^{(j)}(z,u,t)=G^{(j)}(z,u,t){\bf A}_z(t)-\\
{\bf Q}^{(j)}_0(u,t){\bf A}_z(t)
\end{equation*}
with solution
\begin{equation*}
G^{(j)}(z,u,t)=G^{(j)}(z,u,u)\Phi(z,u,t)-\int_u^t{\bf Q}^{(j)}_0(\nu){\bf A}_z(\nu)\Phi(z,\nu,t)d\nu
\end{equation*}
where \(G(z,u,u)=z^j {\bf e}_q\) and \({\bf e}_q\) is a row vector with a one at component \(q\) and zeros elsewhere.  Since 
\({\bf Q}_0^{(j)}(t)=[z^0]G^{(j)}(z,u,t)\), 
\({\bf Q}_0^{(j)}(t)\) solves the Volterra equation of the second kind:
\begin{multline}\label{eq:busy}
{\bf Q}_0^{(j)}(t)={\bf e}_q[z^{-j}]\Phi(z,u,t)
-\int_u^t {\bf Q}^{(j)}_0(\nu)\left({\bf A}_{-1}(\nu)[z^1]\Phi(z,\nu,t)\right.\\
\left.+{\bf A}_0(\nu)[z^0]\Phi(z,\nu,t)+{\bf A}_1(\nu)[z^{-1}]\Phi(z,\nu,t)\right)d\nu.
\end{multline}
The matrix coefficient on $z^n$ in the generating function for the unbounded process that appears (several times) in equation (\ref{eq:busy}) is given in equation (\ref{eq:Phicoeff}).  For example,
\begin{multline*}
[z^{-j}]\left[\Phi(z,u,t)\right]_{(a_1,s_1)(a_2,s_2)}\\
={\rm e}^{-\int_u^t(\lambda(\nu)+\mu(\nu))d\nu}\sum_{\ell=j}^\infty\frac{\left(\int_u^t\mu(\nu)d\nu\right)^{\ell m+s}\left(\int_u^t\lambda(\nu)d\nu\right)^{(-j+\ell)k+a}}{(\ell m+s)!((-j+\ell)k+a)!}
\end{multline*}
with $a=a_2-a_1$, and $s=s_2-s_1$ and $a\ge 0$, $s\ge 0$.
 
\section{Conclusion}\label{sec:conclusion}
In this paper, we developed a method for computing the asymptotic periodic distribution of the level and phase probabilities for a queue with $k$ Erlang arrival phases and $m$ Erlang service phases.  We also showed how to compute the waiting time distribution seen by a customer arriving at any time within the period assuming that the system is in its asymptotic periodic ``steady-state".  This calculation requires computing an integral over a single time-period.  We provide exact Fourier like expansions, but require only finitely many of these terms to compute the level probabilities to arbitrary accuracy.  We compare our results to those obtained by solving a truncated version of the infinite system of differential equations and letting the system run until steady-state is achieved.

The computations require the asymptotic periodic solution for the queue being idle or having a single customer.  These probabilities can be computed using Tikhonov regularization. We also express the busy period as a solution of a Volterra equation of the second kind.

%\paragraph{Paragraph headings} Use paragraph headings as needed.

%\begin{acknowledgements}
%If you'd like to thank anyone, place your comments here
%and remove the percent signs.
%\end{acknowledgements}

\section*{Declarations}

\subsection*{Conflict of interest}
The authors declare that they have no conflicts of interest.
\subsection*{Funding}
Not applicable.
\subsection{Conflicts of interest/Competing interests}
Not applicable.
\subsection*{Code availability}
Not applicable.

% BibTeX users please use one of
%\bibliographystyle{spbasic}      % basic style, author-year citations
\bibliographystyle{spmpsci}      % mathematics and physical sciences
%\bibliographystyle{spphys}       % APS-like style for physics
%\bibliography{}   % name your BibTeX data base
%\bibliographystyle{amsplain} not on ok list
\bibliography{Erlangbib}

\end{document}